\documentstyle{amsppt}
\magnification\magstep2
\parindent=1.5em
\hoffset=-1truecm
\pagewidth{18truecm}
\pageheight{9.2truein}
\TagsOnRight
\CenteredTagsOnSplits
\def\r{\roman}
\topmatter
\title{The true order of the Riemann zeta--function on the
                 critical line}
\endtitle
\author
N.\ V.\ Kuznetsov
\endauthor
\address {2 Komarova Str., Apt. 69, 680000, Khabarovsk, Russia}
\endaddress
\email {kuznet\@eol.ru}
\endemail
\endtopmatter
\pageno=3

\document                                       
\rightheadtext{The true order of the riemann zeta--function}
\leftheadtext{N.\ V.\ Kuznetsov}
\tolerance=500

Аннотация
UDC 511.3+517.43+519.45

For the Riemann zeta-function on the critical line the terminal estimate
have been proved, which had been conjectured by Lindel\"of at the
beginning of this Centure. The proof is based on the authors relations
which connect the bilinear forms of the eigenvalues of the Hecke
operators with sums of the Kloosterman sums. By the way, it is proved
that for the Hecke series (which are associated with the eigenfunctions
of the automorphic Laplacian) the natural analogue of the Lindel\"of
conjecture is true also.
\noindent
Bibl. 14.


\head
\S\,0. PRELIMINARIES
\endhead

\noindent
\subhead
0.1.\  The main result
\endsubhead

One of the two main problems in the theory of the Riemann zeta-function is
the question: what is the order of this function on the critical line?

In this work I prove the following assertion: the Lindel\"of conjecture for
the Riemann zeta-function is true. Moreover the natural analogue of this
conjecture for the Hecke series is true also.

It means that we have the following two theorems.

\proclaim{Theorem 1} Let $\zeta(s)$ be the Riemann
zeta--function which is defined for $\r{Re}\, s>1$ by two equalities
$$ \zeta(s)=\sum_{n=1}^{\infty}\frac{1}{n^s}=
            \prod_p\left(1-\frac {1}{p^s}\right)^{-1}              \tag{0.1}
$$
where $p$ runs over all primes. Then for any $\varepsilon>0$
$$|\zeta(1/2+it)|\ll t^{\varepsilon}                       \tag{0.2} $$
as $t\to +\infty.$
\endproclaim
\proclaim{Theorem 2} Let $\Cal H_j(s)$ be the Hecke series which
corresponds to $j$--th eigenfunction of the automorphic Laplacian for
the case of the full modular group. Then for any fixed $j\geqslant 1$
and for any $\varepsilon>0$ we have
$$ |\Cal H_j(1/2+it)|\ll t^{\varepsilon}                   \tag{0.3} $$
as $t\to +\infty.$
\endproclaim

\head
\S\,1. INITIAL IDENTITIES
\endhead

To prove (0.2) and (0.3) I use the following known facts.
\noindent
\subhead
1.1.\  The fore--traces
\endsubhead

I restrict myself by the case of the full modular group $\Gamma$.

Let $\lambda_0=0<\lambda_1<\ldots\leqslant\lambda_j\leqslant\ldots\ $ are
the eigenvalues of the automorphic Laplacian
$\Cal L=-y^2\left(\frac{\partial^2}{\partial x^2}+
                  \frac{\partial^2}{\partial y^2}\right)$.
It means for $\lambda=\lambda_j$ there is non--zero solution of the
equation
$$    \Cal L u=\lambda u                                   \tag{1.1} $$
with the conditions
$u(\gamma z)\equiv u\left(\frac{az+b}{cz+d}\right)=u(z)$
for any $\gamma\in\Gamma$ and
$$(u,u)\equiv\int\limits_{\Gamma \backslash\Bbb H}|u(z)|^2d\mu(z)<\infty$$
(here $d\mu(z)=y^{-2}dx\,dy$ is the $\Gamma$--invariant measure on the upper
half plane $\Bbb H$; $\Gamma \backslash\Bbb H$ is the fundamental domain of the
full modular group $\Gamma$).

The continuous spectrum of this boundary problem lies on the half--axis
$\lambda\geqslant 1/4$; this spectrum is simple and the corresponding
eigenfunction is the analytical continuation of the Eisenstein series
$E(z,s)$. This series is defined by the equality
$$ E(z,s)=\sum_{\gamma\in\Gamma/\Gamma_{\infty}}
          \left(\r{Im}\,\gamma z\right)^s                          \tag{1.2}
$$
for $z\in\Bbb H$ and $\r{Re} s>1$; here $\Gamma_{\infty}$ is the syclic
subgroup which is generating by the transformation $z\mapsto z+1$.

For all $s$ we have the absolutely convergent series (the Fourier
expansion).

\proclaim{Theorem 1.1} (A.~Selberg, S.~Chowla [1]).
Let $z=x+iy$, $y>0$; then
$$ E(z,s)=y^s+\frac{\xi(1-s)}{\xi(s)}y^{1-s}+\frac 2{\xi(s)}
   \sum_{n\ne 0}\tau_s(n)\r{e}(nx)\sqrt y K_{s-1/2}(2\pi|n|y)        \tag{1.3}
$$
where $\r{e}(x)=\r{e}^{2\pi ix}$,
$$ \tau_s(n)=\sum_{d|n,d>0}\left(\frac{|n|}{d^2}\right)^{s-1/2}, \tag{1.4} $$
$$ \xi(s)=\pi^{-s}\Gamma(s)\zeta(2s)                             \tag{1.5} $$
and $K_{\nu}(\cdot)$ is the modified Bessel function (the Mcdonald
function, which is decreasing exponentially at $+\infty$) of the order $\nu$.
\endproclaim

Each eigenfunction $u_j$ of the discret spectrum has the similar Fourier
expansion, but without zeroth term:
$$ u_j(z)=\sqrt y\sum_{n\ne 0}\rho_j(n)\r{e}(nx)K_{i\varkappa_j}(2\pi|n|y)
                                                                  \tag{1.6} $$

Here $\rho_j(n)$ are the Fourier coefficients of $u_j$ and for
$j\geqslant 1$
$$ \varkappa_j=\sqrt{\lambda_j-1/4}.                               \tag{1.7} $$

Note that in the case of the full modular group
$\lambda_1\thickapprox 91.14$ (it is the result of the computer
calculations; see [14], p. 650-654).

We choose $u_j$ be real and each eigenfunction is even or odd under the
reflection operator
$\big(T_{-1}f\big)(z)=f(-\overline z)$; so we have
$$ T_{-1}u_j=\varepsilon_j u_j                                   \tag{1.8} $$
with $\varepsilon_j=+1$ or $\varepsilon_j=-1$.

Furthermore, it is possible take these eigenfunctions by such way that
they are the eigenfunctions for all the Hecke operators.

Let us define the $n$--th Hecke operator $T(n)$ by the equality
$$ \big(T(n) f\big)(z)=\frac{1}{\sqrt n}\sum\Sb ad=n\\ d>0\endSb
\,\,\,\,   \sum_{b(\bmod d)}f\left(\frac{az+b}d\right);          \tag{1.9} $$
then we have for all integers $n,m\geqslant 1$
$$ T(n)T(m)=\sum_{d|(n,m)}T\left(\frac{nm}{d^2}\right)=T(m)T(n).  \tag{1.10}$$

We take the common system of the eigenfunctions for $\Cal L$ and for all
$T(n)$, $n\geqslant 1$. In this case we have for all $n\geqslant 1$ and
$j\geqslant 1$
$$ \rho_j(n)=\rho_j(1)t_j(n),\ \ \
\rho_j(-n)=\varepsilon_j\rho_j(1)t_j(n),                      \tag{1.11}  $$
where $t_j(n)$ are such (real) numbers that
$$ T(n)u_j=t_j(n)u_j.                                          \tag{1.12}  $$

We have the following splendid expressions for the bilinear forms of
these eigenvalues $t_j(n)$.

\proclaim{Theorem 1.2} (N.~Kuznetsov [2], R.~Bruggeman [3]).
Let $h(r)$ be an even function in $r$ which is regular in the strip
$|\r{Im}\,r|\leqslant\Delta$ for some $\Delta>1/2$ and wich is
$\r{O}\left(|r|^{-2-\delta}\right)$ for some $\delta>0$ when
$r\to\infty$ inside of this strip. Then for any integers $n,m\geqslant 1$
$$\multline
\sum_{j\geqslant 1}\alpha_jt_j(n)t_j(m)h(\varkappa_j)+\frac 1{\pi}
\int\limits_{-\infty}^{\infty}\tau_{1/2+ir}(n)\tau_{1/2+ir}(m)
\frac{h(r)}{|\zeta(1+2ir)|^2}dr=                                        \\
=\frac 1{2\pi}\delta_{n,m}\int\limits_{-\infty}^{\infty}h(u)\,d\chi(u)+
\sum_{c\geqslant 1}\frac 1cS(n,m;c)\varphi
\left(\frac{4\pi\sqrt{nm}}c\right),
\endmultline                                                       \tag{1.13}
$$
where
$$ \alpha_j=(\cosh \pi\varkappa_j)^{-1}|\rho_j(1)|^2,                \tag{1.14} $$
$$ d\chi(u)=\frac 2{\pi}u\,\tanh(\pi u)\,du,                        \tag{1.15} $$
$ S $ denotes the Kloosterman sum,
$$ S(n,m;c)=\sum\Sb a(\bmod c)\\ ad\equiv 1(\bmod c)\endSb
          \r e\left(\frac{na+md}c\right),                          \tag{1.16} $$
and with notation ($J_{\nu}(\cdot)$-- the Bessel function of the
order $\nu$)
$$k_0(x,\nu)=\frac 1{2\cos(\pi\nu)}
\left(J_{2\nu-1}(x)-J_{1-2\nu}(x)\right)                         \tag{1.17} $$
the weight function in the sum of the Kloosterman sums is defined by the
integral transform
$$ \varphi(x)=\int\limits_{-\infty}^{\infty}
             k_0(x,1/2+ir)h(r)\,d\chi(r).                      \tag{1.18} $$
\endproclaim

The identity (1.13) is called "the Kuznetsov trace formula"; I think the
more preferable say "fore--trace". In the reality the famous Selberg
trace formula (for the full modular group and for congruence subgroups)
follows from (1.13). So the set of these identities with all
$n,m\geqslant 1$ may be considered as the set of primary equalities
to construct the trace formulae (and a lot of others identities).

\noindent
\subhead
1.2.\  The regular case
\endsubhead

The first example of the similar identities is the classical Peterson
formula.

\proclaim{Theorem 1.3} (H.~Peterson).
Let $f_{j,k}$ be the orthogonal Hecke basis in the space $\Cal M_k$ of
cusp forms of an even weight $k$, $\nu_k=\r{dim}\,\Cal M_k$ and
$t_{j,k}(n)$ are the eigenvalues of the Hecke operators
$T_k(n):\ \Cal M_k\to\Cal M_k$ under the normalization
$$ \big(T_k(n)f\big)(z)=n^{(k-1)/2}\sum\Sb d>0\\ ad=n\endSb
   \frac 1{d^k}\sum_{b(\bmod d)}f\left(\frac{az+b}d\right).     \tag{1.19} $$
Then we have
$$ \sum_{j=1}^{\nu_k}\|f_{j,k}\|^{-2}t_{j,k}(n)t_{j,k}(m)=
   \delta_{n,m}+2\pi i^{-k}\sum_{c\geqslant 1}
   \frac 1cS(n,m;c)J_{k-1}\left(\frac{4\pi\sqrt{nm}}c\right).   \tag{1.20} $$
\endproclaim

Note that for $k=2,4,6,8,10$ and $14$ the sum on the left side
(1.20) is zero, since $\nu_k=[k/12]$ if $k\not\equiv 2 (\bmod 12)$ and
$\nu_k=\left[\frac k{12}\right]-1$ for $k\equiv 2 (\bmod 12)$
(here $[x]$ denotes the integral part of $x$).

\noindent
\subhead
1.3.\ The sum of the Kloosterman sums
\endsubhead

We can invert (1.13) \big(and the similar identity with
$\varepsilon_j\alpha_j$ instead of $\alpha_j$ which is expessed in terms
of $S(n,-m;c)$\big) and we shall assume that the sum of Kloosterman sums
is given rather than the bilinear form in the Fourier coefficients.

\proclaim{Theorem 1.4} (N.~Kuznetsov [2],[4]).
Let $\varphi\in\r{C}^3(0,\infty),\ \ \varphi(0)=\varphi^{'}(0)=0$ and
assume that $\varphi(x)$ together with its derivatives up to third order
is $\r{O}(x^{-\beta})$ for some $\beta>2$ as $x\to+\infty$. Then, for any
integers $n,m\geqslant 1$, we have
$$ \multline \sum_{c\geqslant 1}\frac 1c
    S(n,m;c)\varphi\left(\frac{4\pi\sqrt{nm}}c\right)=
    \sum_{j\geqslant 1}\alpha_jt_j(n)t_j(m)h(\varkappa_j)+             \\
  +\frac 1{\pi}\int\limits_{-\infty}^{\infty}
   \tau_{1/2+ir}(n)\tau_{1/2+ir}(m)\frac{h(r)}{|\zeta(1+2ir)|^2}dr+    \\
  +\sum_{k\geqslant 6}g(k)\sum_{1\leqslant j\leqslant\nu_{2k}}
   \alpha_{j,2k}t_{j,2k}(n)t_{j,2k}(m),
   \endmultline                                                \tag{1.21} $$
where $h$ and $g$ are defined in terms of $\varphi$ by the
integral transforms
$$ h(r)=\pi\int\limits_{0}^{\infty}
               k_0(x,1/2+ir)\varphi(x)\frac{dx}x,             \tag{1.22} $$
$$ g(k)=(2k-1)\int\limits_{0}^{\infty}
               J_{2k-1}(x)\varphi(x)\frac{dx}x                \tag{1.23} $$
and $\alpha_{j,2k}$ is written instead of $\|f_{j,2k}\|^{-2}$.
\endproclaim

The proof may be found as well as in [4] and
(in more general situation) in [5],[11].

\proclaim{Theorem 1.5} (N.~Kuznetsov[4], M. Huxley [5]).
Assume that to a function $\psi:\ [0,\infty)\to\Bbb C$
the integral transform
$$ h(r)=2 \cosh (\pi r)\int\limits_{0}^{\infty}
               K_{2ir}(x)\psi(x)\frac{dx}x                    \tag{1.24} $$
associates the function $h$ satisfying to the conditions of
Theorem 1.2. Then, for this $\psi$ and for integers $n,m\geqslant 1$,
we have
$$ \multline \sum_{c\geqslant 1}\frac{1}{c}
    S(n,-m;c)\psi\left(\frac{4\pi\sqrt{nm}}c\right)=
    \sum_{j\geqslant 1}\varepsilon_j\alpha_jt_j(n)t_j(m)h(\varkappa_j)+ \\
  +\frac 1{\pi}\int\limits_{-\infty}^{\infty}
   \tau_{1/2+ir}(n)\tau_{1/2+ir}(m)\frac{h(r)}{|\zeta(1+2ir)|^2}dr.
   \endmultline                                                \tag{1.25} $$
\endproclaim

Note that (1.24) is true if $h$ is satisfying to conditions of
Theorem 1.2 and, for a given $h$, the taste function on the left side is
defined by the integral transform
$$ \psi(x)=\frac 4{\pi^2}\int\limits_{-\infty}^{\infty}
               K_{2ir}(x)h(r)\,r\,\sinh(\pi r)dr.                 \tag{1.26} $$

The pair of transforms (1.24) and (1.26), which are inverse to each
other, is the Kontorovich--Lebedev transform.

\noindent
\subhead
1.4.\ The Hecke series
\endsubhead

The Hecke series $\Cal H_{j}(s)$ and $\Cal H_{j,k}(s)$ are defined for
$\r{Re}\,s>1$ as
$$ \Cal H_j(s)=\sum_{n=1}^{\infty}n^{-s}t_j(n),\ \ \
   \Cal H_{j,k}(s)=\sum_{n=1}^{\infty}n^{-s}t_{j,k}(n).       \tag{1.27} $$

Here $t_j(n)$ and $t_{j,k}(n)$ are the eigenvalues of the Hecke
operators in the space of the cusp forms of the weight zero and the even
positive $k$ correspondingly.

The Hecke series which corresponds to continuous spectrum of the
automorphic Laplacian is
$$ \Cal L_{\nu}(s,x)=\sum_{n=1}^{\infty}n^{-s}\tau_{\nu}(n)\r{e}(nx),\ \ \
                            \r{e}(x)=\r{e}^{2\pi ix}.            \tag{1.28} $$

The following assertions are well known.

\proclaim{Theorem 1.6}
The Hecke series $\Cal H_j(s)$ and $\Cal H_{j,k}(s)$ are the entire
functions and these series have the functional equations of the Riemann
type:
$$ \Cal H_j(s)=(4\pi)^{2s-1}\gamma(1-s,1/2+i\varkappa_j)
   \big(-\cos(\pi s)+\varepsilon_j\cosh(\pi\varkappa_j)\big)
   \Cal H_j(1-s),                                                \tag{1.29} $$
$$ \Cal H_{j,k}(s)=-(4\pi)^{2s-1}\gamma(1-s,k/2)
   \cos(\pi s)\Cal H_{j,k}(1-s)                                  \tag{1.30} $$
where
$$ \gamma(u,v)=\frac 1{\pi}2^{2u-1}\Gamma(u+v-1/2)\Gamma(u-v+1/2).\tag{1.31} $$
\endproclaim

\proclaim{Theorem 1.7}
Let $x$ be rational, $x=\frac dc$ with $(d,c)=1,\ \ c\geqslant 1$.
Then for $\nu\ne 1/2$ the only singularities of $\Cal L_{\nu}(s,d/c)$
are simple poles at the point $s_1=\nu+1/2$ and $s_2=3/2-\nu$ with
residues $c^{-2\nu}\zeta(2\nu)$ and $c^{2\nu-2}\zeta(2-2\nu)$; the
function $\Cal L(s,d/c)$ has the functional equation
$$\multline \Cal L_{\nu}\left(s, \frac dc\right)=
    \left(\frac{4\pi}c\right)^{2s-1}\gamma(1-s,\nu)\times        \\
    \times\left\{-\cos(\pi s)\Cal L_{\nu}\left(1-s,-\frac ac\right)+
    \sin(\pi\nu)\Cal L_{\nu}\left(1-s,-\frac ac\right)\right\},
  \endmultline                                                   \tag{1.32} $$
where $a$ is defined by the congruence
$$  ad\equiv 1(\bmod c).                                         \tag{1.33} $$
\endproclaim

\noindent
\subhead
1.5.\ The Ramanujan identities
\endsubhead

Let
$$ \Cal Z(s;\nu,\mu)=\zeta(2s)\sum_{n=1}^{\infty}
     \frac{\tau_{\nu}(n)\tau_{\mu}(n)}{n^s};                         \tag{1.34}
$$
this series converges absolutely if
$\r{Re}\,s>1+|\r{Re}(\nu-1/2)|+|\r{Re}(\mu-1/2)|$.
The Ramanujan identity ([12], 1.3) gives the explicit expression for this
function.

\proclaim{Theorem 1.8}
$$ \Cal Z(s;\nu,\mu)=\zeta(s+\nu-\mu)  \zeta(s-\nu+\mu)
                     \zeta(s+\nu+\mu-1)\zeta(s-\nu-\mu+1).           \tag{1.35}
$$
\endproclaim

The analogues of this identity for the Hecke series are the consequence
of the multiplicative relations (1.10) for Hecke operators.

\proclaim{Theorem 1.9}
Let $\r{Re}\,s>1+|\r{Re}(\nu-1/2)|$; then
$$ \zeta(2s)\sum_{n=1}^{\infty}n^{-s}\tau_{\nu}(n)t_j(n)=
               \Cal H_j(s+\nu-1/2)\Cal H_j(s-\nu+1/2),               \tag{1.36}
$$
$$ \zeta(2s)\sum_{n=1}^{\infty}n^{-s}\tau_{\nu}(n)t_{j,k}(n)=
               \Cal H_{j,k}(s+\nu-1/2)\Cal H_{j,k}(s-\nu+1/2).       \tag{1.37}
$$
\endproclaim

\head
\S\,2. FUNCTIONAL EQUATION FOR THE FOURTH MOMENTS
\endhead

We define the fourth moments of the Hecke series by the equalities
$$\multline Z^{(d)}(s,\nu;\rho,\mu|h_0,h_1)=\sum_{j\geqslant 1}\alpha_j
   \big(h_0(\varkappa_j)+\varepsilon_jh_1(\varkappa_j)\big)\times   \\
   \times\Cal H_j(s+\nu-1/2)\Cal H_j(s-\nu+1/2)
   \Cal H_j(\rho+\mu-1/2)\Cal H_j(\rho-\mu+1/2),
  \endmultline                                                        \tag{2.1}
$$
$$\multline Z^{(r)}(s,\nu;\rho,\mu|g)=\\
=\sum_{k\geqslant 6}
   g(k)\sum_{1\leqslant j\leqslant \nu_{2k}}\alpha_{j,2k}
   \Cal H_{j,2k}(s+\nu-1/2)\Cal H_{j,2k}(s-\nu+1/2)\times          \\
   \times\Cal H_{j,2k}(\rho+\mu-1/2)\Cal H_{j,2k}(\rho-\mu+1/2).
  \endmultline                                                        \tag{2.2}
$$

It would be assumed here and further that $\r{Re}\,\nu=\r{Re}\,\mu=1/2$;
then the first series is the result of double summation
$$\multline Z^{(d)}(s,\nu;\rho,\mu|h_0,h_1)=\zeta(2s)\zeta(2\rho)
   \sum_{n,m\geqslant 1}
   \frac{\tau_{\nu}(n)\tau_{\mu}(m)}{n^sm^{\rho}}\times  \\
   \times\left(\sum_{j\geqslant 1}\alpha_j\left(h_0(\varkappa_j)+
   \varepsilon_jh_1(\varkappa_j)\right)t_j(n)t_j(m)\right),
  \endmultline                                                        \tag{2.3}
$$
if $\r{Re}\,s,\ \ \r{Re}\,\rho>1$; the function $Z^{(r)}(s,\nu;\rho,\mu| g)$
is defined by the similar way for these values of $s,\rho$.

The corresponding moment for the continuous spectrum is
$$ Z^{(c)}(s,\nu;\rho,\mu|h)=\frac 1{\pi}
   \int\limits_{-\infty}^{\infty}\Cal Z(s;\nu,1/2+ir)\Cal Z(\rho;\mu,1/2+ir)
   \frac{h(r)}{|\zeta(1/2+ir)|^2}dr                                   \tag{2.4}
$$
with the condition $\r{Re}\,s,\r{Re}\,\rho<1$
(for $\r{Re}\,\nu=\r{Re}\,\mu=1/2$); here $\Cal Z$ is defined by (1.35).

The same integral (2.4) with $\r{Re}\,s,\r{Re}\,\rho>1$
we denote by $Z^{(c,+)}(s,\nu;\rho,\mu|h)$.

The main purpose of this section is the correction of the proof of
Theorem 15 from [6].

\noindent
\subhead
2.1. The Mellin transform for the integral (1.18)
\endsubhead

Let for a given $h$ the function $\varphi$ is defined by the integral
transform (1.18).

First of all I give the explicit form of the Mellin transform for this
$\varphi$.

\proclaim{Proposition 2.1}
Let $h(r)$ be the even function which is regular in the strip
$|\r{Im}\;r|\leqslant\Delta,\ \ \Delta>1/2,$ and $|h(r)|$ decreases
faster than any fixed degree  of $r$ when $r\to\infty$ in this strip. If
$\varphi$ is defined by the integral (1.18),
then the Mellin transform of $\varphi$ is given by the equality
$$ \widehat\varphi(2w)\equiv\int\limits_0^{\infty}
   \varphi(x)x^{2w-1}dx=\frac 1{\pi}2^{2w-1}\cos(\pi w)
   \int\limits_{-\infty}^{\infty}
   \Gamma(w+iu)\Gamma(w-iu)h(u)\,d\chi(u).                         \tag{2.5}
$$
This function is regular in half--plane $\r{Re}\;w>-\Delta$ excepting
simple poles at $w=-1/2, -3/2, -5/2, \ldots$ and for any fixed
$\r{Re}\;w$ we have
$$ \widehat\varphi(2w)=2^{2w}\cos(\pi w)w^{2w-1}\r{e}^{-2w}
   \int\limits_{-\infty}^{\infty}
   \left(1+\frac{p_1(iu)}{w}+\frac{p_2(iu)}{w^2}+\ldots\right)
   h(u)\,d\chi(u)                                                     \tag{2.6}
$$
as $w\to\infty$; here $p_1, p_2,\ldots$ are polynomials of degree
$2,4,\ldots,$
$$
  p_1(z)=z^2+\frac{1}{6},\ p_2(z)=\frac {1}{2}z^4+\frac {2}{3}z^2+\frac
  {1}{72},
  \ p_3(z)=\frac {1}{6}z^6+\frac {5}{4}z^4+\frac{5}{24}z^2+\frac {1}{1620},
  \ldots                                                              \tag{2.7}
$$
\endproclaim
To prove these assertions we assume firstly that $\r{Re}\;w\in(0,1/4)$.
For this case we can integrate term by term and the known table
integrals give the equality (2.5) for $0<\r{Re}\;w<1/4$.

The regularity $\widehat\varphi(2w)$ for $\r{Re}\;w>0$ is obvious. In
the strip $-1<\r{Re}\;w<0$  we have
$$\multline
   \widehat\varphi(2w)=-\frac 1{\pi}2^{2w+2}w\sin(\pi w)
   \Gamma(2w)h(iw)+                                       \\
  +\frac{2^{2w}}{2\pi}\cos(\pi w)
   \int\limits_{-\infty}^{\infty}
   \Gamma(w+iu)\Gamma(w-iu)h(u)\,d\chi(u).
  \endmultline                                                        \tag{2.7}
$$

Let $|\r{Re}\,w|$ be small and $\r{Re}\,w>0$. Move the path of integration
in (2.5) to the new path $\Cal C $ made up of
\newline
\noindent
$\Cal C_{1}$: line segment $-\infty$ to $-|\r{Im}\,w| -\delta, \delta>
\r{Re}\,w$,
\newline
\noindent
$\Cal C_{2}$: anticlockwise semicircle centre $-|\r{Im}\,w|$ radius $\delta$,
such that the point $-iw$ is lying above the path,
\newline
\noindent
$\Cal C_{3}$: line segment $-|\r{Im}\,w|+\delta$ to $|\r{Im}\,w|-\delta$,
\newline
\noindent
$\Cal C_{4}$: clockwise semicircle centre $|\r{Im}\,w|$ radius $\delta$
(the point $iw$ is lying below the path),
\newline
\noindent
$\Cal C_{5}$: line segment $|\r{Im}\,w| +\delta$ to $+\infty$.

For $\r{Re}\,w >0$ we have
$$
\hat \phi (2w) = -2\pi i \r{Res}_{u=-iw} +2\pi i
\r{Res}_{u=iw}+\int\limits_{\Cal C}(...)du.
$$
Note that the sum of residues at the points $u=\pm iw$ is the first term
in (2.7).

But the integral over $\Cal C$ is regular for $-1< \r{Re}\,w<0$ and for
the case $\r{Re}\,w<0$ we can integrate over the real axis again; it
gives us the equality (2.7).

By the similar way we can receive the analytical continuation in the
strip $-2<\r{Re}\;w<-1$ (if $\Delta>1$) and so on; so
$\widehat\varphi(2w)$ is regular in the half--plane $\r{Re}\;w>-\Delta$
excepting those poles of $\Gamma(2w)$ which are not compensated by
zeroes of $\sin \pi w$.

After that we use the Barnes asymptotic expansion([7], 1.18). Namely, let
$B_n(z)$ are Bernoulli polynomials $(B_0=1,\ \ B_1=z-1/2,\ \
B_2=z^2-z+1/6,\ \ B_3=z^3-\frac 32z^2+\frac 12z,\ldots).$ The Barnes
expansion is the asymptotic series
$$
\log \,\Gamma(w+z)=(w+z-1/2)\log w- w +\frac {1}{2} \log(2\pi)+Q(w,z) \tag{2.8}
$$
with
$$
   Q(w,z)\!=\!\frac{B_2(z)}{1\cdot 2\;w}\!-\!
              \frac{B_3(z)}{2\cdot 3\;w^2}\!+\!
          \ldots\!+\!\frac{(-1)^{n+1}B_{n+1}(z)}{n(n+1)w^n}\!
          +\!\ldots;                                                  \tag{2.9}
$$
this expansion holds if $w\to\infty$, $|\r{arg}\,w|<\pi$ and
$z=o\left(|w|^{1/2}\right)$.

For any fixed $M$ and for $u_0=|w|^{\varepsilon}$
with the arbitrary small (but fixed) $\varepsilon >0$ our integral (2.5)
equals to
$$
\int\limits_{|u|\leqslant u_0}(...)\,d\chi(u) +O(|w|^{-M}),
$$
since we assume the very fast decreasing of $h$.

For $|u|\leqslant u_0$ we can use the Barnes expansion; it gives the
asymptotic equality
$$\multline
     \Gamma(w\!+\!z)\Gamma(w\!-\!z)\!=\!2\pi w^{2w-1}\r{e}^{-2w}
     \exp\Big(Q(w,z)\!+\!Q(w,-z)\Big)\!=                          \\
   =\!2\pi w^{2w-1}\r{e}^{-2w}\!\cdot\!
     \left(\!1\!+\!Q(w,z)\!+\!Q(w,-z)\!+\!\frac 1{2!}\!
     \Big(\!Q(w,z)\!+\!Q(w,-z)\!\Big)^2\!+\!\ldots\!\right).
  \endmultline
$$

Taking the finite number of terms from this expansion we write integral
in the limits ($-\infty,+\infty$) again and come to (2.6).

\subhead
2.2\ The regularization of the initial identity
\endsubhead

It follows from (2.6) that in the general case the Mellin transform
(2.5) of the integral (1.18) is $O(|w|^{2\r{Re}\,w -1})$ as $w\to
\infty$ with the fixed value of  $\r{Re}\, w$, excepting some cases when
$h(u)$ is orthogonal to the degrees $u^m, 0\leqslant m \leqslant m_0$
on the mesure $d\chi(u)$.

Nevertherless, there is the method of "regularization" of this
$\varphi$. We substract from the integral (1.18) the combination of the
Bessel functions; the coefficients of this combination may be determined
so that the difference has the Mellin transform which decreases more
rapidly.

This regularization is very essential for our proof and it will be used
later many times.

\proclaim{Proposition 2.2}
Let $L=\{l_1,\ldots,l_N\}$ be finite set of integers,
$1\leqslant l_1<l_2<\ldots<l_N, \ \ N\geqslant 1$.
Let the function $h$ be taken under assumptions of Proposition 2.1 with
the additional conditions $h(i(l-1/2)) =0$ for $1\leqslant l \leqslant
\Delta -1/2$.

We define $N$ coefficients $c(l)$ from the linear system
$$ \sum_{l\in L}(l-1/2)^{2m}(-1)^lc(l)=
         (-1)^m\int\limits_{-\infty}^{\infty}
         u^{2m}h(u)\,d\chi(u),\ \ \ \ 0\leqslant m\leqslant N-1.       \tag{2.10}
$$
Let
$$ \Phi_N(x)=\varphi(x)-\sum_{l\in L}c(l)J_{2l-1}(x),                \tag{2.11}
$$
where $\varphi$ is defined by the integral transform (1.18) of this $h$;
then the Mellin transform of this difference,
$$ \widehat{\Phi}_N(2w)=\int\limits_0^{\infty}\Phi_N(x)
                             x^{2w-1}\;dx,                           \tag{2.12}
$$
is the regular function in the half plane
$\r{Re}\,w>\max(-\Delta,-l_1+1/2)$ and for any fixed $\r{Re}\,w$ we have
$$
   \left|\widehat{\Phi}_N(2w)\right|\ll|w|^{2\r{Re}\,w-N-1}        \tag{2.13}
$$
as $w\to\infty$.
\endproclaim

Really,
$$\multline
   \int\limits_0^{\infty} J_{2l-1}(x)x^{2w-1}dx=
   2^{2w-1}\frac{\Gamma(l-1/2+w)}{\Gamma(l+1/2-w)}=      \\
 =  \frac 1{\pi}2^{2w-1}(-1)^l\cos(\pi w)
   \Gamma(w+l-1/2)\Gamma(w-l+1/2).
  \endmultline
$$
Now we have for large $|w|$ (and fixed $\r{Re}\,w$)
$$\multline \widehat{\Phi}_N(2w)=2^{2w}\cos(\pi w)
     w^{2w-1}\r{e}^{-2w}\times                            \\
     \times
     \Bigg\{\int\limits_{-\infty}^{\infty}
     \left(1+\frac{p_1(iu)}w+\frac{p_2(iu)}{w^2}+
     \ldots\right)h(u)\,d\chi(u)-                         \\
     -\sum_{l\in L}(-1)^lc(l)\left(1+\frac{p_1(l-1/2)}w+
     \frac{p_2(l-1/2)}{w^2}+\ldots\right)\Bigg\}
  \endmultline                                                       \tag{2.14}
$$
(note that two line segments $|u|\geqslant |w|^{\varepsilon}$ for any
fixed $\varepsilon >0$ give $O(|w|^{-M})$ for any $M\geqslant 6$; so the
Barnes expansion is used for $|u|\leqslant |w|^{varepsilon}$ only).

If $c(l)$ are defined by (2.10) then the terms $w^{-k}$ with
$0\leqslant k\leqslant N-1$ are cancelled. It gives the estimate (2.13).

Now we rewrite the initial identity (1.13).

\proclaim{Proposition 2.3}
Let $h(r)$ be an even function in $r$, regular in the strip
$|\r{Im}\,r|\leqslant \Delta$ with some $\Delta>3/2$ and let
$|h(r)|$ decreases faster than any fixed degree $|r|$ as $r\to\infty$
in this strip; in addition we assume $h\big((l-1/2)i\big)=0$ for
$1\leqslant l\leqslant\Delta+1/2$. Let $\Phi_N$ be defined by (2.11)
with coefficients from (2.10); then for any integers $n,m\geqslant 1$
we have
$$\multline \sum_{j\geqslant 1}\alpha_j t_j(n)t_j(m)h(\varkappa_j)+
    \frac 1{\pi}\int\limits_{-\infty}^{\infty}
    \tau_{1/2+ir}(n)\tau_{1/2+ir}(m)\frac{h(r)}{|\zeta(1+2ir)|^2}dr=    \\
    =\sum_{c\geqslant 1}\frac 1cS(n,m;c)
    \Phi_N\left(\frac{4\pi\sqrt{mn}}c\right)+\\
+ \frac 1{2\pi}\sum_{l\in L}(-1)^lc(l)
    \sum_{1\leqslant j\leqslant\nu_{2l}}\alpha_{j,2l}
    t_{j,2l}(n)t_{j,2l}(m).
  \endmultline                                                       \tag{2.15}
  $$
\endproclaim
This equality is the immediate consequence of the Peterson identity
(1.20). Writing
$$
   \varphi=\Phi_N+\sum_{l\in L}c(l)J_{2l-1},
$$
we get the sum of the Kloosterman sums with the test function $\Phi_N$ and
the finite number sums where the test function is the Bessel function of
odd order. The singular term with $\delta_{n,m}$ contributes
$$
   -\frac 1{2\pi}\delta_{n,m}\sum_{l\in L}(-1)^lc(l)=
   -\frac 1{2\pi}\delta_{n,m}\int\limits_{-\infty}^{\infty}
    h(u)\,d\chi(u).                                                  \tag{2.16}
$$
(the equation (2.10) with $m=0$). So the term with $\delta_{n,m}$
is disappearing and we come to (2.15).

\noindent
\subhead
2.2\ The functional equation
\endsubhead
\proclaim{Theorem 2.2}
Let two parameters $\nu,\mu$ are taken with conditions
$\r{Re}\,\nu=\r{Re}\,\mu=1/2$, $\nu\ne 1/2$, $\mu\ne 1/2$ and two
variables $s,\rho$ are taken inside of the narrow strip
$5/4<\r{Re}\,s,\r{Re}\,\rho<5/4+\varepsilon$ with small (but fixed)
$\varepsilon>0$.

Let $\Phi_N$ for a given $h$, which is satisfying to all conditions of
Proposition 2.3, be defined by (2.11) with $N\geqslant 1$ and
$l_1\geqslant 2$ Then we have
$$\multline
      Z^{(d)}(s,\nu;\rho,\mu|h,0)+Z^{(c,+)}(s,\nu;\rho,\mu|h,0)=    \\
   =\!Z^{(d)}(\rho,\nu;s,\mu|h_0,h_1)\!+\!
      Z^{(c,+)}(\rho,\nu;s,\mu|h_0,h_1)\!+\!
      Z^{(r)}(\rho,\nu;s,\mu|g)+                                    \\
   +\!\frac 1{2\pi}\sum_{l\in L}(-1)^lc(l)z_{2l}(s,\nu;\rho,\mu)\!+\!
      \Cal R_h(s,\nu;\rho,\mu)\!+\!\Cal R_h(s,1-\nu;\rho,\mu)+       \\
   +\!\Cal R_h(\rho,\mu;s,\nu)\!+\!\Cal R_h(\rho,1-\mu;s,\nu)
  \endmultline                                                       \tag{2.17}
$$
where $z_{2l}$ are defined by (2.27), $c(l)$ are taken from (2.10),
$$
   \Cal R_h(s,\nu;\rho,\mu)\!=
     \!2(4\pi)^{2s-2\nu-1}\widehat{\Phi}_N(2\nu\!+\!1\!-\!2s)\!\cdot\!
     \frac{\zeta(2\rho)\zeta(2\nu)}{\zeta(2\rho\!+\!2\nu)}
     \Cal Z(\rho\!+\!\nu;s,\mu) \tag{2.18}
$$
and the coefficients $h_0, h_1, g$ are defined by equalities
$$\multline
   h_0(r)\!\equiv\!h_0(r;s,\nu,\rho,\mu)\!=\!             \\
 =\!-i\!\int\limits_{(\Delta)}\gamma(w,1/2+ir)
   \gamma(\rho-w,\nu)\gamma(s-w,\mu)\cos\pi w\times   \\
   \times\!\Big(\!\cos\pi(\rho-w)\cos\pi(s-w)\!+
   \!\sin\pi\nu\sin\pi\mu\!\Big)
   \widehat{\Phi}_N (2w-2s-2\rho+2)dw,
  \endmultline                                                       \tag{2.19}
$$
$$\multline
   h_1(r)\!\equiv\!h_1(r;s,\nu,\rho,\mu)\!=\!             \\
 =\!-i\!\int\limits_{(\Delta)}\gamma(w,1/2+ir)
   \gamma(\rho-w,\nu)\gamma(s-w,\mu)\r{ch}(\pi r)\times   \\
   \times\!\widehat{\Phi}_N(2w-2s-2\rho+2)
   \Big(\!\cos\pi(s-w)\sin\pi\nu\!+
   \!\cos\pi(\rho-w)\sin\pi\mu\!\Big)dw,
  \endmultline                                                       \tag{2.20}
$$
$$\multline
   g(k)\!\equiv\!g(k;s,\nu,\rho,\mu)\!=\!                \\
 =\!\frac{2(2k-1)}{\pi i}
   \!\int\limits_{(\Delta)}
   \frac{\Gamma(k-1/2+w)}{\Gamma(k+1/2-w)}
   2^{2w-1}\gamma(\rho-w,\nu)\gamma(s-w,\mu)\times       \\
   \times\!\widehat{\Phi}_N(2w-2s-2\rho+2)
   \Big(\!\cos\pi(\rho-w)\cos\pi(s-w)\!+
   \!\sin\pi\nu\sin\pi\mu\!\Big)dw,
  \endmultline
                                                                     \tag{2.21}
$$
where $\Delta$ is taken with condition
$0<\Delta<\min(\r{Re}\,s,\r{Re}\,\rho)$.
\endproclaim
Before proving of this identity I give the following simple estimate.
\proclaim{Proposition 2.4}
Let $s$ and $\nu$ are fixed and $\sigma=\r{Re}\,s\in(1/2,1)$, $\r{Re}\nu
=1/2$; then for  any $d$ with condition $\ (d,c)=1$ we have for any
$\varepsilon>0$
$$
   \left|\Cal L_{\nu}\left(s,\frac dc\right)\right|\ll
   c^{1-\sigma+\varepsilon},                                 \tag{2.22}
$$
as an integer $c \to {+\infty}$.
\endproclaim

If $\r{Re}\,s \geqslant 1+\varepsilon$ then $\Cal
L_{\nu}(s,\frac{d}{c})$ is bounded uniformly in $c$.

For the case $\r{Re}\,s=-\varepsilon$ we have from the functional
equation (1.32)
$$
|\Cal L_{\nu}(s,\frac{d}{c})|\ll_{s,\nu} c^{1+2\varepsilon}.  \tag{2.23}
$$

Now (2.22) follows from the Phragm\'en-Lindel\"of principle.

After this we return to (2.17) and consider the triple sum
$$
  \zeta(2s)\zeta(2\rho)\sum_{n,m\geqslant 1}
  \frac{\tau_{\nu(n)}}{n^s}\frac{\tau_{\mu(m)}}{m^{\rho}}
  \sum_{c\geqslant 1}\frac 1cS(n,m;c)
  \Phi_N\left(\frac{4\pi\sqrt{nm}}c\right)                  \tag{2.24}
$$
where instead of $\Phi_N$ we can write the Mellin integral
$$
   \Phi_N(x)=\frac 1{\pi i}\int\limits_{(\delta)}
   \widehat {\Phi}_N(2w)x^{-2w}dw.                           \tag{2.25}
$$

Under our assumptions $\widehat {\Phi}_N(2w)$ is regular for
$\r{Re}\,w>-\min(\Delta,3/2)$ and we have the estimate (2.13) with
$N\geqslant 1$. As the consequence we have
$|\Phi_N(x)|\ll\min\left(x^3,x^{-1+\varepsilon}\right)$ for any
$\varepsilon>0$ (since $\Delta>3/2$ we can take $\delta=-3/2$ in (2.25)
if $x\to 0$; the integral is absolutely convergent for $\r{Re}\,w<3/2$
and when $x\to\infty$ we can take $\delta=1/2-\varepsilon,\ \
\varepsilon>0$).

The triple sum (2.24) equals to (we use (2.15))
$$
  Z^{(d)}(s,\nu;\rho,\mu|h)+Z^{(c,+)}(s,\nu;\rho,\mu|h)-
  \frac 1{2\pi}\sum_{l\in L}(-1)^lc(l)z_{2l}(s,\nu;\rho,\mu)
                                                                \tag{2.26}
$$
where
$$\multline
    z_{2l}(s,\nu;\rho,\mu)=               \\
 =\!\!\!\!\sum_{1\leqslant j\leqslant\nu_{2l}}\!\!\!\!
    \alpha_{j,2l}\Cal H_{j,2l}(s\!+\!\nu\!-\!1/2)
    \Cal H_{j,2l}(s\!-\!\nu\!+\!1/2)
    \Cal H_{j,2l}(\rho\!+\!\mu\!-\!1/2)
    \Cal H_{j,2l}(\rho\!-\!\mu\!+\!1/2)
  \endmultline
                                                                  \tag{2.27}
$$

On the other side, as the consequence of the Weil estimate, we have for
any $\varepsilon >0 $ the majorant
$$
  \sum_{n,m\geqslant 1}\sum_{c\geqslant 1}n^{-\sigma_1}m^{-\sigma_2}
  d(m)d(n)d(c)\frac{\sqrt{(n,m,c)}}{\sqrt c}
  \min \bigl((\frac{\sqrt{nm}}{c})^3,
  (\frac{c}{\sqrt{nm}})^{1-2\epsilon}\bigr)           \tag{2.28}
$$
for the series in (2.24) $ \big((n,m,c)$ is the greatest common
divisor of $n,m,c;$ and $\sigma_1=\r{Re}\,s,\ \
\sigma_2=\r{Re}\,\rho$  $\big)$. This series converges for
$\sigma_1,\sigma_2>5/4$.

It means the series in (2.24) converges absolutely and we can sum in any
order. Doing the summation over $n,m$ in the first line, we get for our
triple sum the expression
$$
  \frac{\zeta(2s)\zeta(2\rho)}{\pi i}
  \sum_{c\geqslant 1}\frac 1c\sum_{n,m\geqslant 1}
  \frac{\tau_{\nu(n)}\tau_{\mu(m)}}{n^sm^{\rho}}S(n,m;c)
  \int\limits_{(\delta)}\widehat{\Phi}_N(2w)
  \left(\frac{4\pi\sqrt{nm}}c\right)^{-2w}dw.
                                                                  \tag{2.29}
$$
Here for any fixed $c\geqslant 1$ we can take  (and take)
$\delta\in(0,1/2)$; under this condition our integral converges
absolutely and the inner double sum is
$$
   \sum\Sb a(\bmod c)\\ ad\equiv 1(\bmod c)\endSb
   \int\limits_{(\delta)} \left(\frac c{4\pi}\right)^{2w}
   \widehat{\Phi}_N(2w) \Cal L_{\nu}\left(s+w,\frac ac\right)
   \Cal L_{\mu}\left(\rho+w,\frac dc\right)dw.
                                                                    \tag{2.30}
$$

To receive the new expression we move the path of the integration to the
left on the line $\r{Re}\,w=\delta_1$ so that
$\sigma_1+\delta_1,\ \ \sigma_2+\delta_1<0$. Four poles at
$w=1-s\pm(\nu-1/2)$ and $w=1-\rho\pm(\mu-1/2)$ (where
$\Cal L_{\nu}\left(s+w,\frac ac\right)$ and
$\Cal L_{\mu}\left(\rho+w,\frac dc\right)$ have the residues
$c^{-2\pm(2\nu-1)}\zeta\big(1\pm(2\nu-1)\big)$ and
$c^{-2\pm(2\mu-1)}\zeta\big(1\pm(2\mu-1)\big)$ correspondingly)
give us the terms
$$\multline
     2\zeta(2s)\zeta(2\rho)\!\sum_{c\geqslant 1}\!
     \frac{\zeta(2\nu)(4\pi)^{2s-2\nu-1}}{c^{2s}}
     \widehat{\Phi}_N(2\nu\!+\!1\!-\!2s)
     \!\!\!\!\!\!\sum_{ad\equiv 1(\bmod c)}\!\!\!\!\!\!
     \Cal L_{\mu}\left(\!\rho\!+\!\mu\!+\!1/2\!-\!s,\frac dc\!\right)\!+  \\
   + \{\text{the same with $\nu$ replaced by $1-\nu$}\}+                  \\
   + 2\zeta(2s)\zeta(2\rho)\!\sum_{c\geqslant 1}\!
     \frac{\zeta(2\mu)(4\pi)^{2\rho-2\mu-1}}{c^{2\rho}}
     \widehat{\Phi}_N(2\mu+1-2\rho)
     \!\!\!\!\!\!\sum_{ad\equiv 1(\bmod c)}\!\!\!\!\!\!
     \Cal L_{\nu}\left(\!s\!+\!\mu\!+\!1/2\!-\!\rho,\frac ac\!\right)\!+  \\
   + \{\text{the same with $\mu$ replaced by $1-\mu$}\}.
  \endmultline
                                                                \tag{2.31}
$$
Under our assumption we have $|\r{Re}\,(\rho-s)|\leqslant\varepsilon$;
so
$$
   \left|\Cal L_{\mu}\left(\rho+\nu+1/2-s,\frac dc\right)\right|\ll
   c^{\varepsilon}
                                                                 \tag{2.32}
$$
and all series in (2.31) are convergent absolutely and they define the
regular functions in $\rho$ and $s$. If $\r{Re}\,(\rho-s)>0$ then we
have
$$
   \sum_{(d,c)=1}\Cal L_{\mu}\left(\rho+\nu+1/2-s,\frac dc\!\right)=
   \sum_{n\geqslant 1}
   \frac{\tau_{\mu}(n)S(0,n;c)}{n^{\rho+\nu+1/2-s}}
                                                                  \tag{2.33}
$$
and
$$
   \sum_{c\geqslant 1}\frac 1{c^{2s}}\sum_{(d,c)=1}
   \Cal L_{\mu}\left(\rho+\nu+1/2-s,\frac dc\!\right)=
   \frac 1{\zeta(2s)}\Cal Z(\rho+\nu;s,\mu).
                                                                   \tag{2.34}
$$

The analitical continuation of these equalities gives four terms with
$\widehat{\Phi}_N$ on the right side (2.17). In the integral on the line
$\r{Re}\,w=-\delta_1$ we use the functional equation (1.32) and we come
to the expression
$$\multline
   \zeta(2s)\zeta(2\rho)\sum_{c\geqslant 1}\frac 1c
   \sum_{n,m\geqslant 1}
   \frac{\tau_{\nu(n)}}{n^{\rho}}\frac{\tau_{\mu(m)}}{m^s}\times   \\
   \times\left\{S(n,m;c)\varphi_0\left(\frac{4\pi\sqrt{nm}}c\right)+
          S(n,-m;c)\varphi_1\left(n,-m;c\right)\right\}
  \endmultline
                                                                   \tag{2.35}
$$
where for $x>0\ \ \varphi_0$ and $\varphi_1$ are defined by the
integrals
$$\multline
     \varphi_0(x)\equiv\varphi_0(x;s,\nu,\rho,\mu)=
     \frac 1{i\pi}x^{2s+2\rho-2}\int\limits_{(\delta_1)}
     \gamma(1-s-w,\nu)\gamma(1-\rho-w,\mu)\times                    \\
  \times\Big(\cos\pi(s+w)\cos\pi(\rho+w)+\sin\pi\nu\sin\pi\mu\Big)
     \widehat{\Phi}_N(2w)x^{2w}dw,
  \endmultline
                                                                   \tag{2.36}
$$
$$\multline
     \varphi_1(x)\equiv\varphi_1(x;s,\nu,\rho,\mu)=-
     \frac 1{i\pi}x^{2s+2\rho-2}\int\limits_{(\delta_1)}
     \gamma(1-s-w,\nu)\gamma(1-\rho-w,\mu)\times                    \\
  \times\Big(\cos\pi(\rho+w)\sin\pi\nu+\cos\pi(s+w)\sin\pi\mu\Big)
     \widehat{\Phi}_N(2w)x^{2w}dw;
  \endmultline
                                                                     \tag{2.37}
$$
it is possible take here any $\delta_1$ with
$-3/2\leqslant\delta_1<\min(1-\r{Re}\,s,1-\r{Re}\,\rho)$.

Both integrals (2.36) and (2.37) and the triple series (2.35) converge
absolutely if $5/4<\r{Re}\,s,\ \r{Re}\,\rho<5/4+\varepsilon$ for some
small $\varepsilon$ (it is assumed
$\r{Re}\,\nu=\r{Re}\,\mu=1/2,\ \ \nu\ne1/2,\mu\ne1/2$). Really, we have
$$
   \varphi_j(x)=O\Big(\max\big(x^{2\r{Re}\,s},x^{2\r{Re}\,\rho}
     \big)\Big) \quad \text{as} \quad x\to 0
                                                                  \tag{2.38}
$$
and for any $\varepsilon_1>0$
$$
  \varphi_j(x)=O\Big(x^{2\r{Re}\,s+2\r{Re}\,\rho-7+\varepsilon_1}
                \Big), \qquad x\to +\infty.
                                                                 \tag{2.39}
$$

Since (the Weil estimate again)
$|S(n,m;c)|\leqslant c^{1/2}d(c)(m,n,c)^{1/2}$ we have the majorant
$$
  \frac{d(n)d(m)d(c)}{m^{\sigma_1}n^{\sigma_2}}\cdot
  \min\left(\left(\frac{\sqrt{mn}}c\right)^{\sigma_2},
  \left(\frac c{\sqrt{mn}}\right)^{2\sigma_1+2\sigma_2-7+\varepsilon_1}
  \right)
                                                                  \tag{2.40}
$$
where $\sigma_1=\r{Re}\,s,\ \ \sigma_2=\r{Re}\,\rho$ and, for
definiteness, we assume $\sigma_1>\sigma_2$. Note
$0<\sigma_1-5/4,\ \sigma_2-5/4<\varepsilon$ for some small
$\varepsilon>0$; under this condition the triple majorising series with
terms (2.40) converges.

It means we can change the order of summation in (2.35).

Now we consider the inner sum of the Kloosterman sums. We would use
identities (1.13) and (1.25) with the given $h$ on the right side.

The reason of this manner is very simple: the conditions of my foretrace
formulae with a given taste function in front of the Kloosterman sums
are farther from necessary than the similar conditions for the case of
the given taste function on the other side.

\proclaim{Proposition 2.4}
Let $\r{Re}\,s,\r{Re}\,\rho\in(5/4,5/4+\varepsilon)$ for some small
$\varepsilon>0$; let $h_0(r)$ and $g(k)$ are defined by (2.19) and
(2.21); then for any integers $m,n\geqslant 1$ we have for
$\r{Re}\,\nu=\r{Re}\,\mu=1/2$
$$\multline
     \sum_{c\geqslant 1}\frac{1}{c} S(n,m;c)
     \varphi_0\left(\frac{4\pi\sqrt{mn}}c\right)=
     \sum_{j\geqslant 1}\alpha_jt_j(n)t_j(m)h_0(\varkappa_j)+  \\
   + \frac 1{\pi}\!\!\int\limits_{-\infty}^{\infty}\!\!
     \frac{\tau_{1/2\!+\!ir}(n)\tau_{1/2\!+\!ir}(m)}
                {|\zeta(1/2\!+\!ir)|^2}
     h_0(r)\,dr\!+\!\sum_{k\geqslant 6}g(k)\!\!
     \sum_{1\leqslant j\leqslant 6}\!\!
     \alpha_{j,2k}t_{j,2k}(n)t_{j,2k}(m).
  \endmultline
                                                              \tag{2.41}
$$
\endproclaim
The function $h_0(r)$ coincides with the integral transform
$$\multline
      h_0=\pi\!\int\limits_0^{\infty}\!k_0(x,1/2\!+\!ir)
          \varphi_0(x)\frac{dx}x\!=                                    \\
   =  \!i\!\!\int\limits_{(\Delta)}\!\!
      \gamma(s\!+\!\rho\!+\!w\!-\!1,1/2\!+\!ir)
       \gamma(1\!-\!s\!-\!w,\nu)\gamma(1\!-\!\rho\!-\!w,\mu)
       \cos\pi(s\!+\!\rho\!+\!w)\times                                 \\
   \times\!\Big(\!\cos\pi(s\!+\!w)\cos\pi(\rho\!+\!w)\!+
       \!\sin\pi\nu\sin\pi\mu\!\Big)
      \widehat{\Phi}_N (2w)\,dw
  \endmultline
                                                                   \tag{2.42}
$$
(if $|\r{Im}\,r|$ be small and $\Delta$ be taken with condition
$1<\Delta+2\sigma<5/4,\ \ \sigma=\r{Re}\,s=\r{Re}\,\rho$, then we can
integrate in (2.36) under the sign of integral over $w$).

One can see that $h_0(r)$ is regular in $r$ at least for
$|\r{Im}\,r|<2\sigma+\Delta-1$ and we can take any $\Delta$ with
$1-\sigma-\Delta>0$; the strip of regularity is the strip
$|\r{Im}\,r|<\sigma\approx 5/4$. If $r\to\infty$ inside of this strip we
have $|h_0(r)|=O\left(|r|^{-5/2}\right)$ uniformly in
$\r{Im}\,r,\ \ |\r{Im}\,r|<\sigma$.
To find this estimate we move the path of integration to the left on the
line $\r{Re}\,w=\Delta_1,\ \ \Delta=-\sigma_1-\sigma_2+1/4\ (\approx-9/4)$.
The residue at $w=-s-\rho+1-ir$ is $O\left(|r|^{-N-5/2}\right)$ as it
follows from the Stirling expansion; the same is true for the residue at
$w=-s-\rho-ir$. On the line $\r{Re}\,w=\Delta_1$ we have
$$\multline
   \big|\gamma(s+\rho+w-1,1/2+ir)\cos\pi(s+\rho+w)\big|\ll           \\
\ll\Big(\big||r|^2-|w|^2\big|+1\Big)^{\sigma_1+\sigma_2+\Delta_1-3/2}
   \exp\Big(\min\big(0,\pi(|w|-|r|)\big)\Big);
  \endmultline
                                                       \tag{2.43}
$$
so the part of this integral with $|w|>|r|$ is exponentially small and
it is sufficient consider the integral with $|w|\ll |r|$.

If $|w|\leqslant 2|r|$ then the integrand in (2.36) may be estimated as
$$
   \ll\Big(\big||r|^2-|w|^2\big|+1\Big)^{\sigma_1+\sigma_2+\Delta-3/2}
   \big(|w|+1\big)^{1-\sigma_1-\sigma_2-\Delta-N}.
                                                        \tag{2.44}
$$

Let $\Delta_1=\Delta+\sigma_1+\sigma_2\ \
(\sigma_1=\r{Re}\,s, \sigma_2=\r{Re}\,\rho)$; if
$0\leqslant\Delta_1\leqslant 1/4$ then we have
$$\multline
     \big|h_0(r)\big|\ll\int\limits_0^{|r|/2}|r|^{2\Delta_1-3}
     (\eta+1)^{-\Delta_1-1-N}d\eta+                             \\
   + \int\limits_{|r|/2}^{|r|-1}|r|^{-5/2-N}
     \big(|r|-\eta\big)^{\Delta_1-3/2}d\eta+
     \int\limits_{|r|-1}^{|r|+1}|r|^{-5/2-N}d\eta\ll |r|^{-5/2}
  \endmultline
                                                          \tag{2.45}
$$
for any $N\geqslant 1$.

So all conditions of theorem 1.2 are fullfilled and first two terms on
the right side (2.41) are equal to the sum of the Kloosterman sums with
the taste function $\widetilde{\varphi}_0$,
$$
   \widetilde{\varphi}_0(x)=\int\limits_{-\infty}^{\infty}
   k_0(x,1/2+ir)h_0(r)\,d\chi(r),
                                                          \tag{2.46}
$$
and with term which contains $\delta_{n,m}$.

One can see the integral (2.46) gives that component of $\varphi_0$ which is
orthogonal to all Bessel's functions of the odd integer order. Let this
component is $\varphi_0^{(H)}$; then
$$
   \varphi_0=\varphi_0^{(H)}+\Cal N,                        \tag{2.47}
$$
where $\Cal N$ is the corresponding Neumann series,
$$
  \Cal N(x)=\sum_{k=1}^{\infty}2(2k-1)J_{2k-1}(x)
  \int\limits_0^{\infty}\varphi_0(y)J_{2k-1}(y)\frac{dy}{y}=
  \sum_{k=1}^{\infty}g(k)J_{2k-1}(x)
                                                              \tag{2.48}
$$
with $g(k)$ from (2.21) (again we can integrate term by term; after the
calculation we change $w$ by $w-s-\rho+1$).

Taking in (2.21) $\Delta=-1+\varepsilon$ with small $\varepsilon>0$
(it is possible, since $\r{Re}\,(w-s-\rho+1)>-\frac 52$ on this line) we
come to the estimate
$$
   |g(k)|\ll k^{-2+2\varepsilon}.                              \tag{2.49}
$$

It rests use the Peterson identity and we get (2.41), since the singular
term wich $\delta_{n,m}$ disappears (see [2]).

\proclaim{Proposition 2.5}
Under the same assumptions we have
$$\multline
     \sum_{c\geqslant 1}\frac 1cS(n,-m;c)
     \varphi_1\left(\frac{4\pi\sqrt{nm}}c\right)=
     \sum_{j\geqslant 1}\alpha_j\varepsilon_j
     t_j(n)t_j(m)h_1(\varkappa_j)+                               \\
   + \frac 1{\pi}\!\!\int\limits_{-\infty}^{\infty}\!\!
     \tau_{1/2\!+\!ir}(n)\tau_{1/2\!+\!ir}(m)
     h_1(r)\frac{dr}{|\zeta(1+2ir)|^2}
  \endmultline
                                                          \tag{2.50}
$$
with $h_1$ from (2.20).
\endproclaim

We define $h_1$ by the integral transform
$$
   h_1(r)\equiv h_1(r,s,\nu;\rho,\mu)=\pi\int\limits_0^{\infty}
         \varphi_1(x)k_1(x,1/2+ir)\frac{dx}{x},               \tag{2.51}
$$
where for $x>0$
$$
k_1(x,\nu) =\frac{2}{\pi}\,\sin \pi \nu  K_{2\nu -1}(x)        \tag{2.52}
$$

The term by term integration gives for this function the representation
(2.20). The function $h_1$ is satisfying to the same conditions as
$h_0$; it is sufficient to have (2.50) if $\varphi_1$ be defined by the
integral transform (1.26) of $h_1$. But this integral coincides with
$\varphi_1$ if $h_1$ be defined by (2.51).

The union of all previous equalities and two last propositions gives us
(2.21) if $\r{Re}\,s,\ \r{Re}\,\rho>5/4$ and both differences
$(\r{Re} \,s-5/4), (\r{Re}\,\rho-5/4)$ are sufficiently small.

\subhead
2.3. The meromorphic continuation
\endsubhead

Let us introduce the following notation
$$
\Cal R (s,\nu;\rho,\mu|h)=2\frac{\zeta(2s-1)\zeta(2\nu)}{\zeta(2-2s+2\nu)}
   \Cal Z(\rho;\mu,1-s+\nu)h\big(i(s-\nu-1/2)\big)
                                                                 \tag{2.53}
$$
(it is the residue at the point $r=i(s-\nu-1/2)$ of the integrand in the
definition $Z^{(c,+)}(s,\nu;\rho,\mu|h)$ with the additional factor
$2\pi i$).

\proclaim{Proposition 2.6.}
Let $h$ be a regular even function in the sufficiently wide strip
$|\r{Im}\,r| \leqslant\Delta,\ \Delta>1/2$ and let $|h(r)|\ll |r|^{-B}$ for
some $B>5/2$ as $r\to\infty$ in this strip. Then for
$1/2 \leqslant\r{Re}\,s<1,  1/2\leqslant\r{Re}\,\rho<1$ the meromorphic
continuation of $Z^{(c,+)}(s,\nu;\rho,\mu|h)$ is given by the equality
$$\multline
     Z^{(c,+)}(s,\nu;\rho,\mu|h)=Z^{(c)}(s,\nu;\rho,\mu|h)+     \\
   +       \Cal R(s,\nu;\rho,\mu|h)\!+\Cal R(s,\!1\!-\!\nu;\rho,\mu|h)\!+
           \Cal R(\rho,\mu;s,\nu|h)\!+\Cal R(\rho,\!1\!-\!\mu;s,\nu|h).
  \endmultline
                                                               \tag{2.54}
$$
\endproclaim

Really,
$$
   Z^{(c,+)}(s,\nu;\rho,\mu|h)=
   \frac 1{\pi}\int\limits_{-\infty}^{\infty}
   \frac{\Cal Z(s;\nu,1/2+ir)\Cal Z(\rho;\mu,1/2+ir)}
        {\zeta(1+2ir)\zeta(1-2ir)}h(r)\,dr
                                                         \tag{2.55}
$$
and we have the Cauchy integral, because $\zeta$ has only a simple pole.

The product $\Cal Z(s;\nu,1/2+ir)$ has the poles at the points
$r_j,\ 1\leqslant j\leqslant 4$, with
$$
  ir_1=1/2-s+\nu,\quad ir_2=3/2-s-\nu,\quad ir_3=-r_1,\quad ir_4=-r_2.
$$

When $\r{Re}\,s>1$ the points $r_1,r_2$ are lying above the real axis
and if $\r{Re}\,s<1,$ they are below the same. Now one can deform the
path of integration (see picture; the deformation must be so small that
the functions $\zeta(1\pm 2ir)$ have no zeros in side the lines; it is
possible, since the Riemann zeta--function has no zeros on the line
$\r{Re}\,s=1$).

\vskip 5truecm
\centerline{The path of integration $l$.}

Firstly we fix $\rho,\;\;\r{Re}\,\rho>1;$ if $\r{Re}\,s>1$ we have
$$
Z^{(c,+)}(s,\nu;\rho,\mu|h)=\frac 1{\pi}
\int\limits_{(l)} (...)\,dr +\Cal R(s,\nu;\rho,\mu|h)+\Cal R(s,1-\nu;\rho,\mu|h);
                                                      \tag{2.56}
$$
it is the result of the direct calculation of the residues. But for
$\r{Re}\,s<1$ we have $\r{Im}\,r_1,\ \r{Im}\,r_2,<0$ and
$\r{Im}\,r_3,\ \r{Im}\,r_4>0$, so we can integrate over the real axis,
$\int\limits_l(\ldots)\,dr=\int\limits_{-\infty}^{\infty}(\ldots)\,dr$.

Doing the same with a fixed $s,\ \r{Re}\,s<1,$ we get the continuation
in the variable $\rho$ (for what it is sufficient to permute the pair of
variables); it gives the equality (2.54).

As the consequence we can rewrite Theorem 2.2 and get our functional
equation inside of the strip $1/2\leqslant\r{Re}\,s,\,\r{Re}\,\rho<1$.

\proclaim{Theorem 2.3}
Let $h(r)$ be the even regular function in the strip
$|\r{Im}|\,r\leqslant\Delta$ for some $\Delta>5/2,\ \ h(i/2)=h(3i/2)=0$ and
$|h(r)|\ll |r|^{-B}$ for any fixed positive $B$ when $r\to\infty$ inside of
this strip.

Let $\Phi_N$ be defined, for this given $h$, by (2.11) with $N
\geqslant1$ and $l_1\geqslant 1$.

Then for any $s,\nu,\rho,\mu$ with
$\r{Re}\,\nu=\r{Re}\,\mu=1/2,\ \nu,\mu\neq 1/2$ and
$1/2\leqslant\r{Re}\,s,\r{Re}\,\rho<1$ we have
$$\multline
     Z^{(d)}(s,\nu;\rho,\mu|h,0)+Z^{(c)}(s,\nu;\rho,\mu|h,0)=          \\
   = Z^{(d)}(\rho,\nu;s,\mu|h_0,h_1)+Z^{(c)}(\rho,\nu;s,\mu|h_0+h_1)+
     Z^{(r)}(\rho,\nu;s,\mu|g)+                                        \\
   +\frac{1}{2\pi}\sum_{l\in L}(-1)^l c(l) z_{2l}(s,\nu;\rho,\mu)+ \\
   + \!\Cal R_h(s,\nu;\rho,\mu)\!+\!\Cal R_h(s,\!1\!-\!\nu;\rho,\mu)\!+
     \!\Cal R_h(\rho,\mu;s,\nu)\!+\!\Cal R_h(\rho,\!1\!-\!\mu;s,\nu)\!-  \\
   - \!\Big(\!\Cal R(s,\nu;\rho,\mu|h)\!+\!\Cal R(s,\!1\!-\!\nu;\rho,\mu|h)\!+
     \!\Cal R(\rho,\mu;s,\nu|h)\!+\!\Cal R(\rho,\!1\!-\!\mu;s,\nu|h)\!\Big)\!+     \\
   + \!\Cal R(\rho,\nu;s,\mu|h_0\!+\!h_1)\!+
     \!\Cal R(\rho,\!1\!-\!\nu;s,\mu|h_0\!+\!h_1)\!+\\
  +   \!\Cal R(s,\mu;\rho,\nu|h_0\!+\!h_1)\!+
     \!\Cal R(s,\!1\!-\!\mu;\rho,\nu|h_0\!+\!h_1)
  \endmultline
                                                        \tag{2.57}
$$
where $h_0,h_1$ and $g$ are defined by (2.19)--(2.21) with $N\geqslant 1,
\;l_1\geqslant 1$ and $\Cal R_h,\,\Cal R$ are given by the equalities (2.18),
(2.53).
\endproclaim

\head
\S 3. SPECIALIZATION OF THE MAIN FUNCTIONAL EQUATION
\endhead

\subhead
3.1. The choice of variables
\endsubhead

To estimate $|\zeta (1/2 + it)|$ for large positive $t$ we use the
following specialization of the functional equation (2.57):

\noindent
i)we take in this identity
$$
s=\mu=1/2, \nu = 1/2 + i t, \rho = 1/2 + i \tau,         \tag{3.1}
$$
where $t$ and $\tau$ are large, $t \to +\infty$ and $\tau \gg t^4$; we
assume that for some fixed $\delta >0$ the parameter $t$ is larger than
$\tau^{\delta}$ and at the end we will take $t\approx \tau^{1/8}$).
\bigskip
\noindent
ii)we assume that for all real $r$
$$
h(r) > 0                                  \tag{3.2}
$$
and this $h$ satisfies to all conditions of Theorem 2.3;
\bigskip
\noindent
iii) for the choosen $h$ we suppose
$$
\Phi (x) = \int_{-\infty}^{\infty} k_0(x,1/2+iu) h(u)\,d\chi(u)-
\sum_{l \in L} c(l) J_{2l-1}(x),
\tag{3.3}
$$
where $L$ contains 5 elements
$$
L = \{2,3,4,5,7\}                                               \tag{3.4}
$$
(these concret elements are taken since the spaces of the corresponding
cusp forms of the weight $2l$, $l \in L$, are empty; so sums $z_{2l}$ on the
right side (2.57) are zeroes).

The coefficients $c(l)$ in (3.3) are defined by equalities (2.10) with
$N=5$, so the Mellin transform of $\Phi$ is $O(|w|^{2 \r{Re} w - 6})$ if
$w \to \infty$ and $\r{Re} w$ is fixed.

For our specialzation we have on the left side (2.57) the very fast
convergent series
$$
\sum_{j\geqslant 1} \alpha_j |\Cal H_j(\nu)|^2 \Cal H_j^2(\rho)
h(\varkappa_j)                                 \tag{3.5}
$$
and the similar integral over the continuous spectrum.

\subhead
3.2. The averaging
\endsubhead
Let $T_0$ be sufficiently large and $T=T_0^{1-\varepsilon}$ with small
fixed $\varepsilon >0$. We suppose
$$
\omega_T(\rho)=\frac{1}{4}\bigl(\cos \frac{\pi}{2T}(\rho - 1/2
-iT_0)\bigr)^{-1}                                              \tag{3.6}
$$
and consider the average of our functional equation (2.57) with
variables (3.1) over variable $\rho$ with weight $\omega_T(\rho)$.

The average of $Z^{(d)}$ is the integral
$$
\frac{1}{i}\int\limits_{(1/2)} \omega_T(\rho) \sum_{j\geqslant 1} \alpha_j
|\Cal H_j(\nu)|^2 \Cal H_j^2(\rho) h(\varkappa_j) \,d\rho           \tag{3.7}
$$
where we can integrate term by term under our conditions for $h$.

All integrals
$$
\frac{1}{i}\int\limits_{(1/2)}\omega_T (\rho) \Cal H_j^2 (\rho)\,d\rho
                                                                \tag{3.8}
$$
may be calculated over the line $ \r{Re}\rho = 2$, where $\Cal H_j$ is
the absolutely convergent series and we have
$$
\frac{1}{i} \int\limits_{(1/2)} \omega_T(\rho)\Cal H_j^2(\rho)\,d\rho
=T \sum_{n,m \geqslant 1} \frac{t_j(n)t_j(m)}{(nm)^{1/2+iT_0}}
((nm)^T +(nm)^{-T})^{-1}                                     \tag{3.9}
$$
(we use the so called Ramanujan integral here).

The series on the riht side (3.9) equals to
$$
T\left(1+ O(2^{-T})\right)                                    \tag{3.10}
$$
and the main term is positive.

Doing the same in the integral over the continuous spectrum we come to
the following expression for the average on the left side
$$
\multline
\frac{1}{i} \int\limits_{(1/2)} \omega_T(\rho)\lbrace
Z^{(d)}(1/2,\nu;\rho,1/2|h,0)+Z^{(c)}(1/2,\nu;\rho,1/2|h,0)\rbrace\,d\rho=\\
=T\lbrace\sum_{j\geqslant 1}\alpha_j |\Cal H_j(\nu)|^2 h(\varkappa_j)
+\frac{1}{\pi}\int_{-\infty}^{\infty}\frac{|\zeta(\nu+ir)\zeta(\nu
-ir)|^2}{|\zeta(1+2ir)|^2} h(r)\,dr\rbrace +o(1)
\endmultline
                                                         \tag{3.11}
$$

Really, we can integrate over $\rho$ under the sign of integration in
$r$. When we move the path of integration to the line $\r{Re}\,\rho = 2$
the additional terms will be appeared from the poles of
$\zeta^2(\rho+ir)\zeta^2(\rho-ir)$. But the products $h(r)\omega_T(1\pm
ir)$, $h(r) \omega'_T(1\pm ir)$ are $O(T^{-M})$ for any positive $M$.

The upper bound for this expression will be received by the same
averaging on the right side.

Here we have the products $\Cal H_j(\rho+\nu -1/2) \Cal H_j(\rho -\nu +1/2)
\Cal H_j^2(1/2)$ in the sum over the discret spectrum and
$\Cal Z(\rho;\nu,1/2+ir) \Cal Z(1/2;1/2,1/2+ir)$ in the integral, but
coefficients $h_0$ and $h_1$ are depending in $\rho$ also and we must
receive the explicit formulas for them.

For this purpose we will express these coefficients in terms of the
integrals with two hypergeometric functions.

For large values $\tau,t$ and $r$ we can write the full asymptotic
expansions for these hypergeomtric functions. It allows us receive the
sufficiently exact estimates for the average on the right side.

As the result we have the upper bound for (3.11) in the form
$$
\sum_{\varkappa_j\leqslant \sqrt {T_0}}\alpha_j \Cal H_j^2(1/2)
+\sum_{\sqrt {T_0}< \varkappa_j\leqslant 2T_0 }\alpha_j \frac{\log^2
\varkappa_j}{\varkappa_j^2} \Cal H_j^2(1/2) +O(T_0^{1+\varepsilon})
                                                         \tag{3.12}
$$
plus the similar sum over regular cusp forms (see \S\,9).

These mean values are known and it gives the estimate
$O(T^{1+\varepsilon})$ for (3.11); of course, it is sufficient to prove
(0.2) and (0.3).

\head
\S\,4.NEW REPRESENTATIONS FOR THE COEFFICIENTS
\endhead
\subhead
4.1. The definition of the kernels $A_{jk}$
\endsubhead
The Mellin integrals (2.19)--(2.21) are very convenient in general
theory, but these ones are not good for the special cases when one (ore more)
variable is large. There is other form which is more suitable for the
estimation in these cases.

We define four kernels $A_{jk}(r,x;s,\nu),\ \ 0\leqslant j,k\leqslant 1$,
by the equalities
$$\align
     (4\pi)^{2s-1}A_{0j}(r,x;s,\nu)&=\pi x\int\limits_0^{\infty}
      k_0(y,1/2+ir) k_j(xy,\nu)y^{2s-1}dy,\ \ \ j=0,1;       \tag{4.1}\\
     (4\pi)^{2s-1}A_{1j}(r,x;s,\nu)&=\pi x\int\limits_0^{\infty}
      k_1(y,1/2+ir) k_{1-j}(xy,\nu)y^{2s-1}dy,\ \ \ j=0,1;   \tag{4.2}
  \endalign
$$
here $k_0(x,\nu)$ and $k_1(x,\nu)$ are defined by (1.17) and (2.51)
correspondingly.

The integral for $A_{00}$ is absolutely convergent in the strip
$1/2>\r{Re}\,s>|\r{Im}\,r|+|\r{Re}\,(\nu-1/2)|$ and all others are
convergent in the halfplane $\r{Re}\,s>|\r{Im}\,r|+|\r{Re}\,(\nu-1/2)|$.

Note the kernels $A_{0j}$ have a sense for $r=\pm i(l-1/2)$ with an
integer $l$ since
$$
k_0(y,l) = k_0(y,1-l) = (-1)^l J_{2l-1}(y)
$$
and the corresponding integrals (4.1) with $l=\pm i(l-1/2)$ are convergent
for $1/2>\r {Re} s > -l+1/2 +|\r {Re} (\nu-1/2)|$.

\subhead
4.2. The explicit formulas
\endsubhead
Here (and for what will be later on) it is convenient rewrite the
explicit formulas for $A_{jk}$ from [8] (with the additional
representations which are the immediate consequence of the Kummer
relations between the Gauss hypergeometric functions of arguments
$x,\,1-x,\,-\frac{1}{x},...$).

It is very essentially for us that the function
$$
\Cal F(x;s,\nu,r) =|x|^{ir+1/2}\,(1-x)^s\,F(s+\nu-1/2+ir,s-\nu+1/2
+ir;1+2ir;\,x)                                    \tag{4.3}
$$
(here the usual notation for the Gauss hypergeometric function is used;
it is assumed that $x$ be real and $x<1$) is a solution of the
differential equation
$$
\frac{d^2 \Cal
F}{dx^2}+\biggl(\frac{r^2}{x^2\,(1-x)}-\frac{(s-1/2)^2}{x\,(1-x)^2}
+\frac{(\nu-1/2)^2}{x\,(1-x)})+\frac{1-x+x^2}{4x^2 \,(1-x)^2}\biggr)\Cal
F= 0.                                                               \tag{4.4}
$$

It gives the possibility find the asymptotic expansions for $\Cal F$ if
at least one parameter be large (see \S 5).

The first assertion about $A_{jk}$ follows directly from the definition.
\proclaim {Proposition 4.1} Let $A_{jk}$ are defined by two integrals
(4.1) and (4.2). Then these functions for any fixed positive $x$ and for
any fixed $\nu$ with $\r{Re} \nu =1/2$ are regular in $s$ and $r$ if $\r
{Re} s >|\r {Im} r|.$ \endproclaim

The formulas for $A_{jk}$ below give us the meromorphic continuation of
these functions; the result of this continuation we will denote by the
same symbol.

\proclaim{Proposition 4.2} Let us introduce the notations
$$
a=s+\nu-1/2+ir, \,b=s-\nu+1/2+ir, \, c=1+2ir, \tag{4.5}
$$
$$
a'=s+\nu-1/2+ir, \,b'=s-\nu+1/2+ir, \,c'=1-2ir. \tag{4.6}
$$
With these notations for $0 \leqslant x < 1$ we have
$$
\multline
(2\pi)^{2s-1} \,A_{00}(r,x;s,\nu) = \\
=\frac{1}{2\cos \pi \nu}
\biggl(\sin \pi(s+\nu)\frac{\Gamma(a) \Gamma(a')}{\Gamma(2\nu)}
\,x^{2\nu}\,F(a,a';c;\,x^2)+\\
+ \sin \pi (s-\nu)\frac{\Gamma(b) \Gamma(b')}{\Gamma(2-2\nu)}
\,x^{2-2\nu}\,F(b,b';2-2\nu; x^2)\biggr)
\endmultline                \tag{4.7}
$$
\endproclaim

\proclaim{Proposition 4.3} Let $x>1$; then we have
$$
\multline
(2\pi)^{2s-1} A_{00}(r,x;s,\nu)=\frac{ix^{1-2s}}{2 \sinh \pi
r}\biggl(\cos \pi(s+ir) \frac{\Gamma(a) \Gamma(b)}{x^{c-1}
\Gamma(c)}F(a,b;c;\frac{1}{x^2})-\\
-\cos \pi(s-ir)\frac {\Gamma(a')
\Gamma(b')}{x^{c'-1}\Gamma(c')}F(a',b';c';\frac{1}{x^2})\biggr)
\endmultline                         \tag{4.8}
$$
\endproclaim

The next representation gives the union of the previous ones.
\proclaim{Proposition 4.4} For all $x\geqslant 0$
$$
\multline
\noindent
(2\pi)^{2s-1} A_{00}(r,x;s,\nu)=\sin \pi s \, \Gamma(2s-1) |x^2-1|^{1-2s}
x^{2\nu} F(1-b,1-b';2-2s;1-x^2)\\
+\frac{\Gamma(a) \Gamma(a') \Gamma(b) \Gamma(b')}{2\pi \cos \pi s \,
\Gamma(2s)}(\cosh^2 \pi r +\sin^2 \pi \nu -\sin^2 \pi s) x^{2\nu}
F(a,a';2s; 1-x^2)
\endmultline                                  \tag{4.9}
$$
\endproclaim

Two kernels $A_{01}$ and $A_{10}$ are exponentially small when $r \to
\pm \infty$ (and $s,\nu$ are $ o(|r|)$.

\proclaim{Proposition 4.5} For all $x>0$ we have
$$
\multline
(2\pi)^{2s-1} A_{01}(r,x;s,\nu)=\frac {ix^{1-2s} \sin \pi \nu}{2 \sinh
\pi r}\biggl(\frac {\Gamma(a) \Gamma(b)}{x^{c-1} \Gamma(c)} F(a,b;c;
-\frac{1}{x^2})\\
-\frac{\Gamma(a') \Gamma(b')}{x^{c'-1}
\Gamma(c')} F(a',b';c';-\frac{1}{x^2}\biggr)
\endmultline                            \tag{4.10}
$$
and the same kernel equals to
$$
\multline
\frac{1}{2\cos \pi \nu}\biggl( \sin \pi (s+\nu) \frac{\Gamma(a)
\Gamma(a')}{\Gamma(2\nu)} x^{2\nu} F(a,a';2\nu;-x^2)+\\
+\sin \pi (s-\nu) \frac{\Gamma(b) \Gamma(b')}{\Gamma(2-2\nu)} x^{2-2\nu}
F(b,b';2-2\nu;-x^2)\biggr).
\endmultline                                   \tag{4.11}
$$
\endproclaim
\proclaim{Proposition 4.6} For all $x\geqslant 0$ we have
$$
(2\pi)^{2s-1} A_{10}(r,x;s,\nu) =\cosh \pi r \sin \pi \nu
\frac{\Gamma(a) \Gamma(a') \Gamma(b) \Gamma(b')}{\pi \Gamma(2s)} x^{2\nu}
F(a,a';2s;1-x^2);                                   \tag{4.12}
$$
at the same time this kernel for $0\leqslant x <1$ equals to
$$
\frac{\cosh \pi r}{2\cos \pi \nu}\biggl(\frac{\Gamma(a)
\Gamma(a')}{\Gamma(2\nu)}x^{2\nu} F(a,a';2\nu;x^2)- \frac{\Gamma(b)
\Gamma(b')}{\Gamma(2-2\nu)}x^{2-2\nu} F(b.b';2-2\nu;x^2)\biggr)
\tag{4.13}
$$
and
$$
\frac{ix^{1-2s} \sin \pi \nu}{2 \sinh \pi r}\biggl( \frac{\Gamma(a)
\Gamma(b)}{x^{c-1}\Gamma(c)} F(a,b;c;\frac{1}{x^2})- \frac{\Gamma(a')
\Gamma(b')}{x^{c'-1} \Gamma(c')} F(a',b';c';\frac{1}{x^2})\biggr)
\tag{4.14}
$$
if $x>1$.
\endproclaim

Finally, for the last kernel we have
\proclaim{Proposition 4.7} Let $A_{11}$ be defined by (4.2); then for
all $x>0$
$$
\multline
(2\pi)^{2s-1} A_{11}(r,x;s.\nu)=\frac{ix^{1-2s}}{2\sinh \pi
r}\biggl(\cos \pi (s+ir) \frac{\Gamma(a) \Gamma(b)}{x^{c-1} \Gamma(c)}
F(a,b;c;-\frac{1}{x^2}) -\\
-\cos \pi (s-ir) \frac{\Gamma(a') \Gamma(b')}{x^{c'-1} \Gamma(c')}
F(a',b';c';-\frac{1}{x^2})\biggr)
\endmultline                              \tag{4.15}
$$
and this kernel equals to
$$
\multline
\frac{\cosh \pi r}{2 \cos \pi \nu}\biggl(\frac{\Gamma(a)
\Gamma(a')}{\Gamma(2\nu)} x^{2\nu} F(a,a';2\nu; -x^2) -\\
-\frac{\Gamma(b)
\Gamma(b')}{\Gamma(2-2\nu)} x^{2-2\nu}F(b,b';2-2\nu; -x^2)\biggr)
\endmultline
                                                   \tag{4.16}
$$
\endproclaim
\subhead
4.3.The new representations for $h_0$ and $h_1$
\endsubhead
In this subsection we express the integrals (2.19) and (2.20) in terms
of integrals with the kernels $A_{jk}$; it gives us the method to
estimate these coefficients.
\proclaim{Lemma 4.1} Let
$$
0\leqslant |\r{Im} r| <1/2\leqslant \r{Re}s, \r{Re}\rho <1,
\r{\nu}=\r{\mu}=1/2
$$
and let the function $h(r)$ be taken under all conditions of Theorem
2.3.

Then we have the following equalities for two integrals (2.19) and
(2.20):
$$
h_j=\frac{1}{2\pi}\int_{-\infty}^{\infty}B_j(r,u;s,\nu;\rho,\mu)\,h(u)\,d\chi(u)
-\sum_{l\in L}(-1)^l c(l) b_{j,l}(s,\nu;\rho,\mu)              \tag(4.17)
$$
where for $j=0$ and $j=1$
$$
B_j(r,u;s,\nu;\rho,\mu)=B_{j0}(r,u;s,\nu;\rho,\mu)+B_{j1}(r,u;s,\nu;\rho,\mu)
\tag(4.18)
$$
with
$$
B_{0j}=\int_0^{\infty}A_{0j}(r,\sqrt x;\rho,\nu)\,A_{0j}(u,\frac{1}{\sqrt
x}; 1-\rho,\mu)\,x^{\rho-s-1}\,dx,                          \tag{4.19}
$$
$$
B_{1j}=\int_0^{\infty}A_{1j}(r,\sqrt x;\rho,\nu)\, A_{0j}(u,\frac{1}{\sqrt
x};1-\rho,\mu)\,x^{\rho-s-1}\,dx                               \tag{4.20}
$$
and
$$
b_{jl}(r;s,\nu;\rho,\mu)=B_{j}(r,-(l-1/2)i;s,\nu;\rho,\mu).
                                                            \tag{4.21}
$$
Furthermore, there is the symmetry
$$
B_{j}(r,u;s,\nu;\rho,\mu)=B_{j}(r,u;\rho,\mu;s,\nu)
                                                         \tag{4.22}
$$
(the pair $(s,\nu)$ have been replaced by $(\rho,\mu)$ and otherwise).
\endproclaim

The proof is based on the known Mellin's transforms for the Bessel
functions:
$$
\int_0^{\infty}k_0(x,\nu)\,x^{2w-1}\,dx =\gamma(w,\nu)\cos \pi w,\,\,\,
|\r{Re}\nu-1/2|<\r{Re}w<1/2, \tag{4.23}
$$
$$
\int_0^{\infty}k_1(x,\nu)\,x^{2w-1}\,dx=\gamma(w,\nu) \sin \pi \nu,
\,\,\,|\r{Re}\nu-1/2|<\r{Re}w                            \tag{4.24}
$$
(the combinations of the Bessel functions $k_0$ and $k_1$ are defined by
(1.17) and (2.51), $\gamma(w,u)$ is given by (1.31)).

We substitute the definition of $\widehat \Phi$ and change the order of the
integration; it is possible since we assume the fast decreasing of $h$.

We will use the folloing fact from the Mellin theory. We say that
$f_1,f_2 \in L^2$ if
$$
\int_0^{\infty}|f_j|^2\,\frac{dx}{x}<\infty.             \tag{4.25}
$$

If for some $\Delta$ we have $x^{\Delta}f_1\in L^2, x^{\sigma-\Delta}f_2
\in L^2,$ then for any $s$ with the condition $\r{Re}s=\sigma$ we have
$$
\int_0^{\infty}f_1\,f_2 \,x^{s-1}\,dx=\frac{1}{2\pi
i}\int\limits_{(\Delta)} F_1(w)\,F_2(s-w)\,dw,                \tag{4.26}
$$
where the integration is doing over the line $\r{Re}\,w =\Delta$ and
$F_j$ are the Mellin transform $f_j$.

\noindent
{\bf 4.3.1. Function $B_{00}$}

Two functions
$$
f_1(y)=k_0(y,1/2+ir)y^{2c}, f_2(y)=k_0(xy,\nu)y^{2s-2c-1}   \tag{4.27}
$$
for any fixed $x>0$ satisfy to condition (4.25) if for $s$ with $\r{Re}
s>1/2$ the parameter $c$ be taken with condition
$$
\max (0,\sigma -3/4)<\r{Re}c<\min (1/4,\sigma-1/2), \sigma =\r{Re} s.
$$
Using now (4.23) we come to the equality
$$
\multline
\frac{(4\pi)^{2s-1}}{\pi \sqrt x}A_{00}(r,\sqrt x;s,\mu)=\\
=\frac{1}{2\pi
i}\int\limits_{(0)}\gamma(c+w/2,1/2+ir)\gamma(s-c-w/2,\mu)\cos \pi(c+w/2)
\times\\  \times\cos \pi(s-c-w/2) x^{-s+c+w/2}\,dw\\
=-\frac{1}{\pi
i}\int\limits_{(\Delta)}\gamma(1-w,1/2+ir)\gamma(s-1+w,\mu)\cos \pi w
\cos \pi(s+w)\times \\
\times x^{1-s-w}\,dw
\endmultline                                               \tag{4.28}
$$
(here we can take any $\Delta$ with $1>\Delta>1-\r{Re}s $).

This equality means that two functions
$$
2(4\pi)^{2s-2}\,x^{s-3/2}\,A_{00}(r,\sqrt x;s,\mu)
$$
and
$$
F_1(w)=-\gamma(1-w,1/2+ir) \gamma(s-1+w,\mu)\cos \pi w \cos \pi(s+w)
$$
are the Mellin pair.

By the same way we come to the second pair
$$
2(4\pi)^{-2s}\,x^{1/2-\rho}\,A_{00}(u,\frac{1}{\sqrt x}; 1-s,\nu)
$$
and
$$
F_2(w)=-\gamma(w+1-s-\rho,1/2+iu) \gamma(\rho-w,\nu)\cos \pi(\rho-w)
\cos \pi(w-s-\rho)
$$
(in (4.28) we change $\mu,r$ by $\nu,u$ and take $c=1-s-\rho$).

Using (4.26) again we come to the equality
$$
\multline
\frac{1}{2\pi i}\int\nolimits F_2(w)\,F_1(1-w)\,dw  =\\
=\frac{1}{2\pi i}\int\limits_{(\Delta)}\gamma(w,1/2+ir)\gamma(s-w,\mu)
\gamma(w+1-s-\rho,1/2+iu)\gamma(\rho-w,\nu)\times\\
\times \cos \pi w \cos \pi(s-w) \cos \pi(\rho-w) \cos
\pi(w-s-\rho)\,dw\\
=\frac{1}{4 \pi^2}\int_0^{\infty}A_{00}(r,\sqrt
x;s,\mu)A_{00}(u,\frac{1}{\sqrt x};1-s,\nu)\,x^{s-\rho-1}\,dx
\endmultline                                           \tag{4.29}
$$
(here $\max (0,\r{Re}(s+\rho)-1)<\Delta<\min (\r{Re}s, \r{Re}\rho$)).

If we take (4.28) with $\rho,\nu$ instead of $s,\mu$ and the second pair
with $\rho,\mu$ instead of $s,\nu$ then we receive the equality
$$
\multline
\frac{1}{2\pi
i}\int\limits_{(\Delta)}\gamma(w,1/2+ir)\gamma(\rho-w,\nu)\gamma(w+1-s-\rho,1/2+iu)\gamma(s-w,\mu)
\times\\
\times \cos \pi w \cos \pi(\rho-w) \cos \pi(s-w) \cos \pi(w+1-s-\rho)
\,dw=\\
=\frac{1}{4\pi^2}\int_0^{\infty}A_{00}(r,\sqrt
x;\rho,\nu)\,A_{00}(u,\frac{1}{\sqrt x};1-s,\nu)\,x^{\rho -s-1}\,dx
\endmultline                                              \tag{4.30}
$$

It gives us the equality (4.19) for $j=0$ and we see that
$$
B_{00}(r,u;s,\nu;\rho,\mu)=B_{00}(r,u;\rho,\mu;s,\nu).          \tag{4.31}
$$

{\bf 4.3.2. Function $B_{01}$}

Using the same formulas we have
$$
\multline
2(4\pi)^{2\rho-2}\int_0^{\infty}x^{\rho-3/2}A_{01}(r,\sqrt
x;\rho,\nu)\,x^{-w}\,dw=\\
=\gamma(w,1/2+ir)\gamma(\rho-w,\nu)\cos \pi w \sin \pi \nu
\endmultline                               \tag{4.32}
$$
and
$$
\multline
2(4\pi)^{-2\rho}\int_0^{\infty}x^{1/2-s}\,A_{01}(u,\frac{1}{\sqrt
x};1-\rho,\mu)\,x^{w-1}\,dw=\\
=-\gamma(s-w,\mu)\gamma(w+1-s-\rho,1/2+iu)\sin \pi \mu \cos \pi(w-s-\rho)
\endmultline                                          \tag{4.33}
$$

As the consequence we have the equality (4.19) for the case $j=1$; after
the simultaneous change the pair $(\rho,\nu)$ by $(s,\mu)$ in (4.32) and
(4.33) and taking into account (4.31) we come to the equality (4.22)
for the case $j=0$.

\noindent
{\bf 4.3.3. Functions $B_{10}$ and $B_{11}$}

First of all we have two equalities with $\gamma(w,1/2+ir)$:
$$
\multline
2(4\pi)^{2s-2}\int_0^{\infty}x^{s-3/2}A_{11}(r,\sqrt
x;s,\mu)\,x^{-w}\,dw=\\
=\gamma(w,1/2+ir)\gamma(s-w,\mu) \cosh \pi r\cos \pi(s-w) ,
\endmultline                                        \tag{4.34}
$$
$$
\multline
2(4\pi)^{2\rho-2}\int_0^{\infty}x^{\rho-3/2}A_{10}(\rho,\sqrt
x;\rho,\nu)\,x^{-w}\,dx=\gamma(w,1/2+ir)\gamma(\rho-w,\nu)\times\\
\times \cosh \pi r \sin \pi \nu
\endmultline                                         \tag{4.35}
$$

Furthermore,
$$
\multline
2(4\pi)^{-2s}\int_0^{\infty}x^{1/2-\rho}A_{01}(u,\frac{1}{\sqrt
x};1-s,\nu)\,x^{w-1}\,dx
=\\
=-\gamma(w-s-\rho+1,1/2+iu)\gamma(\rho-w,\nu)
 \sin \pi \nu \cos \pi(w-s-\rho),
\endmultline                                               \tag{4.36}
$$
$$
\multline
2(4\pi)^{-2\rho}\int_0^{\infty}x^{1/2-s}A_{00}(u,\frac{1}{\sqrt
x};1-\rho,\mu)\,x^{w-1}\,dx
=\\
=-\gamma(w-s-\rho+1,1/2+iu)
 \gamma(s-w,\mu) \cos \pi (s-w) \cos \pi(w-s-\rho).
\endmultline                                                \tag{4.37}
$$

It follows from these equalities that we have two representations
$$
\multline
i\int\limits_{(\Delta)}\gamma(w,1/2+ir)\gamma(\rho-w,\nu)\gamma(s-w,\mu)
\gamma(w-s-\rho+1)\times\\
\times \cosh \pi r \sin \pi \nu \cos \pi(s-\nu)\cos
\pi(w-s-\rho)\,dw=\\
=\frac{1}{2\pi}\int_0^{\infty}A_{11}(r,\sqrt
x;s,\mu)\,A_{01}(u,\frac{1}{\sqrt x};1-s,\nu)\,x^{s-\rho-1}\,dx\\
=\frac{1}{2\pi}\int_0^{\infty}A_{10}(r,\sqrt
x;\rho,\nu)A_{00}(u,\frac{1}{\sqrt x};1-\rho,\mu)\,x^{\rho-s-1} \,dx
\endmultline
\tag{4.38}
$$

Replacing here the pair $(s,\nu)$ by $(\rho,\mu)$ (and writing $(s,\nu)$
instead of $(\rho,\mu)$ we come to the equality
$$
B_{10}(r,u;s,\nu;\rho,\mu)=B_{11}(r,u;\rho,\mu;s,\nu)        \tag{4.39}
$$
and we have (4.22) for the case $j=1$ also.

All integrals with the Bessel function $J_{2l-1}$ (instead of $J_{2iu}$)
have been considered by the same way and we can omitt the corresponding
computations.

\head
\S 5. ASYMPTOTIC FORMULAS FOR THE KERNELS
\endhead

It is necessary now write the asymptotic expansions for the Gauss
hypergeometric functions in the integral representations for $h_j$
when $s=\mu=1/2,\,\nu=1/2+it$ and positive $t$ and $r$ are sufficintly
large.These formulas are given in subsections 5.1 -- 5.5.

Of course, these formulas must be taken from standard handbooks, but
nobody wrote these ones. The methods of the asymptotic integration
of the ordinary differential equations of the second order with a large
parameter are well known (see, for example, [9]). For our kernels it
rests determine the corresponding normalizing coefficients from the
boundary conditions to write the asymptotic expansion.
\subhead
5.1.The kernel $A_{01}(u,x;1/2,\nu)$
\endsubhead

Let $\nu=1/2+it$ and $t \to +\infty$. We define two asymptotic series
$$
\Cal A(\xi,u;t)=\sum_{n\geqslant 0}\frac{a_n (\xi,u)}{t^{2n}},\,\, \Cal B
(\xi,u;t) =\sum_{n\geqslant 1}\frac {b_n (\xi,u)}{t^{2n}}
$$
by the following recurrent relations: $a_0\equiv 1,\, b_0\equiv 0$ and for
$n\geqslant 0$
$$
b'_{n+1}=\frac{1}{2}(a''_n - f a_n) -(4u^2+1/4)\xi^{-1}(\xi^{-1} b_n)',
\tag{5.1}
$$
$$
a'_n =\frac{1}{2}(- b''_n +  f b_n),                       \tag{5.2}
$$
where $'=\frac{d}{d \xi}$ and
$$
f=u^2(\frac{4}{\xi^2}-\frac{1}{\sinh ^2
\xi/2})+\frac{1}{4}(\frac{1}{\xi^2}-\frac{1}{\sinh^2 \xi}).  \tag{5.4}
$$

To integrate (5.1) and (5.2) we take the initial conditions
$$
b_n(0)=0, a_n(0)=-(1/2+2iu) b'_n (0);                     \tag{5.5}
$$
note that all $a_n$ are even and all $b_n$ are odd functions in $\xi$ and
for $|\xi| < \pi$ they may be represented by the convergent power
series.

After determination $a_n,\,b_n$ from the recurrent relations we define
for an integer $M\geqslant 1$ the following segment of the asymptotic
series:
$$
V_M(\xi;u)=\sqrt \xi J_{2iu}(t\xi)\sum_{0\leqslant n \leqslant
M}\frac{a_n}{t^{2n}} +(\sqrt \xi J_{2iu}(t\xi))'\sum_{1\leqslant n
\leqslant M}\frac {b_n}{t^{2n}}                  \tag{5.6}
$$

With this notations we have the following uniform expansion.
\proclaim{Lemma 5.1} Let $t \to +\infty$; then for real $u, |u|\ll
t^{\varepsilon}$ for some small fixed $\varepsilon >0$ we have uniformly
in $\xi\geqslant0$ for any fixed integer $M\geqslant 1$
$$
\multline
\sqrt{\sinh (\xi)} A_{01}(u,\frac{1}{\sinh \xi/2};1/2,1/2+it)=\\
=\frac{i}{2\sinh \pi u}(\lambda(u,t) V_M(\xi;u) - \lambda(-u,t)
V_M(\xi,-u)+O(\min (\sqrt \xi,\frac{1}{\sqrt t})t^{-2M-2})),
\endmultline                                                 \tag{5.7}
$$
where
$$
\lambda(u,t)=t^{-2iu}\Gamma(1/2+it+iu)\Gamma(1/2-it +iu) \cosh \pi t .
                                                                 \tag{5.8}
$$
\endproclaim

This assertion is the consequence of two facts.

The first one follows from (4.4): function
$$V= \sqrt{\sinh \xi}
A_{01}(u,(\sinh \xi/2)^{-1};1/2,1/2+it)
$$
satisfies to the differential
equation
$$
V'' +(t^2+\frac{u^2}{\sinh^2 \xi/2}+\frac{1}{4\sinh^2 \xi}) V =0
\tag{5.9}
$$
or, the same,
$$
V''+(t^2 +\frac{4u^2}{\xi^2} +\frac{1}{4\xi^2})V = f V    \tag{5.10}
$$
with $f$ from (5.4).

The second fact follows from (4.10): we have
$$
\multline
A_{01}(u,\frac{1}{\sinh \xi/2};1/2,1/2+it)=\\
=\frac{i}{2\sinh \pi
u}\biggl(\frac{\lambda(u,t)}{\Gamma(1+2iu)}
(\frac{t\xi}{2})^{2iu}(1+O(\xi^2))-  \\
-\frac{\lambda(-u,t)}{\Gamma(1-2iu)}(\frac{t\xi}{2})^{-2iu}(1+O(\xi^2))\biggr)
\endmultline                                      \tag{5.11}
$$
as $\xi \to 0$.

The uniform nearness of the solutions of the equation (5.10) to the
corresponding Bessel functions is the well known fact.

So it rests check the coincidence of the asymptotic expansions on the
right sides when $\xi \to 0$. But the initial values $a_n$ and $b_n$
have been taken from this condition namely.

\noindent
{\bf 5.2. The kernels $A_{01}((l-1/2)i,x;1/2,\nu)$}

For all integers $l\geqslant 1$ we have
$$
\multline
A_{01}(l-1/2)i,\frac{1}{\sinh \xi/2};1/2,1/2+it) =\\
=(-1)^l \cosh \pi t \frac{\Gamma(l+it)\Gamma(l-it)}{\Gamma(2l)}(\sinh
\xi/2)^{2l-1} F(l+it,l-it;2l;-\sinh^2 \xi/2);
\endmultline                                            \tag{5.12}
$$
so when $\xi \to 0$ this function equals to
$$
(-1)^l
\frac{\lambda(-(l-1/2)i,t)}{\Gamma(2l)}(\frac{t\xi}{2})^{2l-1}(1+O(\xi^2))
                                                           \tag{5.13}
$$

It means we must take that solution of equation (5.10) which is near to
the function $const\times \sqrt \xi J_{2l-1}(t\xi)$. More exactly,
we have
\proclaim{Lemma 5.2} Let $t$ be sufficiently large and $a_n, b_n$ are
defined by the recurrent relations (5.1)--(5.5) with $u =(l-1/2)i$. Then
for any fixed integer $l\geqslant 1$ we have (uniformly in $\xi
\geqslant 0$ and for any fixed $M\geqslant 1$)
$$
\multline
\sqrt \sinh \xi A_{01}(i(l-1/2),\frac{1}{\sinh \xi/2};1/2,1/2 +it)=\\
=(-1)^l \lambda(-i(l-1/2),t)\biggl(\sqrt \xi
J_{2l-1}(t\xi)\sum_{0\leqslant n\leqslant M}\frac{a_n}{t^{2n}}\\
+(\sqrt \xi J_{2l-1}(t\xi))'\sum_{1\leqslant n \leqslant
M}\frac{b_n}{t^{2n}}+O(\min (\sqrt \xi (t\xi)^{2l-1},\frac{1}{\sqrt
t})t^{-2M-2})\biggr)
\endmultline                                           \tag{5.14}
$$
\endproclaim
\noindent
{\bf 5.3. The kernel $A_{11}(r,x;1/2,1/2)$}

We define the sequence $\{ c_n,d_n\}$ by the following recurrent
relations: $c_0\equiv 1, d_0\equiv 0$ and for $n\geqslant 0$
$$
d'_{n+1}=\frac{1}{2} (c_{n}-\frac{1}{2\xi}(\frac{d_n}{\xi})' - f_1\,
c_n),
                                                            \tag{5.15}
$$
$$
c'_n = \frac{1}{2}(- d''_n + f_1\, d_n),                 \tag{5.16}
$$
where $'=\frac {d}{d\xi}$ and
$$
f_1 =\frac{1}{4}(\frac{1}{\xi^2}-\frac{1}{\sinh^2 \xi})   \tag{5.17}
$$
We suppose for all $n\geqslant 1$
$$
d_n(0)=0                                                   \tag{5.18}
$$
(so $d_n(\xi)$ is an odd function in $\xi$ and $(\xi^{-1}d_n(\xi))'=
\frac{1}{3} d'''(0)+...$ as $\xi \to 0$; for this reason $
\xi^{-1}(\xi^{-1}d_n(\xi))'$ is finite in the neighbourhood of the point
$\xi = 0$).

We take the boundary values for $c_n$ in the form
$$
c_n(0)=-\frac{1}{2}d'_n (0),\,n\geqslant 1,         \tag{5.19}
$$
so we can write instead of (5.15) and (5.16):
$$
d_{n+1}(\xi) =\frac{1}{2} c'_n (\xi)-\frac{1}{4}\int_0^{\xi}(\frac{d_n
(\eta)}{\eta})'\frac{d\eta}{\eta}-\frac{1}{2}\int_0^{\xi}f_1(\eta)c_n
(\eta)\,d\eta, \,n\geqslant 0                         \tag{5.20}
$$
$$
c_n (\xi)=-\frac{1}{2} d'_n (\xi) +\frac{1}{2}\int_0^{\xi}f_1(\eta)d_n
(\eta)\,d\eta, n\geqslant 1.                           \tag{5.21}
$$
With these functions we define the asymptotic expansion
$$
V_{11}(r,\xi)=\sqrt \xi Y_0(r\xi)\sum_{n\geqslant 0}\frac {c_n}{r^{2n}}
+(\sqrt \xi Y_0(r\xi))'\sum_{n\geqslant 1}\frac{d_n}{r^{2n}} \tag{5.22}
$$
and let $V_{11,M}$ denotes sum of the first terms with $n\leqslant M$.

The function $(2\tanh \xi/2)^{-\frac{1}{2}}A_{11}(r,\sinh \xi/2;1/2,1/2)$
is the solution of the differential equation (again this fact is the
consequence of (4.4))
$$
V''+(r^2 +\frac{1}{4\sinh^2 \xi}) V =0           \tag{5.22}
$$
or, the same,
$$
V''+(r^2 + \frac{1}{4\xi^2})V =f_1(\xi)\,V.                      \tag{5.23}
$$

The function $V_{11}$ is the asymptotic solution of this equation (the
recurrent relations (5.20)--(5.21) are taken from this condition).

It is possible that these two solutions are proportional; but for our
purposes the more weak assertion will be sufficient.

We define the second asymptotic solution, supposing
$$
\widetilde V_{11}=\sqrt \xi J_0(r\xi)\sum_{n\geqslant 0}\frac{\tilde
c_n}{r^{2n}}+(\sqrt \xi J_0(r\xi))'\sum_{n\geqslant 1}\frac{\tilde
d_n}{r^{2n}},                                              \tag{5.25}
$$
where $\tilde c_0 \equiv 1$ and $\tilde c_n, \tilde d_n$ for $n\geqslant
1$ are defined by the same recurrent relations (5.15)--(5.16) but with
the following initial conditions.

We take for $n\geqslant 1$
$$
\tilde d_n(0) = 0.                                             \tag{5.26}
$$

Furthermore, let $B_n$ are the Bernoulli polynomials and
$$
\delta_n =d'_n(0) -\frac{(-1)^n}{2n}B_{2n}(1/2)                \tag{5.27}
$$
(the first even Bernoulli polynomials are: $B_0 = 1, B_2(x) =
x^2-x+\frac{1}{6}, B_4(x) = x^4 -2x^3 +x^2 -\frac{1}{30}, ...$ ).

We have $\delta_1 =\delta_2 = 0$ (it is the result of the direct
calculation; may be, the same is true for all $n$).

Let $N\geqslant 3$ be the first integer for which $\delta_N \ne 0$.

If there is such $N, N\geqslant 3$, we take the boundary values for
$\tilde c_n$ from the equation
$$
\tilde c_n(0) + \frac{1}{2} \tilde d'_n(0)=\frac{\delta_{n+N}}{\delta_N},
 n\geqslant 1.                        \tag{5.28}
$$
\proclaim{Lemma 5.3} Let $r $ be positive and sufficiently large. Then
for any fixed integer $M\geqslant 1$ we have uniformly in $\xi \geqslant
0$
$$
\multline
(2 \tanh \xi/2)^{-1/2} A_{11}(r,\sinh \xi/2; 1/2, 1/2) =
-\frac{\pi}{2}V_{11,M}+\frac{\delta_N}{r^{2N}} \widetilde V_{11,M-N} + \\
+O(\min (\sqrt \xi(|\log\frac{1}{r\xi}| +1), \frac{1}{\sqrt r})) r^{-2M-2}
\endmultline                                                \tag{5.29}
$$
where the second term is zero for $M\leqslant N$.
\endproclaim

Really, the unknown solution must be the linear combination of two
known asymptotic solutions; for some $C_0(r),C_1(r)$ we have
$$
(2 \tanh \xi/2)^{-1/2} A_{11}(r,\sinh \xi/2;1/2, 1/2)=C_0(r)W
+C_1(r) \widetilde W.                                    \tag(5.30)
$$

If $\xi \to 0$ we have from (4.16) (after taking the limiting value at
$\nu =1/2$)
$$
\multline
(2\tanh \xi/2)^{-1/2}A_{11}(r,\sinh \xi/2;1/2, 1/2)=-\sqrt \xi
\biggl(\log \frac{r\xi}{2}-\frac
{\Gamma'}{\Gamma}(1)+\\
+\frac{1}{2}\bigl(\frac{\Gamma'}{\Gamma}(1/2+ir)+\frac{\Gamma'}{\Gamma}(1/2-ir)-
 2\log r\bigr) +O(\xi^2 \log \frac{1}{\xi})\biggr)
\endmultline                                              \tag{5.31}
$$

Here the sum of logarithmic derivatives of gamma-function is the known
\newline asymptotic series
$$
\frac{1}{2}\biggl(\frac{\Gamma'}{\Gamma}(1/2+ir)+\frac{\Gamma'}{\Gamma}
(1/2-ir)- 2\log r\biggr) =-\frac {1}{24 r^2}+...+\frac{(-1)^n B_{2n}(1/2)}
{2n r^{2n}}+... .                                               \tag{5.32}
$$

On the other side we have for the case $\xi \to 0$ the asymptotic series
$$
\multline
W=\frac{2}{\pi}\biggl((\log
\frac{r\xi}{2}-\frac{\Gamma'}{\Gamma}(1))\sum_{n\geqslant
0}\frac{2c_n(0)+d'_n(0)}{2 r^{2n}}+\sum_{n\geqslant 1}\frac{d'_n
(0)}{r^{2n}}+\\
+O(\xi^2 \log \frac{1}{\xi})\biggr),
\endmultline                                                   \tag{5.33}
$$
$$
\widetilde W =\sqrt \xi \sum_{n\geqslant 0}\frac{2 \tilde c_n(0)+
\tilde d'_n(0)} {2r^{2n}} +O(\xi^2)                            \tag{5.34}
$$

To have the same asymptotic series in front of $(\log \frac {r\xi}{2}
-\frac {\Gamma'}{\Gamma}(1))$ we must take $C_0=-\frac{\pi}{2}$ in (5.30)
and define the initial value for $c_n(\xi)$ by the equality (5.19).
Then for $C_1(r)=O(r^{-6})$ we have the same coefficients in front of
$r^{-2}$ and $r^{-4}$. If $M=2$ the process is finished. For $M>2$ the
additional term with $\tilde V_{11}$ may be required (I assume that
$\delta_n =0$ for all $n\geqslant 1$ and there is no need to add this
term). In this case we take $C_1(r)=\delta_N\,r^{-2N}$ in (5.30)
(assuming $\delta_N \ne 0$) and after that we can determine $\tilde
c_n(0)$ by the equality (5.28).

As before, we omitt the proof of the known fact of nearness of the
solutions (5.23) to the corresponding Bessel functions.

\noindent
{\bf 5.4. The kernel $A_{10}(r,x;1/2, 1/2)$.}

For this kernel we will give three different asymptotic formulas (for the
intervals $[0,1-\delta), (\delta,1]$ and $[1,\infty)$).

We begin from the differential equations.The function
$$
w = (2 \tan \xi/2)^{1/2} A_{01}(r,\cos \xi/2;1/2, 1/2)
$$
satisfies to the equation
$$
w'' +(-r^2 +\frac{1}{4\sin^2 \xi})w = 0              \tag{5.35}
$$
and the function $\tilde w =(2\tanh \xi/2)^{1/2}A_{10}(r,\cosh
\xi/2;1/2, 1/2)$ is the solution of the equation
$$
\frac{d^2  \tilde w}{d\xi^2}+(r^2 +\frac{1}{4\sinh^2\xi})\tilde w = 0.
                                                       \tag{5.36}
$$

For the case $0\leqslant \xi \leqslant \pi-\delta$ with a fixed (small)
$\delta > 0$ the following simple formula would be sufficient.
\proclaim{Lemma 5.4} Let $r \to +\infty$; then for any fixed $\delta \in
(0,\pi/2)$ we have, uniformly in $0\leqslant \xi \leqslant \pi -\delta$,
$$
(2 \tan \xi/2)^{1/2}A_{10}(r,\cos \xi/2;1/2, 1/2)= \frac{\pi}{\cosh \pi
r}\sqrt \xi I_0(r\xi)(1+O(\frac{1}{r})).                     \tag{5.37}
$$
\endproclaim

The case $0<\delta\leqslant \xi \leqslant \pi$ is more complicated. As
before, we define the sequence $\{a_n ,b_n\}$ by the recurrent
relations: $a_0\equiv 1, b_0\equiv 0,$ and
$$
b'_{n+1}=\frac{1}{4\xi}(\frac{b_n}{\xi})' - \frac{1}{2}(a''_n
-f_2 a_n), n\geqslant 0,                 \tag{5.38}
$$
$$
a'_n =-\frac{1}{2}(b''_n -f_2 b_n), n\geqslant 1, \tag{5.39}
$$
where
$$
f_2 = \frac{1}{4}(\frac{1}{\xi^2} -\frac{1}{\sinh^2 \xi }).  \tag{5.40}
$$

We take the following initial coditions for these equations:
$$
b_n(0)=0, a_n(0)=-\frac{1}{2}b'_n(0)=\frac{(-1)^{n+1}}{4 n}B_{2n}(1/2).
                                                             \tag{5.41}
$$
\proclaim{Lemma 5.5} Let the positive $r$ be sufficiently large; then
for any fixed $M\geqslant 1$ we have for $0\leqslant \xi\leqslant
\pi-\delta$ with any fixed $\delta \in (0,\pi/2)$
$$
\multline
(2\tan \xi/2)^{-1/2}A_{10}(r,\sin \xi/2;1/2, 1/2)=\sqrt \xi
K_0(r\xi)\sum_{0\leqslant n \leqslant M}\frac{a_n}{r^{2n}}+\\
+(\sqrt \xi K_0(r\xi))'\sum_{1\leqslant n \leqslant M}\frac{b_n}{r^{2n}}
+O(\min (\sqrt \xi (|\log \frac{1}{r\xi}|)+1),\frac{1}{\sqrt r}e^{-r\xi})
r^{-2M-2}
\endmultline                                           \tag{5.42}
$$
\endproclaim

Really, the function on the left side must be near to the combination
$$
\sqrt \xi K_0(r\xi)(1+O(\frac{1}{r}))+C(r)\sqrt \xi
I_0(r\xi)(1+O(\frac{1}{r}))                            \tag{5.43}
$$
since for $\xi \to 0$ we have (it is the limiting case (4.13) as $\nu
\to 1/2$)
$$
\multline
(2\tan \xi/2)^{-1/2}A_{10}(r,\sin \xi/2;1/2, 1/2)=-\sqrt \xi \biggl(\log
\frac{r\xi}{2} -
\frac{\Gamma'}{\Gamma}(1)+\\
+\frac{1}{2}\bigl(\frac{\Gamma'}{\Gamma}(1/2+ir)+
\frac{\Gamma'}{\Gamma}(1/2-ir) - 2\log r\bigr)+
O(\xi^2 \log \frac{1}{\xi})\biggr)
\endmultline                              \tag{5.44}
$$
and
$$
K_0(r\xi)= -\log \frac{r\xi}{2}+\frac{\Gamma'}{\Gamma}(1)+O(\xi^2 \log
\frac{1}{\xi}).                              \tag{5.45}
$$

Taking in (5.43) $\xi=\pi-\delta$ and $\xi=\delta$ in (5.34) with a
fixed small $\delta$ we come to the equality
$$
O(e^{-(\pi-\delta)r})=O(e^{-(\pi-\delta)r})
+C(r)e^{(\pi-\delta)r}(1+O(\frac{1}{r})),          \tag{5.46}
$$
so for any fixed $\delta>0$ we have
$$
|C(r)| \ll e^{-2(\pi-\delta)r}.                        \tag{5.47}
$$

It means the second term in (5.43) is exponentially small; now the
comparison of the asymptotic series at $\xi=0$ gives us (5.41).

By the way, we come to the new (and very unexpected) method for the
calculation of the Bernoulli numbers. Since for $n\geqslant 1$
$$
B_{2n}(0)=B_{2n}=\frac{2^{2n-1}}{2^{2n-1}-1}B_{2n}(1/2)  \tag{5.48}
$$
we have for $n\geqslant 1$
$$
B_{2n}=(-1)^n \frac{2^{2n}}{2^{2n-1}-1}\,n\,b'_n(0),    \tag{5.49}
$$
where $b_n$ are defined by the recurrent relations (5.38)--(5.39).

Finally, we consider our kernel when argument is larger than 1. Now we
have that solution of (5.36) which is
$$
\frac{\pi}{\cosh \pi r} \sqrt \xi (1+O(\xi^2))
$$
when $\xi \to 0$. This boundary condition gives us, uniformly in $\xi
\geqslant 0$,
$$
\multline
(2 \tanh \xi/2)^{1/2}A_{10}(r,\cosh \xi/2;1/2, 1/2)=\frac{\pi}{\cosh \pi
r}(\sqrt \xi J_0(r\xi) + \\
+\frac{1}{r}O(\min(\sqrt \xi, \frac{1}{\sqrt r}))).
\endmultline                                         \tag{5.50}
$$

\noindent
{\bf 5.5. The kernel $A_{00}(u,x;1/2,\nu)$}

First of all, the function
$$
v=\sqrt{\sinh \xi} A_{00}(u,\frac {1}{\cosh \xi/2};1/2, \nu)
$$
satisfies to the equation
$$
v''+(t^2 +\frac{1}{4\xi^2})\,v=f_3\,v,                    \tag{5.51}
$$
where
$$
f_3 =\frac{1}{4}(\frac{1}{\xi^2} - \frac{1}{\sinh^2 \xi})
+\frac{u^2}{\cosh^2 \xi/2}.                             \tag{5.52}
$$

Furthermore, taking the limiting case $s=1/2$ in (4.9) we have as $\xi
\to 0$:
$$
\multline
v=-2\sqrt \xi \biggl(\log \frac{t\xi}{2} -\frac{\Gamma'}{\Gamma}(1)+\\
+\frac{1}{4}(\frac{\Gamma'}{\Gamma}(\nu+iu)+\frac{\Gamma'}{\Gamma}(\nu
-iu) +\frac{\Gamma'}{\Gamma}(1-\nu +iu)+\frac{\Gamma'}{\Gamma}(1-\nu
-iu)-4\log t)+O(\xi^2 \log \frac{1}{\xi})\biggr)\\
=-2\sqrt \xi \biggl(\log \frac{t\xi}{2}-\frac{\Gamma'}{\Gamma}(1)+\\
+\sum_{n\geqslant
1}(-1)^{n-1}\frac{B_{2n}(1/2+iu)+B_{2n}(1/2-iu)}{4n\,t^{2n}}+O(\xi^2 \log
\frac{1}{\xi})\biggr)
\endmultline                                                 \tag{5.53}
$$
(here $B_{2n}$ are the Bernoulli polynomials again).

The equation (5.51) has the same type what we had considered in
subsection 5.3. As before, we define the coefficients $p_n,q_n$ by the
relations $p_0\equiv 1, q_0\equiv 0$ and
$$
q'_{n+1}=\frac{1}{2}(p''_n
-\frac{1}{2\xi}(\frac{q_n}{\xi})'- f_3 \,p_n,\,\, n\geqslant 0,
\tag{5.54}
$$
$$
p'_n = -\frac{1}{2}(q''_n - f_3\,q_n),\,\, n\geqslant 1.  \tag{5.55}
$$

The initial values for these relations we take in the form
$$
q_n(0)=0, p_n(0) =-\frac{1}{2}q'_n(0), n \geqslant 1.           \tag{5.56}
$$

Now we have
$$
q'_1 (0)=\frac{1}{2}\,f_3(0) = \frac{1}{4}(B_2(1/2 +iu)+B_2(1/2-iu)),
                                                               \tag{5.57}
$$
$$
q'_2(0)=\frac{1}{6}\,f''_3(0)-\frac{1}{4}\,{f_3}^2 (0)=-\frac{1}{8}(B_4
(1/2 +iu)+B_4 (1/2-iu)).        \tag{5.58}
$$

This coincidence (may be, it will be true for $n\geqslant 3$ also)
allows us, as in 5.3, write the following asymptotic formula.
\proclaim{Lemma 5.5} Let $t \to +\infty$ and let $u$ be real and for any
fixed $\varepsilon>0\,\, |u|\ll t^{\varepsilon}$; then, uniformly in $\xi
\geqslant 0$, we have
$$
\multline
\sqrt \xi A_{00}(u,\frac{1}{\cosh \xi/2};1/2,\,\nu) =\frac{1}{\pi}\sqrt
\xi\, Y_0 (t\xi)(1+\frac{p_1(\xi)}{t^2}+\frac{p_2(\xi)}{t^4})+\\
+\frac{1}{\pi}(\sqrt \xi\,
Y_0 (t\xi))'(\frac{q_1(\xi)}{t^2}+\frac{q_2(\xi)}{t^4}) +O(t^{-6}) \min
(\sqrt \xi (|\log t\xi|+1),\frac{1}{\sqrt t})
\endmultline                                      \tag{5.59}
$$
\endproclaim

Now we consider the same kernel for the case when argument is larger
than 1.

The function $v=\sqrt \sin \xi A_{00}(u,\frac{1}{\sin \xi/2};1/2,\,\nu)$
(and the same function with $\xi$ replaced by $\pi-\xi$) satisfies to
the differential equation
$$
v''+(-t^2 +\frac{u^2}{\sin^2 \xi/2}+\frac{1}{4 \sin^2 \xi})v =0;
\tag{5.60}
$$
if $\xi \to 0$ we find from (4.8)
$$
\multline
v=\frac{\Gamma (\nu+iu) \Gamma (1-\nu +iu)}{2\Gamma(1+2iu)}\sqrt \xi\,
(\xi/2)^{2iu}(1 +O(\xi^2))+\\
+\{\hbox{the same with $(-u)$}\}
\endmultline                                        \tag{5.61}
$$

It is sufficient to assert that we have
\proclaim {Lemma 5.6} Let $t \to +\infty$ and $u$ be real with condition
$|u|\ll t^{\varepsilon}$ for any fixed $\varepsilon >0$. Then we have
for $0 \leqslant \xi \leqslant \pi-\delta$ with any fixed $\xi \in
(0,\frac{\pi}{2})$
$$
\multline
\sqrt \sin \xi\,A_{00}(u,\frac{1}{\sin \xi/2};1/2,\,\nu)=\\
=\r{Re}\{t^{-2iu}\Gamma(\nu+iu)\Gamma(1-\nu+iu)\sqrt
\xi\,I_{2iu}(t\xi)(1+O(\frac{1}{t}))\}
\endmultline                                            \tag{{5.62}}
$$
\endproclaim
Finally, in the neighbourhood $\xi = \pi$ we have other asymptotic
expansion.
\proclaim {Lemma 5.7} Let $\{\tilde p_n,\tilde q_n \}$ are defined by the
relations: $\tilde p_0 \equiv 1, \tilde q_0\equiv 0,$ and
$$
\tilde q'_{n+1}=\frac{1}{2}(\frac{1}{2\xi}(\frac{\tilde
q_n}{\xi})'- \tilde p''_n + f_4 p_n),\,\, n\geqslant 0,
                                                              \tag{5.63}
$$
$$\tilde p'_n =-\frac{1}{2} (\tilde q''_n - f_4 \tilde q_n),\,\,
n\geqslant 1,                                                \tag{5.64}
$$
where
$$
f_4 =\frac{1}{4}(\frac{1}{\xi^2}-\frac{1}{\sin^2 \xi})-\frac{u^2}{\cos^2
\xi/2};                                               \tag{5.65}
$$
let the initial values are taken as
$$
\tilde q_n(0)=0, \tilde p_n(0)=-\frac{1}{2}\tilde
q'_n(0)=\frac{(-1)^n}{4n}(B_{2n}(1/2 +iu)+B_{2n}(1/2-iu)).      \tag{5.66}
$$
Then for any fixed $\delta \in (0,\pi/2)$ we have for $0\leqslant \xi
\leqslant \pi-\delta$ the following asymptotic series
$$
\multline
\sqrt \sin \xi A_{00}(u,\frac{1}{\cos \xi/2};1/2,\nu)=\\
=2\sqrt \xi\,K_0(t\xi)\sum_{n\geqslant 0}\frac{\tilde p_n}{t^{2n}}+2(\sqrt
\xi\,K_0(t\xi))'\sum_{n\geqslant 1}\frac{\tilde q_n}{t^{2n}}.
\endmultline                                          \tag{5.67}
$$
\endproclaim

To receive this expansion we use (5.53) and repeat the considerations of
subsection 5.4.

\noindent
{\bf 5.6. The kernel $A_{10}(u,x;1/2,\nu)$}

Using (4.4) again we see: the function
$$
v = (2 \tanh \xi/2)^{1/2}(\cosh \xi/2)^{1-2\nu}A_{10}(u,\frac{1}{\cosh
\xi/2};1/2, \,\nu)                                           \tag{5.68}
$$
is that solution of the differential equation
$$
v''+(t^2 +\frac{1}{4 \sinh^2 \xi}-\frac{u^2}{\cosh^2 \xi/2})v = 0,
\tag{5.69}
$$
which equals to (we use (4.12) with $s=1/2, \nu=1/2 +it, r=u$)
$$
2\pi \frac{\cosh \pi t \,\cosh \pi u}{\cosh 2\pi t +\cosh 2\pi u} \sqrt
\xi\, (1+O(\xi^2))                                  \tag{5.70}
$$
when $\xi \to 0$.

We have the same equation what we had in 5.5, so we write the full
uniform asymptotic expansion with the main term
$$
4\pi\,e^{-\pi t} \cosh \pi u \,\sqrt \xi\, J_0(t\xi),\,\,\, \xi\geqslant 0.
                                                     \tag{5.71}
$$

For the case when $x>1$ (but $x$ is not very large) we have the exponentially
small function again (with the modified Bessel function instead of
$J_0$).

The function
$$
\tilde v = (2\tan \xi/2)^{1/2} (\cos
\xi/2)^{1-2\nu}A_{10}(u,\frac{1}{\cos \xi/2};1/2,\,\nu)  \tag{5.72}
$$
satisfies to the equation
$$
\tilde v'' +(-t^2 + \frac{u^2}{\cos^2 \xi/2} +\frac{1}{4 \sin^2
\xi})\tilde v= 0.                                 \tag{5.73}
$$

Together with initial conditions at $\xi =0$ it gives us the similar
asymptotic formula for $0\leqslant \xi \leqslant \pi-\delta$
$$
(2\tan \xi/2)^{1/2}(\cos \xi/2)^{1-2\nu}A_{10}(u,\frac{1}{\cos \xi/2};
1/2, \nu) \approx 4\pi\,e^{-\pi t}\,\cosh \pi u \sqrt \xi\,\, I_0(t\xi)
                                                     \tag{5.74}
$$
for any fixed $\delta \in (0,\pi/2)$.

Finally, if $x$ be sufficiently large (it means $0<\delta \leqslant \xi
\leqslant \pi$ in(5.72)) then the solution must be a combination of two
modified Bessel's functions.

The unique possibility which have been agreed with (5.74) is the
Mcdonald function; so we have for $0\leqslant \xi \leqslant \pi-\delta$
with any fixed $\delta \in (0,\pi/2)$
$$
(2\tan \xi/2)^{1/2}(\sin \xi/2)^{1-2\nu}A_{10}(u,\frac{1}{\sin \xi/2};
1/2,\,\nu) \approx 4 \cosh \pi u \sqrt \xi\,\,K_{2iu}(t\xi).   \tag{5.76}
$$

\noindent
{\bf 5.7. The kernel $A_{11}(u,x;1/2,\,\nu)$}

For this case we have the differential equation
$$
w''+(t^2 +\frac{u^2}{\sinh^2 \xi/2} +\frac{1}{4\sinh^2 \xi})w = 0
$$
for the function
$$
w = \sqrt{\sinh \xi} A_{11}(u,\frac{1}{\sinh \xi/2};1/2,\,\nu).
\tag{5.77}
$$

If $\xi \to 0$ we see (it follows from (4.15))
$$
w =\frac{\pi e^{-\pi t}}{2\sinh \pi u}\biggl(\bigl(\frac{t\xi}{2}\bigr)^{2iu}
\frac{(1+O(\xi^2))}{\Gamma(1+2iu)} -\hbox{(the same with $-u$)}\biggr).\tag{5.78}
$$

So this kernel is exponentially small and we have for the main term of
the uniform asymptotic expansion:
$$
\sqrt \sinh \xi A_{11}(u,\frac{1}{\sinh \xi/2}; 1/2, \nu)\approx
\frac{\pi e^{-\pi t}}{2 \sinh \pi u}
\sqrt \xi \,\bigl(J_{2iu}(t\xi)-J_{-2iu}(t\xi)\bigr).   \tag{5.79}
$$
\noindent
{\bf 5.8. The kernel $A_{01}(r,x;1/2,\,1/2)$}

For the function
$$
v=(2\tanh \xi/2)^{-1/2}A_{01}(r,\sinh \xi/2;1/2,\,1/2)    \tag{5.80}
$$
we have the differential equation
$$
v'' +(r^2 + \frac {1}{4\sinh \xi^2})v = 0                    \tag{5.81}
$$
with the initial condition (it follows from (4.11))
$$
v=\frac {\pi}{2\cosh \pi r}\sqrt \xi (1+ O(\xi^2)), \xi \to 0.
$$

It means this function is exponentially small for $r$ large; we have the
following main term for the uniform expansion
$$
(2 \tanh \xi/2)^{-1/2} A_{01}(r,\sinh \xi/2;1/2,\,1/2) \approx
\frac{\pi}{2\cosh \pi r}\sqrt \xi\,\,J_0(r\xi), \xi\geqslant 0.
\tag{5.82}
$$

\noindent
{\bf 5.9. The kernel $A_{00}(r,x;1/2,\,1/2)$}

For this case we have three different expansions.

Firstly we consider the case $0\leqslant \xi <1$. For the function
$$
v = (2 \tan \xi/2)^{-1/2} A_{00}(r,\sin \xi/2;1/2,\,1/2)     \tag{5.83}
$$
we have the differential equation
$$
v'' + (-r^2 + \frac {1}{4\sin^2 \xi}) v =0               \tag{5.84}
$$
and the initial condition (see (4.5) with $s=1/2,\nu=1/2$)
$$
v = \frac{\pi}{2\cosh \pi r}\sqrt \xi (1 + O(\xi^2)), \xi \to 0.
\tag{5.85}
$$

It means this function is exponentially small for $0\leqslant \xi
\leqslant \pi - \delta$ with any fixed $\delta \in (0,\pi/2)$:
$$
(2\tan \xi/2)^{-1/2}A_{00}(r, \sin \xi/2;1/2,\,1/2) \approx \frac{\pi}{2
\cosh \pi r}\sqrt \xi \,\,I_0(r\xi)(1+O(\frac{1}{r})).  \tag{5.86}
$$

The same equation (5.84) we have for the function
$$
\tilde v =(2 \tan \xi/2)^{1/2}A_{00}(r,\cos \xi/2; 1/2,\,\,1/2);
\tag{5.87}
$$
if $\xi \to 0$ we have from (4.9) (as the limiting case $s=1/2$)
$$
\multline
\tilde v =2\sqrt \xi \biggl(-\log \frac{r\xi}{2}
+\frac{\Gamma'}{\Gamma}(1)-\\
-\frac{1}{2}(\frac{\Gamma'}{\Gamma}(1/2 +ir) +\frac{\Gamma'}{\Gamma}(1/2
-ir)-2 \log r)\biggr)+ O(\xi^{5/2}\log\frac{1}{\xi}).
\endmultline                 \tag{5.88}
$$

Furthermore, if $\xi \geqslant \xi_0$ this function must be $O(e^{-\xi_0
r})$. So we have for any $\delta \in (0,\pi/2)$
$$
(2 \tan\xi/2)^{1/2} A_{00}(r,\cos \xi/2;1/2,\,1/2)\approx 2\sqrt
\xi\,\,K_0(r\xi),\,\,0\leqslant \xi \leqslant \pi-\delta. \tag{5.89}
$$

Finally, we have the equation
$$
w'' + (r^2 + \frac{1}{4\sinh^2 \xi})w =0            \tag{5.90}
$$
for the function
$$
w = (2\tanh \xi/2)^{1/2}A_{00}(r,\cosh \xi/2; 1/2,\,1/2).
$$

If $\xi \to 0$ we have the same expansion (5.88) for $w(\xi)$; for this
reason
$$
(2\tanh \xi/2)^{1/2}A_{00}(r,\cosh \xi/2;1/2,\,1/2)\approx -\pi \sqrt
\xi\,\,Y_0(r\xi).                           \tag{5.92}
$$

Of course, it is possible to replace the right side by the more detailed
asymptotic expansion of the form
$$
-\pi \{\sqrt \xi\,\,Y_0(r\xi)\sum_{n\geqslant
0}\frac{a_n}{r^{2n}}+(\sqrt \xi \,\,Y_0(r\xi))'\sum_{n\geqslant
1}\frac{b_n}{r^{2n}}\},               \tag{5.93}
$$
where $a_0\equiv 1, b_0\equiv 0 $ and the recurrent relations for these
coefficients may be written from the condition that this function is the
formal solution of (5.90).

\head
\S 6. THE AVERAGING ON THE RIGHT SIDE
\endhead

In this section we estimate the integrals
$$
\multline
M_{j,k}(T)=\frac {1}{i}\int\limits_{(1/2)}\Cal H_j(\rho+it)
\Cal H_j(\rho-it)\times\\
\times h_k(\varkappa_{j};1/2,\nu;\rho,1/2)\omega_T(\rho)d\rho,\ \ k=0,1.
\endmultline
                                                                   \tag 6.1
$$

The analogous integral for the case of the discret spectrum is
$$
\frac{1}{i}\int\limits_{(1/2)}\Cal
Z(\rho;\nu,1/2+ir)(h_0+h_1)(r;1/2,\nu;\rho,1/2) \omega_T (\rho)\,d\rho.
$$

For some cases the line of integration $\r{Re} \rho=2$ will be taken
instead of the initial line $\r{Re} \rho=1/2$; on the line $\r{Re} \rho
=2$ the result of integration will be expressed in the explicite form.

The contribution of four poles of the function $\Cal Z(\rho;\nu,1/2+ir)$
will be considered separatly in subsection 6.7; for all others there is
no difference in the consideration of the discret and continuous
spectrum.

\subhead
6.1.The integrals $M_{j,0};\varkappa_j\geqslant 2T_0$
\endsubhead

It is sufficient consider the case $\varkappa\leqslant 2T_0$, since for
large $r$ we have the following inequality.
\proclaim{Proposition 6.1}
Let $r\geqslant 2T_0$; then
$$
\vert h_0(r;1/2,\nu;\rho,1/2)\vert\ll\frac {1}{r^6} .    \tag 6.2
$$
\endproclaim

To receive this bound it will be sufficient estimate the absolute value
of the integrand in (2.19).

For the case $s=\mu=1/2$, $\nu=1/2+it$, $\rho=1/2+i\tau$ (we assume
$t=o(\tau^{1/4})$) the integrand in (2.19) is not larger than (we write
$w=1/2+i\eta,\ \eta\gg 1$)
$$
\multline
\exp(-\frac{\pi}{2}(\eta+r+\vert\eta-r
\vert+\vert\eta-\tau+t\vert+\vert\eta-\tau-t\vert)   \\
\text{max}(e^{\pi\vert\tau-2\eta\vert},e^{\pi\tau})
\vert\Hat\varPhi(1-2i\tau+2i\eta)\vert.
\endmultline                                                         \tag 6.3
$$

Here we can assume that $\vert\tau-T_0\vert\ll T\log T_0$ (because of
$\omega_T(\rho)$ in (6.1) is exponentially smoll for
$\vert\rho-iT_0\vert\gg T$).

So for $r\geqslant 2T_0$ the integrand in (2.19) is $O(\exp(-\pi(r-\eta)))$
if $\eta<r$. For this reason we can integrate over the line
$\eta\geqslant r-2\log r$.

But for $\eta\geqslant(1-\delta)r$, $r\geqslant 2T_0$, for any
fixed small $\delta>0$ we have
$$
\vert\Hat\varPhi(1-2i\tau+2i\eta)\vert\ll\eta^{-5},\ \
\vert\gamma(\rho-w,\nu)\vert\ll\eta^{-1},\ \
\vert\gamma(1/2-w,1/2)\vert\ll\eta^{-1}                            \tag 6.4
$$
and we come to the bound (6.2).

It follows from (6.2) that the contribution of the terms with
$\varkappa_j \geqslant 2T_0$ to the average on the right side (2.57) is $o(1)$
for the case of our specialization. Really, from the fuctional equation
one can see that $|\Cal H_j (\rho \pm it)| \ll
\varkappa_j^{1+\varepsilon}$  for the case $\varkappa_j \gg \tau
=\r{Im}\,\rho$. So the result of the averaging of $j$-th term is
$O(T\varkappa_j^{2+2\varepsilon})$. After that we have the sum with
terms $O(T\varkappa_j^{-4+2\varepsilon})\alpha_j \Cal H_j^2(1/2)$ and
$\varkappa_j\geqslant 2T_0$; this sum is $O(T^{-1+2\varepsilon})$.

\subhead
6.2.The terms with $B_{00}$ for $r\ll T_0$
\endsubhead

For this case the new representation (4.17) have been used.

Firstly we consider the integral $B_{00}$ and it would be convenient to
estimate the function $B_{00}(r,u;\rho,1/2;1/2,\nu)$ (we use the
symmetry (4.22) for this case).

\proclaim{Lemma 6.1}
Let $T_0,T,t$ are sufficiently large, $T=T_0^{1-\epsilon},\, t=o(T^{1/4})$;
then we have for $r\ll T_0$
$$
\multline
\Big|\int\limits_{(1/2)}B_{00}(r,u;\rho,1/2;1/2,\nu)\Cal H_j(\rho+it)
\Cal H_j(\rho-it)\omega_T(\rho)\,d\rho\Big|\ll                         \\
\ll\cases
\log^2T,\,\,\,\,\,\, \ r\ll\sqrt{T_0}                               \\
T\,r^{-2}\log^2r,\,\,\, r\gg\sqrt{T_0}
\endcases
\endmultline                                                         \tag 6.5
$$
\endproclaim

For the beginning we write in (4.17) (with the variables
$(\rho,1/2;1/2,\nu)$)
$$
B_{00}=\int_0^{1/\sqrt2}+\int_{1/\sqrt2}^1+\int_1^{\infty}=
B^{(1)}+B^{(2)}+B^{(3)},
$$
where (after the understandable change of the variable)
$$
B^{(1)}=\int_0^{\pi/2}f^{(1)}(\xi;r,u,t)
(\sin \xi/2)^{-2\rho}\,d\xi,                                         \tag 6.6
$$
$$
B^{(2)}=\frac12\int_0^{\pi/2}f^{(2)}(\xi;r,u,t)
(\cos \xi/2)^{-2\rho}\,d\xi,                                         \tag 6.7
$$
$$
B^{(3)}=\frac12\int_0^{\pi/2}f^{(3)}(\xi;r,u,t)
(\text{ch} \xi/2)^{-2\rho}\,d\xi.                                    \tag 6.8
$$

In this equalities
$$
f^{(1)}=(2\text{tg}\xi/2)^{-1/2}A_{00}(r,\sin\xi/2;1/2,1/2)
\sqrt{\sin\xi}A_{00}(u,\frac 1{\sin\xi/2};1/2,\nu),                  \tag 6.9
$$
$$
f^{(2)}=(2\text{tg}\xi/2)^{1/2}A_{00}(r,\cos\xi/2;1/2,1/2)
\sqrt{\sin\xi}A_{00}(u,\frac 1{\cos\xi/2};1/2,\nu),                  \tag 6.10
$$
$$
f^{(3)}=(2\text{th}\xi/2)^{1/2}A_{00}(r,\text{ch}\xi/2;1/2,1/2)
\sqrt{\text{sh}\xi}A_{00}(u,\frac 1{\text{ch}\xi/2};1/2,\nu);        \tag 6.11
$$
it is essential that these functions are not depending in $\rho$.

Using (5.86), (5.89), (5.62) and (5.67) we see that for any fixed
$\delta>0$ we have
$$
B^{(1)}+B^{(2)}=\frac12\int\limits_0^{\delta}f^{(2)}(\xi;r,u,t)
(\cos\xi/2)^{-2\rho}\,d\xi+O(\exp(-\delta(r+t)))                    \tag 6.12
$$
(in reality $B^{(1)}$ is $O(e^{-\frac{\pi}{4}(r+t)})$).

The remainder term may be omitted (it gives $o(1)$ into the final result).
For the remained integrals we integrate firstly over the variable $\rho$.
Using the Ramanujan integral again we have
$$
\multline
\frac {1}{i}\int\limits_{(1/2)}B_{00}(r,u;\rho,1/2;1/2,\nu)
\Cal H_j(\rho+it)\Cal H_j(\rho-it)\omega_T(\rho)\,d\rho=\\
=\frac{1}{i}\int\limits_{(2)}(...)\,d\rho =\\
=T\,\sum_{n,m\geqslant
1}\frac{t_j(n)t_j(m)}{(nm)^{1/2+iT_0}}(\frac{n}{m})^{it} \times\\
\times \Big\{\int\limits_0^{\delta}f^{(2)}(\xi;r,u,t)
\psi(nm\cos^2\xi/2)(cos \xi/2)^{-1-2iT_0}\,d\xi+\\
+\int\limits_0^{\infty}f^{(3)}(\xi;r,u,t)\psi(nm\cosh^2\xi/2)
(\cosh \xi/2)^{-1-2iT_0}
\,d\xi\Big\}+O(T_0e^{-\delta(r+t)}),
\endmultline                                            \tag 6.13
$$
where
$$
\psi (x)=\left(x^T+x^{-T}\right)^{-1}.
                                                        \tag 6.14
$$

In (6.13) we have $x=\cos^2\xi/2$ with small $\xi$ or
$x=\cosh^2\xi/2\geqslant 1$. For both cases we have the inequality
$x\geqslant 2-\delta $ with small $\delta$ if $nm\geqslant 2$.
So for $nm\geqslant 2$
$$
|\psi(nm x)|\ll a ^{-T}, 3/2< a < 2-\delta                 \tag 6.15
$$
and all these terms give $o(1)$ in the final result.

It means we can omit all terms with $nm\geqslant 2$ on the right side
(6.13), and take only one term with $n=m=1$ from this sum.

Now we come to the integrals with the main terms (we use (5.59), (5.67),
(5.89) and (5.92))
$$
T\int_0^{\delta}\xi K_0(r\xi)K_0(t\xi)
\psi(\cos^2 \xi/2)(\cos \xi/2)^{-1-2iT_0} \,d\xi,        \tag 6.16
$$
$$
T\int_0^{\infty}\xi Y_0(r\xi)Y_0(t\xi)
\psi(\cosh^2 \xi/2) (\cosh\xi/2)^{-1-2iT_0}\,d\xi.        \tag 6.17
$$

It is sufficient estimate these integrals; the next terms from the
corresponding asymptotic expansions for the kernels have the same nature
but with additional multipliers $r^{-2m}$ or $t^{-2n}$ with
$n,m\geqslant 1$.

In both integrals we can integrate in the interval
$0\leqslant\xi\leqslant\delta(T)$ with $\delta(T)=T^{-1/2}\log T$ only
(since the part of our integrals with $\xi\geqslant\delta(T)$
contributes $O(T\exp(-\frac14\log^2T))$ and this part may be omitted).

We write for $x>1$
$$
\psi(x)=x^{-T}\Big\{\sum_{0\leqslant m\leqslant M-1}
x^{-m}+x^{-M}(1+x^{-2T})^{-1}\Big\}
$$
and for $x<1$
$$
\psi(x)=x^{T}\Big\{\sum_{0\leqslant m\leqslant M-1}
x^m+x^M(1+x^{2T})^{-1}\Big\}.
$$

Let $\mu_m=1+2mT+2iT_0$; we write the power series
$$
(\cos\xi/2)^{\mu_m}=e^{-\mu_m\xi^2/8}
\Big(1+\sum_{k\geqslant 2}a_k(\mu_m)\xi^{2k}\Big).                   \tag 6.18
$$
Here we can assume $\xi\leqslant(mT)^{-1/2}\log T$ and for this reason we
have
$$
|a_{2k}(\mu_m)\xi^{4k}|\ll T_0^k(mT)^{-2k}(\log T)^{4k}\ll
T_0^{-k(1-2\varepsilon)+\varepsilon}                                $$
(note that $T=T_0^{1-\varepsilon}$); it means we can omit all terms with
$k\geqslant 4$ in (6.18).

Furthermore, we write the power series for $K_0(t\xi)$ in (6.16); since
$t\xi\ll tT^{-1/2}\log T\ll T^{-1/4+\varepsilon}$ it is sufficient take
the finite number terms $(<4)$ from this series.

As the result we come to the table integrals
$$
T\int_0^{\infty}K_0(r\xi)e^{-\mu_m\xi^2/8}\xi^{2 l+1}\,d\xi,\ \
l=0,1,2,\dots                                                      \tag 6.19
$$
or
$$
T\frac{\partial}{\partial l}\int_0^{\infty}K_0(r\xi)
e^{-\mu_m\xi^2/8}\xi^{2 l+1}\,d\xi\Big\vert _{l=0,1,2,\dots}       \tag 6.20
$$

The integral (6.19) equals to
$$
\frac{T}{2r}(\frac{8}{\mu_m})^{l+1/2} \Gamma^2 (l+1)
\exp\left(\frac{r^2}{\mu_m}\right)W_{-l-1/2,0}
\left(\frac{2r^2}{\mu_m}\right)                                      \tag 6.21
$$
where $W_{-l-1/2,0}(x)$ denotes the Whittaker function (which
exponentially decreseas when $x\to+\infty$).

If $r^2\gg|\mu_m|$ then
$$
\Big\vert \exp(\frac{r^2}{\mu_m})W_{-l-1/2,0}
\left(\frac{2r^2}{\mu_m}\right)\Big\vert\ll
\left(\frac{|\mu_m|}{r^2}\right)^{l+1/2}                           \tag 6.22
$$
so for the main term with $l=0$ we have
$$
T\Big|\int_0^{\infty}K_0(r\xi)e^{\mu_m\xi^2/8}\xi\,d\xi\Big|\ll
\frac T{r^2},\ \ r^2\gg |\mu_m|.                                 \tag 6.23
$$
At the same time for $r^2\ll |\mu_m|$ we have
$$
\Big|W_{-l-1/2,0}\left(\frac{2r^2}{\mu_m}\right)\Big|\ll
\frac r{\sqrt{|\mu_m|}}\log\frac r{|\mu_m|}                          \tag 6.24
$$
and it gives the bound
$$
T\Big|\int_0^{\infty}K_0(r\xi)e^{\mu_m\xi^2/8}\xi\,d\xi\Big|\ll
\log r,\ \ \ r^2\ll |\mu_m|;                                    \tag 6.25
$$
the differentiation with respect to $r$ gives the additional multiplier
$\log r$.

Exactly by the same way we come to the integrals
$$
T\int_0^{\infty}Y_0(r\xi)\Big\{\log\frac{t\xi}2(1+c_1 \,(t\xi)^2+\dots)
+b_1 \,(t\xi)^2+\dots)\Big\}e^{-\mu_m\xi^2/8}\xi^{2 l+1}\,d\xi         \tag 6.26
$$
(we use the power series for $Y_0(t\xi)$ in (6.17)). The leading term is
$$
\multline
T\,\int_0^{\infty} Y_0(r\xi) \exp (-\frac{\mu_m \xi^2}{8}) \xi^{2 l+1}\,d\xi =\\
=\frac{T}{r}(\frac{8}{\mu_m})^{l+1/2}\frac{\exp(-\frac{r^2}{8\mu_m})}{\sin
\pi \nu}\bigl(\Gamma(l+1) \cos \pi l
M_{l+1/2,0}(\frac{2r^2}{\mu_m})-W_{l+1/2}(\frac{2r^2}{\mu_m})\bigr) \\
=\frac{T}{\pi r}(\frac{8}{\mu_m})^{l+1/2} \exp
(-\frac{r^2}{8\mu_m})
\bigl((\Gamma'(l+1)-2\Gamma(l+1)
\frac{\Gamma'}{\Gamma}(1))M_{l+1/2,0}(\frac{2r^2}{\mu_m})+ \\
+\frac{\partial}{\partial
\varepsilon}M_{l+1/2,\varepsilon}(\frac{2r^2}{\mu_m})\Bigg
\vert_{\varepsilon=0})
\endmultline                             \tag 6.27
$$

The last equality follows as the limiting case of the relation between
$W_{\varkappa,\varepsilon}$ and $M_{\varkappa,\pm\varepsilon}$: if
$2\varepsilon$ is not integer we have
$$
W_{\varkappa,\varepsilon}(z)=
\frac{\Gamma(-2\varepsilon)}{\Gamma(1/2-\varkappa-\varepsilon)}
M_{\varkappa,\varepsilon}(z)+
\frac{\Gamma(2\varepsilon)}{\Gamma(1/2-\varkappa+\varepsilon)}
M_{\varkappa,-\varepsilon}(z).                                       \tag 6.28
$$
Now the assertion of Lemma 6.1 follows after the using of the expansion
of $M_{l+1/2,0}(z)$ for small and large $|z|$.
Since we have the similar asymptotic expansions for the kernels
$A_{00}(u,-i(l-1/2);1/2,\nu)$ we come to the same estimates for the
integrals wiht $b_{0,l}$.

\proclaim{Lemma 6.2}
Under the same assumptions what we had in the previous lemma we have for
$l\in L$
$$
\multline
\Big|\int\limits_{(1/2)}b_{0,l}(r;1/2,\nu;\rho,1/2)\Cal H_j(\rho+it)
\Cal H_j(\rho-it)\omega_T(\rho)\,d\rho\Big|\ll                         \\
\ll\cases
\log^2T,\ \ r\ll\,\,\, T_0^{1/2+\varepsilon}                                 \\
T\, r^{-2}\log^2r,\,\,\,\, r\gg T_0^{1/2+\varepsilon}.
\endcases
\endmultline                                                         \tag 6.29
$$
\endproclaim

\subhead
6.3.The integrals $ M_{j,0}$; terms with $B_{01}$
\endsubhead

The kernel $A_{01}(r,\sinh \xi/2;1/2,1/2)$ is exponentially small
for all $\xi\geqslant 0$ (see (5.82)); so all terms with
$\varkappa_j\geqslant t^{\varepsilon}$ with arbitrary small $\varepsilon >0$
give the contribution $o(1)$. For this reason the case
$r\ll t^{\varepsilon}$ have been considered only.

Integrating over the variable $\rho$ in the first line we get again
$$
\multline
\frac{1}{i}\int\limits_{(1/2)}\Cal H_j(\rho+it) \Cal H_j(\rho-it) B_{01}
(r,u;\rho,1/2;1/2,\nu)\omega_T(\rho)\,d\rho=\\
T\,\sum_{N,m\geqslant 1}\frac{t_j (N)\tau_{\nu}(N)}{(Nm^2)^{1/2+iT_0}}\times\\
\times\int_0^{\infty}A_{01}(r,\sinh \xi/2;1/2,1/2)\times\\
\times \sqrt{\sinh \xi}A_{01}(u,\frac{1}{\sinh \xi/2};1/2,\nu)
\psi(Nm^2\sinh^2 \xi/2)\frac{d\xi}{\sinh \xi/2}
\endmultline                        \tag{6.30}
$$
($ \psi$ is defined by (6.14)).

If $u=-i(l-1/2), l \geqslant 2$, we have (see (5.14))
$$
|\sqrt{\sinh \xi}A_{01}(-i(l-1/2),\frac{1}{\sinh \xi/2};1/2,\nu)|\ll
\min(t^3 \xi^{7/2},t^{-1/2}).               \tag{6.31}
$$

Exactly the same bound we have for the integral
$$
\int_{-\infty}^{\infty}\sqrt{\sinh \xi}A_{01}(u,\frac{1}{\sinh
\xi/2};1/2,\nu)h(u)\,d\chi(u)
$$
because of we can integrate terms with $J_{\pm 2iu}$ from (5.7) over the
line $\r{Im} u=\mp 3/2$.

Furthermore (see (5.82)),
$$
|(2\tanh \xi/2)^{-1/2}A_{01}(r,\sinh \xi/2; 1/2,1/2)|\ll \sqrt \xi
e^{-\pi r}.                     \tag{6.32}
$$

These inequalities give the following bound for the integrals in (6.30)
with $u=-i(l-1/2)$:
$$
\multline
\ll e^{-\pi r}\int_0^{\infty}\min \Biggl((t\xi)^3,(t\xi)^{-1/2}\Biggr)
\psi(Nm^2\sinh^2 \xi/2)\,d\xi\\
\ll e^{-\pi r}(T
\sqrt{Nm^2})^{-1}\min\Biggl((\frac{t}{\sqrt{Nm^2}})^3,(\frac{m\sqrt
N}{t})^{1/2}\Biggr)
\endmultline                         \tag{6.33}
$$
(we use the fact of the exponential smallness $\psi(x)$ for
$|x-1|\geqslant \frac{1}{T}$).

The same estimate we have for the integral with $h(u)$.

As the result we have
\proclaim{Lemma 6.3} Under assumptions of Lemma 6.1 the integral on the
left side (6.30) is $O(\log^2 t)e^{-\pi r}$ for $u=-i(l-1/2), l\geqslant
2$; the same estimate we have for the result of the integration of this
function over $u $ with the weight $h(u)\,d \chi(u)$.
\endproclaim

Really, for both cases we have the bound
$$
\sum_{N,l\geqslant 1}\frac{|t_j(N)\tau_{\nu}(N)|}{Nl^2}\,e^{-\pi
r}\min((\frac{t}{\sqrt{Nl^2}})^3, (\frac{l\sqrt N}{t})^{1/2}) \tag{6.34}
$$
and it rests to note that $|t_j(N)|$ are bounded in average for
$\varkappa\ll t^{\varepsilon}$.

\subhead
6.4.The integrals $B_{10}$
\endsubhead

\proclaim{Lemma 6.4}
Under the same conditions we have
$$
\Big|\int\limits_{(1/2)}B_{10}(r,u;\rho,1/2;1/2,\nu)\Cal H_j(\rho+it)
\Cal H_j(\rho-it)\omega_T(\rho)\,d\rho\Big|\ll\,\,\, e^{-\pi r/2} \log T
                                                        \tag 6.35
$$
\endproclaim

We write $B_{10}$ as $I_1+I_2+I_3$, where
$$
\multline
I_1=\int_0^{\pi/2}\sqrt{\sin\xi}A_{00}\left(u,
\frac 1 {\sin\xi/2};1/2,\nu\right)(2\text{tg}\xi/2)^{-1/2}\times\\
\times A_{10}(r,\sin\xi/2;1/2,1/2)(\sin\xi/2)^{-2\rho}\,d\xi,
\endmultline
\tag 6.36
$$
$$
\multline
I_2=\int_0^{\pi/2}\sqrt{\sin\xi}A_{00}\left(u,
\frac 1 {\cos\xi/2};1/2,\nu\right)(2\text{tg}\xi/2)^{1/2}\times\\
\times A_{10}(r,\cos\xi/2;1/2,1/2)(\cos\xi/2)^{-2\rho}\,d\xi,
\endmultline
\tag 6.37
$$
$$
\multline
I_3=\frac 12\int_0^{\infty}\sqrt{\sinh\xi}A_{00}\left(u,
\frac 1 {\cosh\xi/2};1/2,\nu\right)(2\tanh\xi/2)^{-1/2}\times\\
\times A_{10}(r,\cosh\xi/2;1/2,1/2)(\cosh\xi/2)^{-2\rho}\,d\xi.
\endmultline                                                 \tag 6.38
$$

Using (5.62) and (5.42) we see that for all $r$ the first integral $I_1$
is $O(\exp(-\frac{\pi}{2}t))$. At the same time for $r$'s very large we can move
the line integration over $u$ to ensure the estimate $O(r^{-3}
\exp(-\frac{\pi}{2} t)$;
so the contribution of this term will be exponentally small.

For the second integral we use the asymptotic formulas (5.67) and (5.37).
The part of this integral with $\xi\geqslant t^{-1+\varepsilon}$ for any
fixed (small) $\varepsilon>0$ gives the exponentially small contribution
again. In the interval $\xi\leqslant t^{-1+\varepsilon}$ (where
$\cos\xi/2=1+O(1/t)$) we integrate over $\rho$ under the sign of the
integration over $\xi$. As before, we come to the integral (6.7) and
after that we have the integral
$$
\frac{T}{\cosh(\pi r)}\int_0^{\delta}\xi K_0(t\xi)I_0(r\xi)
e^{-\frac{\pi T}8\xi^2}\,d\xi\ll\log T e^{-(\pi-\delta)r}.    \tag 6.39
$$

We consider the third integral by the same way; the final integral
(we use (5.59) and (5.50)) is
$$
\frac{T}{\cosh(\pi r)}\int_0^{\infty}\xi Y_0(t\xi)J_0(r\xi)
e^{-\frac{\pi T}8\xi^2}\,d\xi\ll\,\, e^{-\pi r}\log T.          \tag 6.40
$$

\subhead
6.5.Terms with
$B_{11}(r,u;\rho,1/2;1/2,\nu)\ (=B_{10}(r,u;1/2,\nu;\rho,1/2))$
\endsubhead

The estimation of this last integral is more complicated; it would be
realized by the following way.

Firstly we consider $B_{11}(r,u;\rho,1/2;1/2,\nu)$ for
$r\leqslant T_0-2 T^{1+\varepsilon}$; for this case $B_{11}$ is
exponentially small.

For $r\gg T_0$ the integral $B_{10}$ instead of $B_{11}$ will be
considered; as before we integrate over $\rho$ in the first line.

In contrast with the previous integrals the result of this integration
is the "long" sum. But we have the opportunity use "the convolution
formulas" (that is the replacement of the summation over $\varkappa_j$
by the integration over this variable). Asymptotically the last integrand
is the product of two $\delta$-functions with the different supports;
it allows us to finish our estimates for the sums over discret and continuous
sprectrum.

\subhead
6.5.1.The exponential case
\endsubhead

\proclaim{Lemma 6.5}
Let $T_0,\,\tau,\,t$ are sufficiently large, $T=T_0^{1-\varepsilon}$,
$t=o(T^{1/4})$; then we have for $r\leqslant T_0-2T^{1+\varepsilon}$
$$
\Big|\int\limits_{(1/2)}B_{11}(r,u;1/2,\nu;\rho,1/2)\Cal H_j(\rho+it)
\Cal
H_j(\rho-it)\omega_T(\rho)\,d\rho\Big|\ll\exp(-T^{\varepsilon}).
                                                           \tag 6.41
$$
\endproclaim

First of all, we can assume here
$|\tau-T_0|\leqslant T^{1+\varepsilon}\ (=T_0^{1-\varepsilon^2})$,
since $|\omega_T(\rho)|\ll\exp\left(-\frac{|\tau-T_0|}T\right)$.

Now we use (4.20) (with $s=\mu=1/2,\ \rho=1/2+i\tau,\ \nu=1/2+it$).

In this expression (see (4.13) and (4.14)) for $0\leqslant x<1$

$$
\multline
(2\pi)^{2\rho-1}A_{10}(r,\sqrt x;\rho,\nu)=   \\
=\frac{\cosh(\pi r)}{2\cos(\pi\nu)}
\Bigg\{\frac{\Gamma(\rho+it+ir)
\Gamma(\rho+it-ir)}{\Gamma(2\nu)}\times \\
\times x^{\nu}F(\rho+it+ir,\rho+it-ir;2\nu;x)-                \\
-\frac{\Gamma(\rho-it+ir)\Gamma(\rho-it-ir)}{\Gamma(2-2\nu)}\times  \\
\times x^{1-\nu}F(\rho--it+ir,\rho-it-ir;2-2\nu;x)\Bigg\}
\endmultline                                                \tag 6.42
$$
and for $x>1$
$$
\multline
(2\pi)^{2\rho-1}A_{10}(r,\sqrt x;\rho,\nu)=      \\
=\frac{ix^{1/2-\rho}\sin(\pi\nu)}{\sinh(\pi r)}
\Bigg\{\frac{\Gamma(\rho+it+ir)
\Gamma(\rho-it-ir)}{\Gamma(1+2ir)}\times \\
\times x^{-ir}F(\rho+it+ir,\rho-it+ir;1+2ir;\frac 1x)-        \\
-\frac{\Gamma(\rho+it-ir)\Gamma(\rho-it-ir)}{\Gamma(1-2ir)}\times \\
\times x^{ir}F(\rho+it-ir,\rho-it-ir;1-2ir;\frac 1x)\Bigg\}.
\endmultline                                           \tag 6.43
$$

After the multiplication by $\omega_T(\rho)$ we have in these
expressions the multipliers in front of the hypergeometric functions
which are not larger than
$$O(\frac 1{\sqrt
t})\exp(-\pi(\tau-t-r+\frac{|\tau-T_0|}{T})).         \tag{6.44}
$$

If $r\leqslant T_0-2T^{1+\varepsilon}$ (note that
$T^{1+\varepsilon}=T_0^{1-\varepsilon^2}=o(T_0)$) then this bound is not
larger than $O(\exp(-\pi T^{\varepsilon}))$; so it is sufficient to see
that these hypergeometric functions are bounded.

For the function
$$
F=(1-x)^{\rho}x^{\nu}F(\rho+\nu-1/2+ir,\rho+\nu-1/2-ir;2\nu;x)   \tag 6.45
$$
we have the differential equation
$$
F''+\left(\frac{p(x)}{x^2(1-x)^2}+
\frac{1-x+x^2}{4x^2(1-x)^2}\right)F=0                        \tag 6.46
$$
with
$$
p(x)=t^2(1-x)+(\tau^2-r^2(1-x))x;                           \tag 6.47
$$
this $p$ is positive for all $x\in[0,1]$ if $r\leqslant \tau$. For this
reason we have for $F$ the Liouville--Green approximation
$$
F\cong\sqrt t\frac{x^{1/2}(1-x)^{1/2}}{p^{1/4}}\exp(i\xi),\ \
\xi=\int_0^x\left(\frac{\sqrt{p(y)}}{1-y}-t\right)\frac{dy}{y}+t\log x.
                                                             \tag 6.48
$$

As the consequence we have (we use the Stirling expansion) for $r<\tau-t$
$$
|A_{10}(r,\sqrt x;\rho,\nu)|\ll e^{-\pi(\tau-t-r)}x^{1/2},\ \ x<1.\tag 6.49
$$

By the similar way we find for $x>1$
$$
|A_{10}(r,\sqrt x;\rho,\nu)|\ll e^{-\pi(\tau-t-r)}.            \tag 6.50
$$

Using again the corresponding differential equation we see from (4.10)
and (4.11) that
$$
\Bigg|A_{01}\left(u,\frac 1{\sqrt x};1-\rho,1/2\right)\Bigg|\ll
\cases
1,\ \ \ \ \ \ \ \ \ \ x\ll 1                \\
\frac 1{\sqrt x}\log x,\ x\gg 1.
\endcases                                                       \tag 6.51
$$

Now (6.41) follows.

\subhead
6.5.2.The integration over $\rho$
\endsubhead

\proclaim{Proposition 6.2}
$$
\multline
\frac1i\int\limits_{(1/2)}B_{10}(r,u;1/2,\nu;\rho,1/2)
\Cal H_j(\rho+it)\Cal H_j(\rho-it)\omega_T(\rho)\,d\rho=\\
=\frac T{2\pi}\sum_{N,l\geqslant 1}
\frac{t_j(N)\tau_{\nu}(N)}{(Nl^2)^{1/2+iT_0}} f_{N,l}(r,u)
\endmultline                                      \tag{6.52}
$$
where with $\psi$ from (6.14)
$$
f_{N,l}(r,u)= \int\limits_0^{\infty}A_{11}(r,\sqrt x;1/2,1/2)
A_{01}(u,\frac 1{\sqrt x};1/2,\nu)\psi(Nl^2x)x^{-1/2-iT_0}\,dx.
                                                      \tag 6.53
$$
\endproclaim

On the line $\r{Re} \rho=1/2$ we can change the order of integration.
After that we have the integral
$$
\int\limits_{(1/2)}\Cal H_j(\rho+it)\Cal H_j(\rho-it)x^{-\rho}
\omega_T(\rho)\,d\rho
$$
with $x>0$. Here we can integrate over the line $\r{Re}\rho=2$,
where the Hecke series are absolutely convergent. Using the Ramanujan
integral again we come to the expression
$$
\sum_{n,m\geqslant 1}\frac{t_j(n)t_j(m)}{(nm)^{1/2+iT_0}}
x^{-1/2-iT_0}\psi(nmx)\left(\frac{n}{m}\right)^{it}.         \tag 6.54
$$
Now we have
$$
t_j(n)t_j(m)=\sum_{l\backslash(n,m)}t_j\left(\frac {nm}{l^2}\right)
                                                            \tag 6.55
$$
and after the corresponding change $n$ by $nl$, $m$ by $ml$ and
$mn$ by $N$ we get (6.48).

\subhead
6.5.3.The large $N$'s in (6.52)
\endsubhead

First of all we replace the series in (6.52) by the finite sum with
$N\leqslant N_0(T)$.

\proclaim{Proposition 6.3}
Let
$$
h_{N,l}(r)=\int\limits_{-\infty}^{\infty}f_{N,l}(r,u)h(u)\,d\chi(u);
                                                             \tag 6.56
$$
then for $N_0=t^3\sqrt T $ we have
$$
\sum_{j\geqslant 1}\alpha_j \Cal H^2_{j}(1/2)\sum_{Nl^2\geqslant
N_0}\frac{t_j(N)\tau_{\nu}(N)}{(Nl^2)^{1/2+iT_0}}
h_{N,l}(\varkappa_j)=O(T^{\varepsilon})
                                                               \tag{6.57}
$$
and the same estimate is valid for the corresponding integral over
the continuous spectrum.
\endproclaim

To estimate $|h_{N,l}|$ we rewrite the definition $f_{N,l}$.
We do the change of the variable $x\mapsto \frac{1}{Nl^2}
\exp(\frac{x}{T})$ in (6.52).

Furthermore, we write instead of the kernel $A_{01}(u,...)$  the function
$$
\frac{i\sin \pi \nu}{\sinh \pi
u}\frac{\Gamma(\nu+iu)\Gamma(1-\nu+iu)}{\Gamma(1+2iu)}x^{iu}
F(\nu+iu,1-\nu+iu;1+2iu;-x);
$$
(this replacement is possible since $h(u)$ is the even function) and we
integrate this function on the line $\r{Im}\, u=-3/2$. Now we have the
representation
$$
\multline
T(Nl^2)^{1/2-iT_0}\,h_{N,l}(r)=\int\limits_{\r{Im}\,u=-3/2}\frac{i\sin
\pi \nu}{\sinh \pi
u}\frac{\Gamma(\nu+iu)\Gamma(1-\nu)}{(Nl^2)^{iu}\Gamma(1+2iu)}h(u)\times \\
\times\int_{-\infty}^{\infty}A_{11}(r,\frac{1}{l \sqrt
N}\exp(\frac{x}{2T});1/2,1/2)\times   \\
\times F(\nu+iu,1-\nu+iu;1+2iu;-\frac{1}{Nl^2}\exp(-\frac{x}{2T}))\times
\\
\exp(\frac{1+2iu-2iT_0}{2T}x) \frac{dx}{\cosh x}\,d\chi(u)
\endmultline                                    \tag{6.58}
$$

The part of this integral with $|u|\geqslant T^{\varepsilon}$ or
$|x|\geqslant T^{\varepsilon}$ gives $o(1)$ in the final result.

Since $Nl^2>t^3\,\sqrt T$ we can use the power series for the hypergeometric
functions in (6.58).

For the kernel $A_{11}$ we use the asymptotic expansion (5.29). Let
$|x|\leqslant T^{\varepsilon}$ and $N$ be large; we define the positive
$\xi$ by the equation
$$
\sinh \xi/2 =\frac{1}{l\sqrt N}\exp(\frac{x}{2T})=\frac{1}{l\sqrt
N}(1+\frac{x}{2T}+\frac{x^2}{8T^2}+... )      \tag{6.59}
$$
From (5.29) we find $A_{11}=O(\xi |\log (r\xi)|)$ for $r\leqslant l\sqrt
N$, $A_{11}=O(\sqrt{\frac{\xi}{r}})$ for $Nl^2\leqslant r \leqslant
T^{1+\varepsilon} Nl^2$ and for $r\geqslant T^{1+\varepsilon} Nl^2$ we
have the asymptotic series with the main term
$$
A_{11}(r,\sinh \xi/2;1/2,1/2) \approx -\frac{1}{\sqrt \pi}(\tanh
\xi/2)^{1/2}\sinh (r\xi -\frac{\pi}{4})  \tag{6.60}
$$

If $r\gg T^{1+\varepsilon}_0 l\sqrt N$ we can integrate by parts any
times. Each integration gives the additional multiplier $(r^{-1})Tl\sqrt
N)$; as the result we have for any $M\geqslant 2$
$$
\multline
\frac{1}{l\sqrt N}|h_{N,l}|\ll \frac{t^3}{TN^3 l^6}\log(rNl^2),\,\,r \ll
l\sqrt N,  \\
\ll \frac{t^3}{T \sqrt r (Nl^2)^{11/4}},\,\,\,\,  l\sqrt N \ll r \ll
T^{1+\varepsilon}_0 l \sqrt N,  \\
\ll \Biggl(\frac{T l \sqrt N}{r}\Biggr)^M \frac{t^3}
{T \sqrt r (Nl^2)^{11/4}},\,\, r\gg T^{1+\varepsilon}_0.
\endmultline                     \tag{6.61}
$$

We have the known estimates ([13],\, [2]):
$$\sum_{\varkappa_j \leqslant P}\alpha_j \Cal H^4_j(1/2)\ll
P^{2+\varepsilon},                    \tag{6.62}
$$
$$
\sum_{\varkappa_j \leqslant P}\alpha_j t^2_j(N)\ll P^2
+N^{1/2+\varepsilon};                   \tag{6.63}
$$
for this reason we have for $P\gg N^{1/4}$
$$
\sum_{\varkappa_j\leqslant P}\alpha_j |\Cal H^2_j(1/2)t_j(N)|
\frac{1}{\sqrt{\varkappa_j}} \ll
\sqrt{P^{1+\varepsilon}(P^2+N^{1/2+\varepsilon})}\ll
P^{3/2+\varepsilon}.
                                 \tag{6.64}
$$

It means the sum in (6.57) is not larger than $O(\frac{t^3
T^{1/2+\varepsilon}}{N_0})$ for any $\varepsilon >0$; taking $N_0=t^3
T^{1/2}$ we come to the assertion of Proposition 6.3.

\subhead
6.5.4. The convolution formula
\endsubhead

Let for a given function $f$ and for an integer $N\geqslant 1$
$$
Z_N(f)=\sum_{j\geqslant 1}\alpha_j t_j(N)\Cal
H^2_j(1/2)f(\varkappa_j)+\frac{1}{\pi}\int_{-\infty}^{\infty}
\tau_{1/2+ir}(N)\frac{|\zeta(1/2+ir)^4}{|\zeta(1+2ir)|^2}\,f(r)\,dr
                                               \tag{6.65}
$$

This quadratic form in the Hecke series had been expressed in terms of
the convolution of the Fourier coefficients of the Eisenstein series
(Theorem 3.1 in [8]; here we take the special case $s=\nu=1/2$).
\proclaim{Theorem 6.1} Let $f(r)$ be the even regular function in the
strip $|\r{Im}\,r|\leqslant 5/2, f(\pm i/2)=0$ and $|f(r)|\ll |r|^{-B}$
for some $B>4$ when $r\to \infty$ in this strip. Then for any integer
$N\geqslant 1$ we have
$$
\multline
Z_N(f)=\frac{1}{\pi}\frac{d(N)}{\sqrt
N}\int_{-\infty}^{\infty}(-\log \frac{N}{2\pi} - 2\frac{\Gamma'}{\Gamma}(1/2+ir)-
2\frac{\Gamma'}{\Gamma}(1)+\zeta'(0))f(r)\,dr -\\
-\frac{1}{4\pi}\int_{-\infty}^{\infty}f(r)\frac{d\chi(r)}{\cosh \pi r}+
\frac{1}{\sqrt N}\sum_{n\ne
N}d(n)d(n-N)w_0(\sqrt{\frac{n}{N}})+ \\
+\frac{1}{\sqrt N}\sum_{n\geqslant 1}d(n)d(n+N)w_1(\sqrt {\frac{n}{N}}),
\endmultline                                              \tag{6.66}
$$
where $d(n)$ is the number of the positive divisors of $|n|$ and two
functions $w_j$ are defined by the equalities
$$
w_j(x)=\frac{1}{\pi
x}\int_{-\infty}^{\infty}A_{0j}(r,x;1/2,1/2)f(r)\,d\chi(r),\,\,j=0,1
                                                \tag{6.67}
$$
\endproclaim

\subhead
6.5.4.The construction of the suitable taste function
\endsubhead

For the case $\varkappa_j\geqslant T_0 - 2 T^{1+\varepsilon}$ we use the
representation (6.50).

Let us define
$$
\Omega(r)=\frac{2\beta}{\pi}\int_{T_0/2}^{\infty}\bigl((\cosh
\beta(r-u))^{-1}+(\cosh \beta
(r+u))^{-1}\bigr)\,du,\,\,\beta=T^{-1+\varepsilon}_0. \tag{6.68}
$$

This even function is regular in the wide strip $|\r{Im}| r
<\frac{\pi}{2\beta}$, for real $r$ this function is
$O(\exp(-T_0^{\varepsilon}))$ for $|r|\leqslant T_0/4$ and for the
positive $r$ with $r\geqslant 3T_0/4$ we have
$$
1-c_1 \exp(-c_2 T_0^{\varepsilon})\leqslant \Omega(r)\leqslant 1
                                                         \tag{6.69}
$$
with the fixed positive constants $c_1,\,c_2$.

Furthermore, let
$$
q(r)=\frac{(r^2+1/4)^2 (r^2+9/4)^2}{(r^2+1/4)^2 (r^2 +9/4)^2 + M^2},
                                                             \tag{6.70}
$$
where the  fixed positive $M$ be so large that $q$ has no poles in the
strip $|\r{Im} r|\leqslant 5/2$.

With these notations we define
$$
F_{N,l}(r,u)=q(r)\Omega(r)\,f_{N,l}(r,u)                 \tag{6.71}
$$
($f_{N,l}(r,u)$ is defined by (6.51)).

It is obvious this function satisfies to all conditions of Theorem 6.1.

The main characteristic of this function is the following fact: the
result of the summation over the discret and continuous spectrum of the
quantities (6.50) we can replace by the sum
$$
\frac{T}{2\pi}\sum_{Nl^2\leqslant N_0}\frac{
\tau_{\nu}(N)}{(Nl^2)^{1/2+iT_0}}\,Z_N(F_{N,l}), \,\,N_0=t^3
\sqrt T.                                                  \tag{6.72}
$$

Really, the difference
$$
F_{N,l}(r,u)-f_{N,l}(r,u)=(q(r)\Omega(r)-1)f_{N,l}(r,u)
$$
after the summation over $N,l$ gives the exponentially small
contribution for $r\leqslant T_0/4$ (because of $\Omega$ is very small
for these values of $r$ and we have case which had been considered in
6.5.1.). But for $r\gg T_0$ this difference is $O(r^{-8})|f_{N,l}(r,u)|$.
In any case we have for the integral (6.58) the rough estimate
$h_{N,l}(r)=O((rl\sqrt N)^{-1})$; with the additional multiplier
$r^{-8}$ it is sufficient to assert the smallness of the full sum with this
difference.
\subhead
6.5.6. The first sum with $w_0$ (the case $n\geqslant N+1$)
\endsubhead

We use identity (6.66) to estimate $Z_{N,l}(F_{N,l})$; first of all we
consider the integral
$$
\frac{1}{\sqrt N}w_0(\sqrt {\frac{n}{N}}) =\frac{1}{\pi\sqrt
n}\int_{-\infty}^{\infty}A_{00}(r,\sqrt{\frac{n}{N}};1/2,1/2)\,
F_{N,l}(r,u)\,d\chi (r).                   \tag{6.73}
$$

After the replacement $f_{N,l}$ by the integral representation (6.51)
(where we do the change of the variable $x\mapsto \sinh^2 \eta/2$) we
come to the double integral (we write $\cosh^2 \xi/2$ instead of
$\frac{n}{N}, \,\xi>0$)
$$
\multline
\frac{1}{\pi \sqrt n}\int_{-\infty}^{\infty}
A_{00}(r,\cosh \xi/2;1/2,1/2)\int_0^{\infty}A_{11}(r,\sinh
\xi/2;1/2,1/2)\times   \\
\times A_{01}(u,\frac{1}{\sinh \eta/2};1/2,\nu)\psi (Nl^2 \sinh^2 \eta/2)
(\sinh \eta/2)^{-2iT_0}\cosh \eta/2 \,d \eta \times   \\
\times q(r)\Omega(r)\,d\chi(r).
\endmultline                                       \tag{6.74}
$$

Let $\delta =T^{-1+\varepsilon}$ with small $\varepsilon >0$; we write
the integral on the right side (6.74) as $v_1+v_2+v_3$, where
$$
v_1(n,l)=\int_{-\infty}^{\infty}\int_0^{\xi-\delta}(...)\,d\eta\,d\chi(r),
v_2=\int_{-\infty}^{\infty}\int_{\xi-\delta}^{\xi+\delta} (...)\,d\eta\,d\chi(r),
                                                            \tag{6.75}
$$
and in $v_3$ the integration is doing on $\eta\geqslant \xi+\delta$.

\proclaim{Proposition 6.4} For $j=1$ and $j=3$ we have
$$
\sum_{n\geqslant N+1}d(n)d(n-N)|v_j(n,l)| \ll \exp(-T^{\varepsilon})
                                                          \tag{6.76}
$$
\endproclaim

Firstly we note that for all $n\geqslant N+1 $
$$
\xi(n)\gg T^{-1/4}t^{-3/2}                        \tag{6.77}
$$
(if $n$ be near to $N$ we have $\cosh \xi/2 = 1+ \frac{1}{8}\xi^2
+... =1+\frac{n-N}{N}+...$; so $\xi(n)\gg \frac{1}{\sqrt N}\gg
N_0^{-1/2}$).

Now for $\eta \geqslant \xi+\delta$ we have $Nl^2\sinh^2\eta/2
\geqslant l^2(n-N+\delta \sqrt N(n-N))$.

For this reason for $\eta\geqslant \xi+\delta$
$$
\psi(Nl^2 \sinh^2 \eta/2)\ll \exp(-l^2(n-N+\delta \sqrt N(n-N))T),
                                                       \tag{6.78}
$$
so the series over $n$ is convergent and this sum is
$O(\exp(-T^{\varepsilon}))$.

For the case $j=1$ we have the same result if $\xi(n)=O(1)$.

For large $n$ we have $\xi(n)\geqslant \log \frac{n}{N}$; writing
$Y_0(r\xi)$ as $\frac{1}{2i}(H_0^{(1)}-H_0^{(2)})$ we integrate the
Hankel function $H_0^{(1)}$ on the line $\r{Im}\, r=\Delta >1/2$ and the
other function on the line $\r{Im}\, r=-\Delta<-1/2$ (we assume $\Delta$ be
near to 1/2). It gives the convergence of our series and the smallness
$\psi(Nl^2\sinh^2 \eta/2)$ for $\eta\leqslant \xi-\delta$ gives the same
result.

It rests estimate the sum with $v_2$.
\proclaim{Proposition 6.5}
$$
\sum_{l^2(n-N)\geqslant 2}d(n)d(n-N)v_2(n,l) \ll \exp(-T^{\varepsilon})
                                                               \tag{6.79}
$$
\endproclaim

Really, in this integral we have $|\xi-\eta|\leqslant
T^{-1+\varepsilon}$, so $Nl^2\sinh^2 \eta/2$ is near to $l^2(n-N)$.

If $l^2(n-N)\geqslant 2$ we have $\psi(Nl^2\sinh^2\eta/2)\ll
\exp(-Tl^2(n-N))$ and any rough estimate for the result of integration
over $r$ and $\eta$ gives the desired result.

Now we consider the integral with $n=N+1$ and $l=1$.

\proclaim{Proposition 6.6} Let $B(r;\xi,\eta)$ denotes the integrand
in (6.74); then
for $n=N+1$ and $l=1$ we have
$$
\lim_{P\to
\infty}\int_0^P\int_{\xi-\delta}^{\xi+\delta} B(r;\xi,\eta)\,d\eta\,d\chi(r)
=\frac{1}{\pi}A_{01}(u,\sqrt
N;1/2,\nu)N^{iT_0}(1+\frac{1}{N})^{1/2} +O(\frac{N}{T^2})   \tag{6.80}
$$
\endproclaim

We change the order of the integration and integrate over $r$ in the
first line. The part of this integral with $r\leqslant T$ gives the
exponentially small contribution and may be omitted.

In the interval $T\leqslant r \leqslant P$ the quantity $r\xi$ is large
since $\xi\gg T^{-1/4} t^{-3/2}$. For this reason both expansions (5.93)
and (5.29) may be rewritten in the trigonometric form. Using the asymptotic
expansions for the Bessel functions we have
$$
\multline
A_{00}(r, \cosh \xi/2;1/2,1/2)=-(\pi r \tanh \xi/2)^{-1/2}
\bigl((1+\Cal E_1) \sin(r\xi-\frac{\pi}{4}) +     \\
+\frac{\Cal E_2}{r\xi}\cos(r\xi-\frac{\pi}{4})\bigr),
\endmultline                                           \tag{6.81}
$$
$$
\multline
A_{11}(r,\sinh \eta/2;1/2,1/2)=   \\
=-(\frac{1}{\pi r})^{1/2}\bigl((1+\Cal D_1)\sin
(r\eta)-\frac{\pi}{4}) +   \\
+\frac{\Cal D_2}{r \eta} \cos(r \eta -\frac{\pi}{4}) \bigr)
\endmultline                                       \tag{6.82}
$$
where $\Cal E_j$ and $\Cal D_j$ are the asymptotic series in
the inverse degrees of $r^2, (r\xi)^2$ and $r^2,(r\eta)^2$
correspondingly.

In these expansions we have $\Cal E_1=O((r\xi)^{-2})$ and $\Cal
E_2=O(1)$ and the similar order have $\Cal D_j$.

As the result we have
$$
\multline
A_{00}(r,\cosh \xi/2;1/2,1/2) A_{11}(r,\sinh
\eta/2;1/2,1/2)=         \\
=\frac{1}{2\pi r}(\frac{\tanh \eta/2}{\tanh
\xi/2})^{1/2}\Biggl\{\cos r(\xi-\eta)\,\,(1+\Cal E_3)+       \\
+\frac{\sinh r(\xi-\eta)}{r}((1+\Cal E_1)\frac{\Cal
D_2}{\eta}-(1+\Cal D_1)\frac{\Cal E_2}{\xi}+...\Biggr\},
\endmultline                                    \tag{6.83}
$$
where $\Cal E_3$ is the asymptotic series in the inverse degrees of
$r^2,(r\xi)^2, (r\eta)^2$ ($\Cal E_3 = O((r\xi)^{-2})+(r\eta)^{-2}$) and
the unwritten terms contain $\cos r(\xi+\eta)$ or $\sin r(\xi +\eta)$.

These unwritten terms contribute $O(T^{-M})$ for any fixed $M\geqslant
4$ (since we can integrate these terms by parts any times).

For large positive $r$ we have $d\chi(r)=\frac{2}{\pi}\,r\,(1+O(e^{-\pi
r}))\,dr$ and
$$
\multline
\int_T^P A_{00}(r,\cosh \xi/2;1/2,1/2)A_{11}(r,\sinh
\eta/2;1/2,1/2)q(r)\Omega(r)\,d\chi(r)=   \\
=\frac{1}{\pi^2}(\frac{\tanh \eta/2}{\tanh \xi/2})^{1/2}\Biggl\{\frac{\sin
P(\xi-\eta)}{\xi-\eta}(1+\Cal E_3(P))-  \\
-\int_T^P\frac{\sin
r(\xi-\eta)}{\xi-\eta}((1+\Cal E_3)q \Omega \tanh \pi r)'\,dr+ \\
+\int_T^P \frac{\sin r(\xi-\eta)}{\xi-\eta}(\frac{(1+\Cal E_1)\Cal
D_2}{\eta}-\frac{(1+\Cal D_1)\Cal E_2}{\xi})q\Omega\tanh \pi r\,dr+  \\
+O(\exp(-T))\Biggr\}.
\endmultline                                   \tag{6.84}
$$

After multiplication by $A_{01}(u,...)$ the first term gives in the
limit $P\to \infty$
$$
\multline
A_{01}(u,\frac{1}{\sinh \xi/2};1/2,\nu) \cosh\xi/2 (\sinh
\xi/2)^{-2iT_0}\psi(N\sinh^2\xi/2)= \\
=A_{01}(u,\sqrt N;1/2,\nu) N^{iT_0}(1+\frac{1}{N})^{1/2}
\endmultline                                      \tag{6.85}
$$

Note that all asymptotic series $\Cal E_1, \Cal E_2,...$ allow the term
by term differentiation (because of it is sufficient take the finite
number terms in asymptotic expansions for our kernels and for expansions
of the Bessel function this fact is known very well). In particular,
$\Cal E'_3$ is the asymptotic series with the main term
$-\frac{9}{64}\frac{1}{r^3\xi^2}$ and $\Cal E'_1, \Cal D'_2 $ are
the asymptotic series with the main terms of the order $r^{-3}\xi^{-2},
r^{-3}\eta^{-2}$ correspondingly. The derivatives of $q$ and $\Omega$
have the smaller order.

Using the fact that $\int_P^r u^{-1}\sin au \,du$ is bounded for all
positive $P,r$ and $a$, we find the estimate $O(T^{-2}\xi^{-3})$ for the
last integral in (6.85). So the result of the integration over $\eta$
gives not larger than
$$
\multline
\max_{|\xi-\eta|\leqslant T^{-1+\varepsilon}}|A_{01}(u,\frac{1}{\sinh
\eta/2};1/2,\nu)|\times\frac{1}{T^2
\xi^3}\int_{\xi-\delta}^{\xi+\delta}\psi(N\sinh^2 \eta/2)\,d\eta  \ll
\\ \ll  \frac{1}{T^3\xi^3 \sqrt N}.
\endmultline
$$

Finally, the second integral is $O(\frac{1}{T^2\xi^2})=O(\frac{N}{T^2})$;
it gives us (6.80).
\subhead
6.5.7. The second sum with $w_0$ ($n\leqslant N-1$)
\endsubhead

\proclaim{Proposition 6.7} Let $w_0$ be defined by (6.73); then
$$
\sum_{n\leqslant N-1, l\geqslant 1}d(n)d(N-n)w_0(\sqrt{\frac{n}{N}}) \ll
\exp(-T^{3/4} t^{-3/2})                \tag{6.86}
$$
\endproclaim

Really, we have (see (5.89))
$$
A_{00}(r,\cos \xi/2;1/2,1/2)\approx
\sqrt{\frac{2\xi}{\tanh\xi/2}}\,K_0(r\xi),\,\,0\leqslant \xi\leqslant
3\pi/4.                \tag{6.87}
$$
Let $\xi$ be defined by the equation $\cos\xi/2=\sqrt{\frac{n}{N}}$; then
we have $\xi \gg \sqrt{\frac{N-n}{N}}\gg N^{-1/2}$ and it means that
$\xi \gg T^{-1/4} t^{-3/2}$ since $N\leqslant N_0
=T^{1/2} t^3$.
       
We can assume $r\gg T$; then $r\xi\gg T^{1/4-\varepsilon}t^{-3/2}$ and
(6.86) follows after any rough estimate of the result of summation and
integration.
\subhead
6.5.8. The sum with $w_1$
\endsubhead

\proclaim{Proposition 6.8} Let $w_1$ be defined by (6.67) with $F_{N,l}$
instead of $f$; then
$$
\sum_{n\geqslant 1}d(n)d(n+N) w_1(\sqrt{\frac{n}{N}})  \ll \exp(-T).
                                                          \tag{6.88}
$$
\endproclaim

This estimate follows from (5.82) and the definition of $\Omega(r)$.
We have $\Omega(r)\ll \exp(-T)$ for $r \leqslant T$ and at the same time
$$
A_{01}(r,\sinh \xi/2;1/2,1/2)\approx \frac{\pi}{2\cosh \pi r}\sqrt \xi
J_0(r\xi) \ll \exp(-\pi r).           \tag{6.89}
$$

It is sufficient for $n\leqslant NT^M$ for any fixed $M\geqslant 2$. For
very large $n$ we integrate the first term from (4.10) on the line
$\r{Im} \,=-1/2-\varepsilon_0$ and the second one on the line
$\r{Im}\,r=\,1/2+\varepsilon_0$ with small positive $\varepsilon_0$. Due
to the multipliers $(\frac{n}{N})^{\pm ir}$ in the denominator it
ensures the convergence and we come to (6.86) again.

\subhead
6.5.9. The first term of the convolution
\endsubhead

\proclaim {Proposition 6.9} Let $W_{N,l}$ denotes the first integral on
the wright side (6.66) for the case $f=F_{N,l}(r)$; then we have
$$
|W_{N,l}| \ll \frac{\log T}{T^3}.               \tag{6.90}
$$
\endproclaim

The part of this integral with $|r|\leqslant T$ is exponentially small
because of the smallness $\Omega$. For $r\gg T$ we use the asymptotic
expansion (5.29),
$$
A_{11}(r,\sinh \xi/2;1/2,1/2)=-\frac{\pi}{2}\sqrt{2\tanh
\xi/2}\bigl\{\sqrt \xi Y_0(r\xi)\,\Cal E_1 + (\sqrt \xi
Y_0(r\xi))'\,\Cal E_2\bigr\},                 \tag{6.91}
$$
where $\Cal E_j$ are the polynomials in $r^{-2}$ of the degree
$\leqslant 2$ (all other terms give $o(1)$ in the final result).

Firstly we consider the integral over $r$ in the finite interval
$(T,P)$.

We write
$$
rY_0(r\xi)=\frac{1}{\xi}\frac{\partial}{\partial r}rY_1(r\xi),
\frac{\partial}{\partial \xi}\sqrt \xi
Y_0(r\xi)=\sqrt{\frac{r}{\xi}}\frac{\partial}{\partial r}\sqrt r
Y_0(r\xi)                                         \tag{6.92}
$$
and integrate by parts from $T$ up to $P$.

The integrated terms at $r=T$ are exponentially small and at $r=P$ we
have the main term
$$
(c_1(N)+\frac{\Gamma'}{\Gamma}(1/2+iP))(2\xi^{-1}\tanh
\xi/2)^{1/2}PY_1(P\xi)\Cal E_1(P)q(P)\Omega(P)\tanh \pi P,    \tag{6.93}
$$
where $c_1(N)$ is the linear function in $\log N$.

Here
$$
PY_1(P\xi)=-\frac{\partial}{\partial \xi}Y_0(P\xi)
                                                     \tag{6.94}
$$
and we can integrate by parts over $\xi$. The integrated terms at
$\xi=0$ and at $\xi=\infty$ are zeros (due to $\psi(Nl^2\sinh^2 \xi/2)$
in the integrand) and the new integral will contain
$Y_0(P\xi)=O(\frac{1}{\sqrt{P\xi}})$. So this term give zero at the limit
$P\to{+\infty}$ and it rests estimate the convergent integral
$$
\multline
\int_T^{\infty}\int_0^{\infty}rY_1(r\xi)\frac{\partial}{\partial
r}\bigl\{(c_1(N)+\frac{\Gamma'}{\Gamma}(1/2+ir)+\frac{\Gamma'}
{\Gamma}(1/2-ir))\Cal E_1 q \Omega \tanh \pi r \bigr\}\times  \\
A_{01}(u.\frac{1}{\sinh \xi/2};1/2,\nu)\psi(Nl^2\sinh^2\xi/2)(\sinh
\xi/2)^{-2iT_0}(\frac {\sinh \xi}{2\xi})^{1/2}\,d\xi\,dr
\endmultline                                              \tag{6.95}
$$
(we change the variable $x$ in the definition (6.54) by $\sinh \xi/2$).

Here the derivative in respect to $r$ equals
$\frac{2}{r}-\frac{a_1}{r^3}+...$, so after one additional integration
by parts we come to the expression
$$
\multline
\xi^{-1}Y_0(r\xi)\frac{\partial}{\partial r}r\frac{\partial}{\partial r}
\bigl\{(c_1(N)+\frac{\Gamma'}{\Gamma}(1/2+ir)+\frac{\Gamma'}{\Gamma}(1/2-ir))
\Cal E_1 q \Omega \tanh \pi r \bigr\}= \\
=\xi^{-1}Y_0(r\xi)(\frac{4a_1}{r^3}+...).
\endmultline                                        \tag{6.96}
$$

Since the kernel $A_{01}$ is bounded here (for the real $u$ and
$\nu=1/2+it, t\gg 1$), the result integration over $\xi $ is not larger than
$$
\frac{\log T}{T^2}\int_0^{\infty}\psi(Nl^2\sinh^2 \eta/2)\frac{d\xi}{\xi}\ll
\frac{\log T}{T^3}.                      \tag{6.97}
$$

Exactly by the same way we estimate the integral with $(\sqrt \xi
Y_0(r\xi))'$; using the second equality from (6.92) and integrating by
parts we come to the integral with integrand
$$
\sqrt{\frac{r}{\xi}}Y_0(r\xi)\frac{\partial}{\partial r}\bigl\{\sqrt r
(c_1(N)+\frac{\Gamma'}{\Gamma}(1/2+ir)+\frac{\Gamma'}{\Gamma}(1/2-ir))\Cal
E_2 q \Omega \tanh \pi r \bigr\}.         \tag{6.98}
$$

Here $\xi^{-1}\Cal E_2 =\frac{b_1(\xi)}{\xi r^2}+\frac{b_2(\xi)}{\xi
r^4}+...$,so we have the considered integral with the additional
multiplier $r^2$ in the denominator.

It is obviously that the same estimate is valid for the second integral
in (6.66) for the case $f =F_{n,l}$.

\subhead
6.6. The union of estimates
\endsubhead
\proclaim{Lemma 6.6} Let $\rho,\nu$ and the initial function $h$ are
taken with conditions of subsection 3.1; then we have
$$
\multline
\int\limits_{(1/2)}\bigl(Z^{(d)}(\rho,\nu;1/2,1/2)\vert
h_0,h_1)+Z^{(c)}(\rho,\nu;1/2,1/2\vert
h_0+h_1\bigr)\omega_T(\rho)\,d\rho =   \\
O(T^{1+\varepsilon}) +
\text{\{contribution of poles $\Cal Z(\rho;\nu,1/2+ir)\}$}
\endmultline                                             \tag{6.99}
$$
\endproclaim

For our specialization the coefficient $h_1$ is exponentially small for
$\varkappa_j\leqslant T_0-2T^{1+\varepsilon}$ (lemmas 6.5 and 6.4
together with representation (4.17)).

The explicit averaging and the replacement of the inner sum over large
$\varkappa_j$ by the convolution gives for this sum the estimate (the
main term of this convolution is (6.80))
$$
T\,\,\, O(\sum_{N\leqslant N_0}\frac{d^2(N)}{N}|\int\nolimits A_{01}(u,\sqrt
N;1/2,\nu)h(u)\,d\chi(r)|)\ll T \log^4 T.          \tag{6.100}
$$
All other terms from the convolution give $O(\frac{\log^2
T}{T^3}\sum_{N\leqslant N_0}\frac{d^2(N)}{N})\ll T^{-3+\varepsilon}$
(the union of (6.76),(6.79),(6.86),(6.88) and (6.90)). So the main
contribution gives the sum with the coefficient $h_0$.

For this sum we have (lemmas 6.2,\,\,6.2) the estimate
$$
\multline
|\int\limits_{(1/2)}Z^{(d)}(\rho,\nu;1/2,1/2\vert
h_0,0)\omega_T(\rho)\,d\rho|\ll  \\
\ll \log^2 T\sum_{\varkappa_j\leqslant R}\alpha_j \Cal H^2_j(1/2)
+T\sum_{R\leqslant \varkappa_j \leqslant 2T_0}\alpha_j \frac{\log^2
\varkappa_j}{\varkappa^2_j}\Cal H^2_j(1/2)
\endmultline                                   \tag{6.101}
$$
with $R=T_0^{1/2+\varepsilon}$. Using the asymptotic formula (58) from
[8] we see that the right side is $ O(T^{1+\varepsilon})$ for any fixed
$\varepsilon >0$.
\subhead
6.7.The contribution of the poles $\Cal Z(\rho;\nu,1/2)$
\endsubhead

Estimating the integrals
$$
\multline
\int\limits_{(1/2)}Z^{(c)}(\rho,\nu;1/2,\,1/2 | h^j) \omega_T
(\rho)\,d\rho =     \\
=\frac{1}{\pi}\int\nolimits \int\nolimits \Cal Z(\rho;\nu,1/2+ir) \Cal
Z(1/2;1/2,1/2+ir)\frac{h_j(r)}{|\zeta(1+2ir)|^2}\omega_T
(\rho)\,dr\,d\rho,
\endmultline
                                                  \tag{6.102}
$$
we can assume that for the case $j=0$ the integration over $r$ have been
done on the segment $|r|\leqslant 2T_0$ and for the case $j=1$ --
on the line segment $r\geqslant T_0 - 2 T^{1+\varepsilon}$, since $h_0$
is exponentially small for $r\geqslant 2T_0$ and $h_1$ is small for
$r\leqslant T_0 - 2 T^{1+\varepsilon}$.

For the similar reason we can assume that the integration over $\rho$
have been done on the segment $(1/2+iT_0-iT_1,\, 1/2+iT_0 +iT_1)$ with
$T_1=T^{1+\varepsilon}$.

Integrating over $\rho$ in the first line we change this segment by the
new path $\Cal C$ which made up of
\newline
\noindent
$\Cal C_1$: segment $(1/2 + iT_0 -iT_1, \,2+iT_0 -iT_1),
T_1=T^{1+\varepsilon}$;
\newline
\noindent
$\Cal C_2$: segment $(2+iT_0 -iT_1, \,2+iT_0 +iT_1)$;
\newline
\noindent
$\Cal C_3$: segment $(2+iT_0 +iT_1,\,1/2 +iT_0 +iT_1)$.

Two integrals over $\Cal C_j$ with $j=1$ and $j=3$ are exponentially
small; the integral over $\Cal C_2$ have been considered early.

Now we consider the contribution of four residues of the integrand at
the points $\rho=\nu+1/2 \pm ir$ and $\rho=3/2-\nu \pm ir$ where $\Cal
Z(\rho;\nu,1/2+ir)$ has the simple pole.

These residues are the following integrals:
$$
\multline
\zeta(2\nu)\int\limits_{|t \pm r-T_0| \leqslant T_1 } \zeta(2\nu\pm 2ir)
\zeta(1\pm 2ir) \frac{|\zeta(1/2+ir)|^4}{|\zeta(1+2ir)|^2}\times    \\
\times h_j(r;1/2,\nu;\nu+1/2 \pm ir, 1/2)\omega_T(\nu+1/2 \pm ir)\,dr
\endmultline                                            \tag{6.103}
$$
and the same integrals with $\nu$ replaced by $1-\nu$.

Using the known estimates of zeta-function on the unit line we see that
these integrals are not larger than
$$
o((\log T)^{10})\int_{T_0-T_1}^{T_0+T_1}|\zeta^4(1/2+ir)
h_j(r;1/2,\nu;\nu+1/2\pm ir,1/2) \omega_T(\nu+1/2 \pm ir)| \,dr.
                                                         \tag{6.104}
$$

It follows from the initial representations (2.19) and (2.20) that both
functions $h_j(r;1/2,\nu;\nu+1/2 \pm ir,1/2)$ are bounded for
$T_0-T_1\leqslant r \leqslant T_0+T_1$. So we have the fourth moment of
zeta-function and all these integrals are $O(T^{1+\varepsilon})$.

\head
\S 7. SUM OVER REGULAR CUSP FORMS
\endhead
\noindent
\subhead
7.1.The integral representations for the coefficients
\endsubhead

\proclaim{Lemma 7.1} Let the coefficient $g(k)$ be defined by the
equality (2.21) with the choice (3.1) for the parameters
$s,\nu,\rho,\mu$ and with $\Phi$ from (3.3) instead of $\Phi_N$ in
(2.21);then we have
$$
g(k)\equiv g(k;\nu,\rho)=(-1)^{k-1}g_0(k;\nu,\rho) +g_1(k;\nu,\rho),
                                                         \tag{7.1}
$$
where
$$
\multline
g_0(k;\nu,\rho)=\frac{2k-1}{\pi}\int_0^{\infty}A_{00}(i(k-1/2),\sqrt
x;1/2,1/2)\times   \\
\times\bigl(\int_{-\infty}^{\infty}A_{00}(u,\frac{1}{\sqrt
x;1/2,\nu})h(u)\,d\chi(u)-    \\
-\sum_{l\in L}(-1)^l c(l)
A_{00}(i(l-1/2),\frac{1}{\sqrt x};1/2,\nu)\bigr)\,x^{-\rho-1/2}\,dx
\endmultline                                 \tag{7.2}
$$
and
$$
\multline
g_1(k;\nu,\rho)=\frac{2k-1}{i\pi}\cosh \pi t \times       \\
\times\int\limits_{(1/2)}2^{2w-1}
\frac{\Gamma(k-1/2+w)}{\Gamma(k+1/2-w)}\gamma(\rho-w,\nu)
\gamma(1/2-w,1/2) \hat\Phi(2w-2\rho+1)\,dw
\endmultline                                      \tag{7.3}
$$
\endproclaim

The first representation follows from the equalities
$$
\multline
\frac{1}{2\pi}\int_0^{\infty}A_{00}(-i(k-1/2),\sqrt x;1/2,1/2)
x^{-w-1}\,dx=\\
=\gamma(w,k) \gamma(1/2-w,1/2)\cos \pi w \sin \pi w ,
\endmultline                                                   \tag{7.4}
$$
$$
\multline
\frac{1}{2\pi}\int_0^{\infty}A_{00}(u,\frac{1}{\sqrt x};1/2,\nu)
x^{-\rho+w-1/2}\,dx=\\
=\gamma(\rho-w,\nu) \gamma(w-\rho+1/2,1/2+iu) \cos \pi(\rho-w) \cos \pi
(w-\rho+1/2);
\endmultline                                     \tag{7.5}
$$
here $k\geqslant 2$ is an integer and the last equality holds for
$u=-i(l-1/2)$ if $l\geqslant 1$ be an integer.

The function $g_1$ is the part of (2.21) with $s=\mu=1/2$.

\noindent
\subhead
7.2. The coefficient $g_1$
\endsubhead
\proclaim{Proposition 7.1}
$$
|g_1(k;\nu,\rho)|\ll k^{-3} e^{-\tau}.           \tag{7.6}
$$
\endproclaim

Really, the integrand in (7.5) on the line $\r{Re}\, w=-3/2$ is
$O(k^{-4} \exp (-\pi (\tau +2|w|))\,|\rho-w|^5)$ for $\r{Im}\, w<0$,
$O(k^{-4} e^{-\pi \tau} \tau^5)$ for $0\leqslant \r{Im}\, w\leqslant
\tau-t$ and $O(k^{-4}e^{-\pi(|w|+t)} |w|^5)$ for $\r{Im}\, w>\tau-t$.

The sum over regular cusp forms with this coefficient $g_1$ gives $o(1)$
into the final result.
\noindent
\subhead
7.3. The kernel $A_{00}(i(k-1/2),x;1/2,1/2)$
\endsubhead

We substitute in (4.8) $r=i(k-1/2)$ with an integer $k$,$\nu=1/2$ and
$s=1/2+\varepsilon$; taking the limit $\varepsilon\to 0$ we come to the
definition
$$
\multline
A_{00}(i(k-1/2),x;1/2,1/2)=2x(\frac{\Gamma'}{\Gamma}(1)-
\frac{\Gamma'}{\Gamma}(k)) F(k,1-k;1;1-x^2)-\\
-x\frac{\partial}{\partial \varepsilon}(|1-x^2|^{\varepsilon}
F(k+\varepsilon,1-k+\varepsilon;1+2\varepsilon;1-x^2))\big\vert_{\varepsilon=0}
\endmultline                                \tag{7.7}
$$

It follows from this equality that two functions
$$
\multline
v=\sqrt{2\tan \xi/2}A_{00}(i(k-1/2),\cos \xi/2;1/2,1/2),\\
\tilde v=\sqrt{2\tanh \xi/2}A_{00}(i(k-1/2),\cosh \xi/2;1/2,1/2)
\endmultline                           \tag{7.8}
$$
are the solutions of the differential equations
$$
v''+((k-1/2)^2 +\frac{1}{4\sin^2\xi})v=0,\tilde v''+(-(k-1/2)^2
+\frac{1}{4\sinh^2 \xi})\tilde v=0.            \tag{7.9}
$$
\proclaim{Lemma 7.2} Let $k$ be sufficiently large; then for
$0\leqslant\xi\leqslant3\pi/4$ we have the asymptotic series
$$
\multline
-\pi \sqrt {2\tan \xi/2}A_{00}(i(k-1/2),\cos \xi/2;1/2,1/2)=\\
=\sqrt \xi\,Y_0((k-1/2)\xi)\sum_{n\geqslant
0}\frac{a_n(\xi)}{(k-1/2)^{2n}}+   \\
(\sqrt \xi\,Y_0((k-1/2)\xi))'\sum_{n\geqslant
1}\frac{b_n(\xi)}{(k-1/2)^{2n}},
\endmultline                               \tag{7.10}
$$
where $a_0\equiv 1$ and all $a_n,\,\, \xi^{-1}b_n$ with $1\leqslant n
\leqslant M$ are bounded smooth functions for any fixed $M$;furthermore,
we have the expansion of the same form for $(2\tan \xi/2)^{-1/2}
A_{00}(i(k-1/2),\sin \xi/2;1/2,1/2), 0\leqslant \xi\leqslant 3\pi/4$.
\endproclaim
\proclaim{Lemma 7.3} Let $k$ be sufficiently large;then for
$\xi\geqslant 0 $ we have the uniform asymptotic formula
$$
\sqrt{2\tan \xi/2} A_{00}(i(k-1/2),\cosh \xi/2;1/2,1/2)=2\sqrt \xi
K_0((r-1/2)\xi)(1+O(\frac{1}{k})).               \tag{7.11}
$$
\endproclaim

The proofs may be omitted for both lemmas since these ones are exactly
the same what we had in \S 5.

\noindent
\subhead
7.4.The coefficients $g_0(k)$ with large $k$'s
\endsubhead
\proclaim{Lemma 7.4} Let $k\geqslant k_0(T) =T^4$; then we have
$$
\int\limits_{(1/2)}\omega_T(\rho)\sum_{k\geqslant k_0}
g_0(k)\sum_{j}\alpha_{j,2k}\Cal H_{j,2k}(\rho+it) \Cal H_{j,2k}(\rho-it)
\Cal H^2_{j,2k}(1/2)\,d\rho \ll T_0            \tag{7.12}
$$
\endproclaim

To get this estimate we use the initial definition (2.21); we integrate
on the line $\Delta= \r{Re}\,w =-2$. The integrand on this line is not
larger than
$$
|k+w|^{-5}\,(|w-\rho|+1)^{-6}\,(|w|+1)^4
$$
(we use the Stirling expansion and the estimate $|\hat \Phi(2w)|\ll
|w|^{2\r{Re}\,w -6}$ ).

Consequently, we have
$$
|g_0(k)|\ll \frac{{T_0}^4}{k^4}                   \tag{7.13}
$$

It follows from the functional equation that for $|\rho| =o(k)$
$$
|\Cal H_{j,2k}(\rho+it) \Cal H_{j,2k}(\rho-it)| \ll k^2.   \tag{7.14}
$$

It gives us the bound
$$
\ll \sum_{k\geqslant k_0} \frac{T_0^5}{k^2}\sum_{j}\alpha_{j,2k}\Cal
H^2_{j,2k}                                       \tag{7.14}
$$
for the sum (7.12). It is known [10] that
$$
\sum_{k\leqslant L}\sum_{j}\alpha_{j,2k}\Cal H^2_{j,2k}(1/2) \ll L\,\log
L \tag{7.16}
$$
and we come to (7.12).
\noindent
\subhead
7.5. Sum with $k\leqslant {T_0}^4$
\endsubhead
\proclaim {Proposition 7.2}
$$
\multline
\frac{1}{4i}\int\limits_{(1/2)}\omega_T(\rho)
Z^{(r)}(\rho,\nu;1/2,1/2\vert g_0)\,d\rho =\\
=\sum_{N,l\geqslant 1}\sum_{k\leqslant T^4_0}
G(k,N,l)\sum_{j}\alpha_{j,2k}\frac{t_{j.2k}(N)\tau_{\nu}(N)}{(Nl^2)^{1/2+iT_0}}
\Cal H^2_{j,2k}(1/2) +O(T_0)
\endmultline                        \tag{7.17}
$$
where $\psi$ is defined by (6.14) and
$$
\multline
G(k,N,l)=\frac{(2k-1) T}{\pi}\int_0^{\infty}A_{00}(i(k-1/2),\sqrt
x;1/2,1/2)\times\\
\times\Biggl(\int_{-\infty}^{\infty}A_{00}(u,\frac{1}{\sqrt x};1/2,\nu)
h(u)\,d\chi(u)-\\
-\sum_{l\in L}(-1)^l c(l)A_{00}(-i(l-1/2),\frac{1}{\sqrt
x};1/2,\nu)\Biggr)\psi(Nl^2x)x^{-1/2-iT_0}\,dx
\endmultline                                      \tag{7.18}
$$
\endproclaim

This representation follows from (7.2) by the same way as in Proposition
6.2; the part of $Z^{(r)}$ with $k\geqslant T_0^4$ gives $O(T_0)$.
\noindent
\subhead
7.5.1. Case $Nl^2\geqslant 2$
\endsubhead
\proclaim{Lemma 7.5} For any $\varepsilon >0$ the part of the sum in
(7.17) with $Nl^2\geqslant 2$ is not larger than
$O(\exp(-t)+\exp(-T^{\varepsilon}))$.
\endproclaim

Really, the function $\psi(x)$ is $O(\exp(-(NT)^{\varepsilon}))$ if
$|x-1|\geqslant (NT)^{\varepsilon} T^{-1}$. So the part of our
integral with $|Nl^2 x -1|\geqslant \frac{(NT)^{\varepsilon}}{T}$ gives
the contribution $O(\exp(-T^{\varepsilon}))$ into (7.17).

Under the condition $Nl^2\geqslant 2$ we have in the remaining part
$x\leqslant 1/2$; but for these values of $x$ we have (see (5.62))
$$
|A_{00}(u,\frac{1}{\sqrt x};1/2,\nu)| \ll \exp(-\pi t/2)
$$
and the same estimate is valid if $u=-i(l-1/2)$ with the fixed integer
$l$.

So for $N \ll T^B$ for some fixed $B$ we have for our sum the estimate
$O(e^{-t})$. But for $N \gg T^B$ we can integrate over the variable $u$
on the line $\r{Im}\,u =-3/2$ when the integrand contains
$I_{2iu}(t\xi) $(here $ \sin \xi/2 =\sqrt x$; in terms with $I_{-2iu}$ we
change $u$ by $-u$). For $\xi\ll(Nl^2)^{-1/2}$ it gives the inequality
$$
|\int_{-\infty}^{\infty}A_{00}(u,\frac{1}{\sin
\xi/2};1/2,\nu)\,h(u)\,d\chi(u)| \ll \frac{t^3}{(Nl^2)^{3/2}}\,e^{-\pi
t}.                                          \tag{7.20}
$$

It means the series in $N, l$ is convergent and we get $O(e^{-t})$
again.
\noindent
\subhead
7.5.2. Case $N=l=1$
\endsubhead
\proclaim{Lemma 7.6} For any fixed $\varepsilon >0$ we have
$$
|G(k,1,1)|\ll
\cases
T\, k^{-1} \log^2 k,\,\,k\geqslant
T^{1/2+\varepsilon},\\
k \,\log^2 T, \,\,k\leqslant T^{1/2+\varepsilon}.
\endcases
\tag{7.21}
$$
\endproclaim

It is sufficient estimate that part of integral where $|x-1|\ll
\frac{\log T}{T}$.

Doing the change of the variable $x\mapsto \cos^2 \xi/2$ for $x\leqslant
1$ and $x\mapsto \cosh^2 \xi/2$ for $x\geqslant1$ and using the
asymptotic formulas (5.67),(5.59),(7.11) and (7.10) we come to the
integral with the main term
$$
\multline
\int_0^{\delta}Y_0((k-1/2)\xi) K_0(t\xi) \psi(\cos^2 \xi/2)(\cos
\xi/2)^{-2iT_0}\xi\,d\xi +\\
+\frac{1}{2}\int_0^{\delta}K_0((k-1/2)\xi) Y_0(t\xi) \psi(\cosh^2 \xi/2)
(\cosh^2 \xi/2)^{-2iT_0}\xi\,d\xi
\endmultline                                   \tag{7.22}
$$
where $\delta =\delta (T) =\frac{\log T}{\sqrt T}$.

In this integral $K_0(t\xi)$ and $Y_0(t\xi)$ may be replaced by the
corresponding power series since $t\xi\ll T^{-1/4}$.

If $k\delta(T)\ll T^{\varepsilon}$ both integrals here are
$O(\frac{\log^2 T}{T})$; it gives the second inequality in (7.21).

Let now $k\delta(T)\gg T^{\varepsilon}$. The second integral in (7.22)
is
$$
O(\int_0^{\infty}K_0((k-1/2)\xi)\xi (|\log \xi| +1)\,d\xi)=O(\frac{\log
k}{k^2}).
$$
In the first integral the small interval $\xi\leqslant k^{\varepsilon
-1}$ contributes $O(\frac{\log^2 k}{k^2})$.

If $k^{\varepsilon -1}\leqslant \xi \leqslant \delta(T)$ we use the
asymptotic expansion for $Y_0((k-1/2)\xi)$. After that we have the
integrals of the form
$$
\frac{1}{\sqrt k}\int\nolimits\exp(\pm i(k-1/2)\xi -2iT_0 \log \cos
\xi/2 +\log \psi(\cos^2 \xi/2))\sqrt \xi\,d\xi.       \tag{7.23}
$$

The point of the stationary phase is near to $(2k-1)T^{-1}_0$; for
$k\geqslant T^{1/2+\varepsilon}$ this point lies outside of the interval
$(0,\delta(T))$. So we can integrate by parts any times and this part of
our integral is $O(k^{-2} \log^k )$ also.
\noindent
\subhead
7.5.3.Sum over regular cusps
\endsubhead
\proclaim{Lemma 7.7} For any fixed $\varepsilon >0$ we have for our
specialization ((3.1), (3.2) and the conditions for $h$ in 3.1)
$$
\int\limits_{(1/2)}\omega_T(\rho) Z^{(r)}(\rho,\nu;1/2,1/2\vert
g_0)\,d\rho \,\,\ll T^{1+\varepsilon}.            \tag{7.24}
$$
\endproclaim

This estimate is the direct consequence of (7.17), (7.21), Lemma 7.5 and
(7.16). The integral on the left side (7.24) is not larger than
$$
\sum_{k\leqslant k_0} k\,\log^2 T\,\,\sum_{j}\alpha_{j,2k} \Cal
H^2_{j,2k}(1/2) +\sum_{k_0\leqslant k \leqslant
T^4}\frac{T}{k}\sum_{j}\alpha_{j,2k}\Cal H^2_{j,2k}(1/2)\ll
T^{1+2\varepsilon}\log^3 T                       \tag{7.25}
$$
(it follows from (7.16)); this inequality concludes the proof.

\head
\S 8. $\Cal R$ AND $\Cal R_h$ -- FUNCTIONS
\endhead

To finish the estimation it rests consider the average of 12 terms with
$\Cal R_h$ and $\Cal R$ on the right side (2.57) for the case of our
specialization.

Here some terms have the poles at $\mu=1/2$ and $s=1/2$, but the full
sum must be regular since all other sums and integrals in this identity
are regular at this point.
\noindent
\subhead
8.1. Terms without pole at $s=1/2$
\endsubhead

Firstly we consider the simplest case when there are no poles at
$s=1/2$.
\proclaim{Proposition 8.1} For $\nu=1/2+it,\rho=1/2+i\tau,\,\,t,\tau
\to\infty, t=o(\tau^{1/4})$ we have
$$
\Cal R_h(1/2,\nu;\rho,1/2)\equiv 2(4\pi)^{-2\nu}\hat\Psi(2\nu)
\frac{\zeta(2\rho)\zeta(2\nu)}{\zeta(2\rho+2\nu)}\zeta^4(\rho+\nu)
=O(t^{-6}\log^6 \tau)
                                             \tag{8.1}
$$
\endproclaim

This equality is the result of the substitution $s=\mu=1/2$ in the
definition (2.18); after that we use the trivial estimate
$\zeta (1+it)=O(\log t)$.
\proclaim{Proposition 8.2} For any fixed $M\geqslant 2$ we have
$$
\Cal R(1/2,\nu;\rho,1/2\vert
h)=-\frac{\zeta(2\nu)\zeta^2(\rho+\nu)\zeta^2(\rho-\nu)}{\zeta(1+2\nu)}\,
h(-i\nu)=O(\frac{1}{t^M}).    \tag{8.2}
$$
\endproclaim

It follows from the definition and our assumptions about $h$.

\proclaim{Proposition 8.3} For any fixed $M\geqslant 2$ we have
$$
\lim_{\mu \to 1/2}\Biggl(\Cal R(\rho,\mu;1/2,\nu\vert h)+\Cal
R(\rho,1-\mu;1/2,\nu\vert h)\Biggr) =O(\tau^{-M}) \tag{8.3}
$$
\endproclaim

Writing $\zeta(2\mu)=\frac{1}{2\mu-1}+ \gamma +...$ (here $\gamma$ is
the Euler constant) we see that the left side on (8.3) equals to
$$
\multline
2\gamma\,h(i(\rho-1))\frac{\zeta(2\rho-1)}{\zeta(3-2\rho)}\zeta(\rho-\nu)
\zeta(1-\rho+\nu)\zeta(\rho+\nu-1)\zeta(2-\rho-\nu)+\\
+2\zeta(2\rho-1)\frac{\partial}{\partial \mu}\frac{h(i(\rho-\mu-1/2))}
{\zeta(2-2\rho+2\mu)}\Cal Z(1/2;\nu,\rho-\mu)\vert_{\mu=1/2}
\endmultline
$$

Now the assertion follows since $h(r)$ decreases more rapidly than any
fixed degree of $r$.

\proclaim{Proposition 8.4} Under our assumptions for $h$ we have for any
fixed $M\geqslant 4$
$$
\Cal R(1/2,\nu;\rho,1/2)=-\frac{\zeta(2\nu)\zeta^2(\rho+\nu)\zeta^2(\rho-\nu)}
{\zeta(2\nu+1)}\,h(-i\nu)=O(t^{-M}).           \tag{8.4}
$$
\endproclaim

It is the consequence of the definition and the fast decreasing of $h$.

\proclaim{Proposition 8.5} Let $\tilde h_0(r),\tilde h_1(r)$ are the
values of two integrals (2.19) and (2.20) for the case
$\r{Im}\,r<-\Delta$; then for $-1-\Delta<r<-\Delta$ we have the
analitical continuation
$$
\multline
h_0(r)=\tilde h_0(r)+
2^{2ir}\Gamma(2ir)\gamma(\rho-ir,\nu)\gamma(s-ir,\mu)\cosh \pi r
\times\\
\times\Biggl(\cos \pi(\rho-ir) \cos \pi(s-ir) +\sin \pi \nu \sin \pi
\mu)\Biggr)\,\hat\Psi(2ir-2s-2\rho+2),
\endmultline                                     \tag{8.5}
$$
$$
\multline
h_1(r)=\tilde h_1(r) +
2^{2ir}\Gamma(2ir)\gamma(\rho-ir,\nu)\gamma(s-ir,\mu) \cosh \pi r
\times\\
\times\Biggl(\cos \pi(s-ir) \sin \pi \nu +\cos \pi(\rho -ir) \sin \pi
\mu\Biggr)\,\hat\Psi(2ir -2s -2\rho +2).
\endmultline                                        \tag{8.6}
$$
\endproclaim

Really, these integrals (for $0<\Delta<1/2$) define the even regular
functions in the strip $\r{Im}\,r<\Delta$. Writing
$$
\gamma(w,1/2+ir)=\frac{2^{2w-1}}{\pi}\frac{\Gamma(w+ir+1)\Gamma(w-ir+1)}
{(w+ir)(w-ir)}                      \tag{8.7}
$$
we come to the Cauchy type integral. We have two simple poles at $w=\pm
ir$. When $\r{Im}\,r$ is near to $-\Delta$ we deform the path of the
integration so that the point $ir$ lies to right on the path. When
$\r{Im}\,r<-\Delta$ we can take the initial line as the path of the
integration; it gives us the relation
$$
h_0(r)=\tilde h_0(r)+2\pi i \r{Res}_{w=ir} \Biggl(-i(\text {the
integrand})\Biggr)
                                                         \tag{8.8}
$$
and we get (8.5) and (8.6).
\proclaim{Proposition 8.6} Under assumptions of subsection 3.1 we have
$$
\Cal R(\rho,\nu;1/2,1/2\vert h_0+h_1) =O(t^{-5}\log^2 \tau),\tag{8.9}
$$
$$
\Cal R(\rho,1-\nu;1/2,1/2\vert h_0+h_1)=O(t^{-5}\log^2 \tau).\tag{8.10}
$$
\endproclaim

First of all (we use (2.49) with $s=\mu=1/2$)
$$
\multline
\Cal R(\rho,\nu;1/2,1/2\vert
h_j)=2\frac{\zeta(2\rho-1)\zeta(2\nu)}{\zeta(2-2\rho+2\nu)}\zeta^2(1-\rho+\nu)
\zeta^2(\rho-\nu)\,h_j(i(\rho-\nu-1/2))\\
=O(\tau^{3/2} \log^2 \tau)\,|h_j(i(\rho-\nu-1/2))|.
\endmultline                          \tag{8.11}
$$

Now for $h_j$, $j=0$ or $j=1$, we use (8.5)--(8.6). In these equalities
we have for difference $h_j-\tilde h_j$ the bound $O(\tau^{-3/2})
|\hat \Psi(2\nu-4\rho+2)|=O(\tau^{-13/2})$ since $|\hat \Psi(2w)|\ll
|w|^{2\r{Re} w -6}$.

So it rests estimate the integrals $\tilde h_j$ at the point
$r=i(\rho-\nu-1/2)$.

We have
$$
\multline
\tilde h_j(i(\rho-\nu-1/2))=\frac{i}{\pi^2}\int\limits_{(\Delta)}
\frac {\Gamma(w+\rho-\nu-1/2)\Gamma(\rho-w+\nu-1/2)}{\cos \pi(\rho-w-\nu)}
\Gamma^2(1/2-w)\times\\
\times S_j (w,\rho,\nu)\hat \Phi(2w-2\rho+1)\,dw,
\endmultline                                               \tag{8.12}
$$
where it is assumed $0<\Delta<1/2$ and
$$
S_0=-2^{2\rho-2w-4}\Biggl(\sin \pi(\rho+w) +\sin \pi(3w-\rho) +2\sin
\pi(\nu+w) +2\sin \pi(\nu-w)\Biggr),        \tag{8.13}
$$
$$
\multline
S_1=2^{2\rho-2w-4}\Biggl(\sin \pi(w+\rho-2\nu) + \sin \pi(w+2\nu-\rho) +
\sin \pi(2\rho-\nu-w)-\\
-\sin \pi(w+\rho) -\sin \pi(w-\rho) +\sin \pi(w-\nu) \Biggr).
\endmultline                                 \tag{8.14}
$$

Firstly we consider the integral $h_0$.

Let us write $w=\Delta+i\eta$; then for $\rho=1/2+i\tau$, $\nu=1/2+it$
we have the following estimate for the integrand:
$$
\multline
\ll \exp(-\frac{\pi}{2}\Cal M)\times\\
\times(|\eta+\tau-t|+1)^{\Delta-1}(|\tau-\eta+t|+1)^{-\Delta}
(|\eta|+1)^{-2\Delta}(|\eta-\tau|+1)^{2\Delta-6}
\endmultline                                     \tag{8.15}
$$
where
$$
\Cal
M=|\eta+\tau-t|+|\tau-\eta+t|+2|\eta|-2|\tau-\eta-t|
-2\max(|\tau-3\eta|,|\tau+\eta|))
$$
(it follows from the Stirling expansion and from our construction of
$\Phi$).

If $\eta \leqslant -2\tau$ or $\eta\geqslant 2\tau$ this estimate gives
not larger than $O(|\eta|^{-7})$ and this part of our integral is
$O(\tau^{-6})$.

Now for $-\tau+t\leqslant \eta \leqslant 0$ we have the exponential
multiplier $\exp(-\pi(\tau-t+\eta))$, for $0\leqslant \eta \leqslant
\tau-t$ there is $\exp(-\pi(\tau-t-\eta))$ in the integrand, for
$\tau-t\leqslant\eta\leqslant \tau$ this multiplier is
$O(\exp(-\pi(\eta-\tau+t)))$ and for $\tau\leqslant\eta\leqslant \tau+t$
we have the multiplier $\exp(-\pi(\tau+t-\eta))$. This exponential
multiplier equals to 1 for $\eta\leqslant -\tau+t $ and for
$\eta\geqslant \tau+t$; it is exponentially small for $\eta\in
(-\tau+t,\tau+t)$. So the parts of our integral with $-\tau+t+2\log
\tau\leqslant \eta\leqslant \tau-t-2\log \tau$ or $\tau-t+2\log
\tau\leqslant \eta\leqslant\tau+t-2\log \tau$ give the small
contribution. The part of this integral with $\eta\approx -\tau$ is
$O(\tau^{-5})$.

Finally, for $\eta\approx \tau$ we move the path of
integration on the line $\r{Re}\,w=2$. The residues at
$w=\rho\pm(\nu-1/2)$ give $O(\tau^{-3/2}\, t^{-11/2})$ and the result of
integration over the line $\r{Re}\,w=2$ is $O(\tau^{-3}\,t^{-2})$.

We have the similar situation for the case $j=1$. The integrand contains
the multiplier $\exp(-\pi |\eta|)$ for $\eta\leqslant
\frac{1}{2}(\tau-t)$, the multiplier $\exp(-\pi(\tau-t-\eta))$ for
$\frac{1}{2}(\tau-t)\leqslant \eta \leqslant \tau+t$ and $\exp(-2\pi
(\eta-\tau))$ for $\eta\geqslant \tau+t$.

It means this integral is determined by two intervals: $|\eta|\leqslant
2 \log \tau$ and $|\eta -\tau+t|\leqslant 2\log \tau$. In the case
$w\approx \rho-\nu$ we integrate over the line $\r{Re}\,w=2$ and the
pole at $w=\rho-\nu+1/2$ contributes $O(\tau^{-3/2} t^{-11/2})$ again.

Together with (8.11) the estimates of $|\tilde h_j|$ give (8.9) and
(8.10).
\noindent
\subhead
8.2. Terms with pole at $s=1/2$
\endsubhead
\proclaim{Proposition 8.7} Under our assumptions for $h$ we have
$$
\multline
\lim_{s,\nu \to 1/2}\Biggl(\Cal R_h(\rho,\mu;s,\nu))+\Cal
R_h(\rho,1-\mu;s,\nu)+\\
+\Cal R(s,\mu;\rho,\nu\vert h_0+h_1)+\Cal R(s,1-\mu;\rho,\nu\vert
h_0+h_1\Biggr)=O(t^{-5} \log^2 \tau)
\endmultline                                  \tag{8.16}
$$
\endproclaim

Here two first terms have the pole at $s=1/2$, since
$$
\Cal R_h(\rho,\mu;s,\nu)=2\zeta(2s)(4\pi)^{2\rho-2\mu-1}\hat
\Phi(2\mu+1-2\rho)\frac{\zeta(2\mu)}{\zeta(2s+2\mu)}\Cal
Z(s+\mu;\rho,\nu)                    \tag{8.17}
$$
(sum $\Cal R_h(\rho,\mu;s,\nu)+\Cal R_h(\rho,1-\mu;s,\nu)$ has no pole
at $\mu=1/2$).

To see that sum in (8.16) has no pole at $s=1/2$ we transform the
expression for $\Cal R$. By the definition,
$$
\Cal R(s,\mu;\rho,\nu\vert
h_0+h_1)=2\frac{\zeta(2s-1)\zeta(2\mu)}{\zeta(2-2s+2\mu)}\Cal
Z(\rho;\nu,s-\mu)(h_0+h_1)(i(s-\mu-1/2)). \tag{8.18}
$$
First of all we have (we use (8.6),(8.7) and the definition (1.31))
$$
\multline
(h_0+h_1)(i(s-\mu-1/2))=(\tilde h_0+\tilde h_1)(i(s-\mu-1/2))+\\
+2^{2\rho+4s-2\mu-4}\frac{\Gamma(2s-1)}{\pi \cos
(s-\mu)}\Gamma(\rho+\nu+s-1)\Gamma(\rho-\nu+s-\mu)\times\\
\times \bigl(\sin \pi(2s-\mu) +\sin \pi \mu)(\sin \pi(\rho+s-\mu) +\sin \pi
\nu \bigr) \hat \Phi(2\mu+3-4s-2\rho)
\endmultline                               \tag{8.19}
$$
(the functional equation for gamma-function have been used here).

The last term has the pole at $s=1/2$ also; it compensates the pole from
$\Cal R_h$.

The following identity follows from the Riemann functional equation and
the doubling formula for the gamma-function.
\proclaim{Proposition 8.8}
$$
\multline
(\sin \pi(\rho+s-\mu)+\sin \pi
\nu)\Gamma(\rho+\nu+s-\mu-1)\Gamma(\rho-\nu+s-\mu)\Cal
Z(\rho;\nu,s-\mu)=\\
=\frac{1}{2}(2\pi)^{2\rho+2s-2\mu-1}\Cal Z(1-s+\mu;\rho,\nu)
\endmultline                                               \tag{8.20}
$$
\endproclaim

Really, changing zeta-functions of arguments $\rho-\nu+s-\mu$ and
$\rho+\nu+s-\mu-1$ by the functions at $1-\rho+\nu+\mu-s$ and
$2-\rho-\nu+\mu-s$ correspondingly we get
$$
\multline
\Cal
Z(\rho;\nu,s-\mu)=\zeta(\rho+\nu-s+\mu)\zeta(\rho-\nu+s-\mu)
\zeta(\rho+\nu+s-\mu-1)\zeta(\rho-\nu-s+\mu+1)\\
=\Cal
Z(1-s+\mu;\rho,\nu)\times\\
\times\pi^{2\rho+2s-2\mu-2}
\frac{\Gamma((1-\rho+\nu+\mu-s)/2)\Gamma((2-\rho-\nu+\mu-s)/2)}
{\Gamma((\rho-\nu+s-\mu)/2)\Gamma((\rho+\nu+s-\mu-1)/2)}
\endmultline                                    \tag{8.21}
$$
Now the doubling formula gives
$$
\multline
\Gamma(\rho+\nu+s-\mu-1)\Gamma(\rho-\nu+s-\mu)=\pi^{-1}2^{2\rho+2s-2\mu-3}
\Gamma((\rho+\nu+s-\mu-1)/2) \\
\times\Gamma((\rho+\nu+s-\mu)/2)
\Gamma((\rho-\nu+s-\mu)/2)\Gamma((\rho-\nu+s-\mu+1)/2)
\endmultline                                        \tag{8.22}
$$
and we have (using the functional equation for gamma-function again)
$$
\multline
\Gamma(\rho+\nu+s-\mu-1)\Gamma(\rho-\nu+s-\mu)\Cal Z(\rho;\nu,s-\mu)=\\
=\frac{(2\pi)^{2\rho+2s-2\mu-1}}{\sin \pi(\rho+\nu+s-\mu)/2 \,\,\sin
\pi(\rho-\nu+s-\mu+1)/2 } \Cal Z(1-s+\mu;\rho,\nu);
\endmultline                                         \tag{8.23}
$$
this equality coincides with (8.20).

We see now that the coefficient in front of $\zeta(2s)$ equals to
$$
2\frac{\zeta(2\mu)}{\zeta(1+2\mu)}\Biggl((4\pi)^{2\rho-2\mu-1}
\hat\Phi(2\mu+1-2\rho)\Cal Z(\mu+1/2;\rho,\nu)+
(s-1/2)r_1(\mu,\nu)...\Biggr)                       \tag{8.24}
$$
and the coefficient in front of $\Gamma(2s-1)$ equals to (note that
$\zeta(0)=-1/2$)
$$
\multline
2\frac{\zeta(2\mu)}{\zeta(1+2\mu)}\Biggl(-(4\pi)^{2\rho-2\mu-1}
\hat\Phi(2\mu+1-2\rho)\Cal Z(\mu+1/2;\rho,\nu)\sin \pi \mu
+\\
+(s-1/2)r_2(\mu,\nu)...\Biggr),
\endmultline                     \tag{8.25}
$$
where $r_1(\mu,\nu),\,r_2(\mu,\nu)$ are regular at $\mu=1/2$.

It means there is the limit $s \to 1/2$ (8.16) and it has the form
$$
  \zeta(2\mu)\Cal F(\mu,\rho) +\zeta(2-2\mu)\Cal F(1-\mu,\rho);
$$
here $\Cal F(\mu,\rho)$ is regular at $\mu=1/2$ and contains the product of
four zeta-functions on the unit line (and the derivative of this
product). 

Now the existence of the limit $\mu\to 1/2$ is obvious; there is no need
to calculate the explicit form of this limit. One can see that the
result contains the derivatives of zeta-functions of the order
$\geqslant 2$ and the derivatives of $\hat\Phi$ at the point $2-2\rho$
of the order 0,1,2. Since these last derivatives have the same order as
$\hat \Phi$ (each differentiation gives the additional multiplier $O(\log
\tau)$ only) we get the estimate $O(\tau ^{-6}\log^6 \tau)$ for this
limit.

It rests estimate two terms with $(\tilde h_0+\tilde h_1)(i/2)$. We have
$\tilde h_1(i/2)=0 $ and (for $s=\mu=1/2$)
$$
\multline
\tilde h_0(i/2)=\frac{2^{2\rho-2w-1}}{i
\pi}\int\limits_{(\Delta)}\frac{\cos \pi(\rho-w) \sin \pi w +\sin \pi
\nu}{(2w-1) \cos \pi w}\times\\
\times\Gamma(\rho-w+\nu-1/2)\Gamma(\rho-w-\nu+1/2)\hat\Phi(2w-2\rho+1)\,dw.
\endmultline                                                    \tag{8.26}
$$

The part of this integral with $|w|\leqslant (1-\delta)\,\tau$ with some
fixed (small) $\delta >0$ is $O(\tau^{-6} \log \tau)$. For
$\eta\geqslant (1-\delta)\,\tau$ (here $\eta=\r{Re}\,w$) we have
$$
\multline
|\tau-\eta+t|+|\tau-\eta-t|-2|\tau-\eta|=0,\,\,(1-\delta)\tau\leqslant
\eta\leqslant \tau-t,\\
=2(\eta-(\tau-t)),\,\,\tau-t\leqslant \eta\leqslant \tau,\\
=2(\tau+t-\eta),\,\,\tau\leqslant \eta\leqslant \tau+t,\\
=0,\,\,\,\,\,\,\,\,\,\,\,\,\eta\geqslant \tau+t.
\endmultline
$$

It means thar our integral is determined by the interval
$|\tau-t-\eta|\leqslant 2\log \tau$ ang by halfaxis $\eta\geqslant
\tau+t-2\log \tau$. For this $\eta$ we can integrate over the line
$\Delta =0$; since for both cases $|w-\rho|\geqslant t$ we have
$$
|\tilde h_0(i/2)|\ll \tau^{-1}\,t^{-5}.               \tag{8.27}
$$
Taking into account that
$$
|\Cal Z(\rho;\nu,0)|\ll \tau\,\,\log^2 \tau,           \tag{8.28}
$$
we come to the final estimate (8.16).

\head
\S9. THE LINDEL\"OF CONJETCTURE
\endhead

It rests collect together all previous estimates to finish the proof
of (0.2) and (0.3).

As it had been shown in subsection 3.1 the result of the averaging of
the left side of our main functional equation (2.57) for the special
case $\rho=1/2+i\tau,\nu=1/2+it,s=\mu=1/2$ for any positive (small)
$\delta$  and for any fixed $M\geqslant 1$ is not smaller than
$$
\multline
T\Biggl(\sum_{j\geqslant 1}\alpha_j |\Cal H_j(\nu)|^2 \,h(\varkappa_j) +\\
\frac{1}{\pi}\int_{-\infty}^{\infty}\frac{|\zeta(\nu+ir)\zeta(\nu-ir)|^{2}}
{|\zeta(1+2ir)|^2} h(r)\,dr\Biggr)+o(1)\gg\\
\gg T/4 \sum_{1\leqslant j\leqslant M}\alpha_j |\Cal H(\nu)|^2
h(\varkappa_j) +\frac{T}{4\pi}\int_{-\delta}^{\delta}
\frac{|\zeta(\nu+ir)\zeta(\nu-ir)|^2}{|\zeta(1+2ir)|^2}\,h(r)\,dr
\endmultline                                               \tag{9.1}
$$

If $\delta\ll(\log t)^{-2}$ the last integral is not smaller than
$$
( h(0)/3\pi)\,T\, \delta^3 |\zeta(1/2+it)|^4.             \tag{9.2}
$$

On the other side this average is not larger than
$$
\log^2 T\sum_{\varkappa_j\leqslant \sqrt T_0}\alpha_j \Cal H^2_j(1/2) +
T\sum_{\sqrt T_0\leqslant \varkappa_j\leqslant T^4_0}\alpha_j
\frac{\log^2 \varkappa}{\varkappa^2}\Cal H^2_j(1/2)+T^{1+\varepsilon}+
$$
(Lemmas 6.1,\,\,6.2,\,\,6.6)
$$
+\log^2 T\sum_{k\leqslant k_0}k\sum_{j}\alpha_{j,2k}\Cal
H^2_{j,2k}(1/2)+T\sum_{k_0\leqslant k \leqslant k_1}|\frac{\log^2
k}{k}\sum_{j}\alpha_{j,2k}\Cal H^2_{j,2k}(1/2)+
$$
(Lemma 7.6; here $k_0=T_0^{1+\varepsilon}, k_1=T^4$)
$$
+ O(T\frac{\log^5}{t^5})+
$$
(Propositions 8.6,\,\,8.7)

+\{terms of the smaller order\}.

All listed sums are $O(T^{1+\varepsilon})$ for any fixed
$\varepsilon>0$. Consequently (we take $\delta =(\log t)^{-2}$ in (9.2)),
$$
|\zeta^4(1/2+it)|\ll T^{\varepsilon}           \tag{9.3}
$$
and for every fixed $j\geqslant 1$
$$
|\Cal H_j(1/2+it)|^2 \ll T^{\varepsilon}.         \tag{9.4}
$$

The unique condition for $t$ is $t=o(T^{1/4})$; so, supposing
$t=T^{1/8}$ we come to (0.2) and (0.3).

\head
\S 10. SOME CONCLUDING REMARKS
\endhead

I had been writing this work many years. The plans of the proof had been changed
many times. Some parts had retained, the others were written anew.

Now the potential reader has the definitive text. May be, it would be
useful to indicate the crucial points of the proposed proof.

Of course, the base of the whole construction is the the fore-trace
formulas (1.13) and (1.25) which are well known today.

The first crucial point is the regularization (2.11) where the
coefficients are defined by (2.10). This simple trick took off many
problems with the convergence and it is hardly understand why this
method was not used early.

The idea of the special averaging of the main functional equation (2.61)
(when the unknown quantity -- namely, $ |\zeta(1/2+it)|^4$ -- is the
coefficient in front of large parameter) is blowed by the proof of the
Dirihlet formula for number of classes.

To realize this idea the new representations (4.17) are necessary
(together with the initial expressions (2.19) and (2.20)).

The whole technical \S 5 must be outside of this work. This section may
be considered as the supplement to the Erd\'elyi-Bateman handbook; one
can compare our approach with text on the pages 602-643 in [14].

The real crucial step is the replacement of the summation over the
discret spectrum (with very large $\varkappa_j$) by the integration
over this variable. The replacement of the sum of the quatities (6.52)
over the full spectrum by the sum (6.72) allows us use my convolution
formulas and, as the result, to finish estimates.

\newpage
\head
REFERENCES
\endhead
\bigskip
1. A. Selberg, S. Chowla, {\it On Epstein's zeta-function}, J. reine und
angew. Math.,227,1967,p. 86-110
\bigskip
2. Кузнецов Н.В.,{\it Гипотеза Петерсона для параболических форм веса
нуль и гипотеза Линника,I, Суммы сумм Клоостермана}, Мат. сборник, III
(153), N 3, 1980, с. 334-383.
\bigskip
3. R.W. Bruggeman, {\it Fourier coefficients of cusp forms}, Invent.
Math., 445, 1978, p. 1-18.
\bigskip
4. Кузнецов Н.В.,{\it Гипотеза Птерсона для форм веса нуль и гипотеза
Линника}, препринт, Хабаровск, 1977.
\bigskip
5. M.N. Huxley, {\it Introduction to Kloostermania}, in:Elementary and
analytic theory of numbers, Banach center publications, vol. 17,
Warszawa, 1985, p. 217-306.
\bigskip
6. N.V. Kuznetsov,{\it Sums of the Kloosterman sums and the eighth power
moment of rhe Riemann zeta-function}, in: Number theory and related
topics, Tata Inst. of fund. Research, Bombay and Oxford Univ. press,
1989, p. 57-117.
\bigskip
7. H. Bateman, A. Erd\'elyi, {\it Higher transcendental functions}, vol.
1, N.-Y.--Toronto--London, 1953.
\bigskip
8. Кузнецов Н.В.{\it Свертка коэффициентов Фурье рядов
Эйзенштейна--Мааса}, в: Записки научных семинаров ЛОМИ, т.129,
Автоморфные функции и теория чисел, "Наука", 1983, с. 43-83.
\bigskip
9.Кузнецов Н.В.,{\it О собственных функциях одного интегрального
уравнения}, Записки научных семинаров ЛОМИ, 17:3, 1970, с. 66-149.
\bigskip
10. Кузнецов Н.В.,{\it О среднем значении рядов Гекке}, препринт
06-1994, ИПМ ДВО РАН, 35 с.
\bigskip
11. Deshoillers J.-M. and Iwaniec H.,{\it Kloosterman sums and Fourier
coefficients of cusp forms},Invent. Math.,70(1982),p. 219-288.
\bigskip
12. Титчмарш Е.К., {\it Теория дзета-функции Римана}, М., 1953.
\bigskip
13. Виноградов А.И., Тахтаджян Л.А., {\it Дзета-функция аддитивной
проблемы делителей и спектральное разложение автоморфного Лапласиана},
Записки научных семинаров ЛОМИ, 1984, т. 134, с. 84-116.
\bigskip
14. D.A.Hejhal, {\it The Selberg trace formula for $PSL(2,\Bbb R)$},
Springer-Verlag, vol.1001, 1983.

\enddocument